%% file: paper2-via-Morse-theory.tex
\documentclass[a4paper,leqno]{amsart}
\input{pre.tex}

\title[Gorenstein distance algebras of lattices]{Distance algebras of finite distributive lattices are Gorenstein}
\author{Aaron Chan and Bohan Xing}
\date{\today}

\newcommand{\Addresses}{{
  \bigskip
  \footnotesize

  A. Chan, \textsc{Department of Mathematics, Nagoya University, Chikusa-ku, Nagoya 464-8602, Japan}\par\nopagebreak
  \textit{E-mail address}: \texttt{aaron.kychan@gmail.com}

  \medskip

  B. Xing, \textsc{School of Mathematical Sciences, Laboratory of Mathematics and Complex Systems, Beijing Normal University, Beijing 100875, P.R. China}\par\nopagebreak
  \textit{E-mail address}: \texttt{bhxing@mail.bnu.edu.cn}
}}

\begin{document}
\begin{abstract}
Distance algebras were introduced by Asao and Ivanov to realize magnitude homology and cohomology of quasi-metric spaces as Tor- and Ext-groups respectively. 
Let $\Sigma_L$ be the distance algebra of the Hasse graph of a finite distributive lattice $L$, where the metric is given by the lengths of shortest paths in the graph.  
We show that $\Sigma_L$ is Gorenstein over any Gorenstein commutative noetherian ring $R$, and give bound of its self-injective dimension. 
The proof constructs, for each indecomposable (relative) injective module, a projective resolution that is in general non-minimal.  When $R=\Bbbk$ is a field, we use algebraic Morse theory to reduce it to a minimal resolution, so that the Betti numbers can be obtained as the dimensions of homologies of certain order complexes.
As a consequence, we show that the self-injective dimension of $\Sigma_L$ is the maximal Castelnuovo--Mumford regularity of certain Stanley--Reisner rings.  We also show that the first two terms of these minimal resolutions are multiplicity-free, whereas arbitrarily large multiplicities occur in every higher degree.  In the appendix, we include a more abstract formulation of the non-minimal projective resolution, namely, a simplicial replacement of a diagram of projective modules.
\end{abstract}
\maketitle
\tableofcontents

\section{Introduction}
\label{sec:introduction}


Magnitude homology was introduced by Hepworth and Willerton as a categorification of the magnitude of a graph~\cite{HW17}. Magnitude was introduced by Leinster, who extended the Euler characteristic of finite categories~\cite{Lei08} to enriched categories, including metric spaces~\cite{LW13,LM17}, graphs~\cite{Lei19}, and finite-dimensional algebras~\cite{CKL16}; it had appeared independently as a measure of biological diversity~\cite{SP94}, and has since been applied, for example, to machine learning~\cite{LASR24}. 
Asao and Ivanov realized magnitude (co)homology as bigraded Tor- and Ext-groups over a certain $\mathbb{R}$-graded algebra~\cite{AI24}, which we call the \emph{distance algebra}.\footnote{This term is not widely used in the literature yet, but is Asao's preferred terminology.} 
This opens up the possibility of using representation-theoretic techniques to study magnitude (co)homology; for instance, Asao and Wakatsuki constructed minimal projective resolutions of the simple modules over the distance algebras of geodetic metric spaces~\cite{AW24} -- we remark that these are precisely the monomial distance algebras.

As representation theorists, we find it natural to ask for a characterization of quasi-metric spaces whose distance algebras exhibit good representation-theoretic and homological properties. 
In our companion paper~\cite{CX26}, we initiated such an investigation for the distance algebra $\Sigma_G$ of a finite connected graph $G$, where the metric is given by lengths of geodesic paths. 
In the present paper, we focus on the following problem.

\begin{que}\label{que:main}
Determine when $\Sigma_G$ is Iwanaga--Gorenstein, and compute its self-injective dimension from $G$.
\end{que}

Since the global dimension of $\Sigma_G$ is almost always infinite, the Iwanaga--Gorenstein property is the natural homological finiteness condition to ask for. The homological study of Gorenstein rings has its foundations in Bass's work on finite self-injective dimension~\cite{Bas63} and Auslander--Bridger's theory of Gorenstein dimension~\cite{AB69}.  In the noncommutative setting, Iwanaga studied two-sided noetherian rings of finite self-injective dimension on both sides~\cite{Iwa79}; these are now called Iwanaga--Gorenstein rings. 
For finite-dimensional algebras, this class contains both self-injective algebras and algebras of finite global dimension.  Our interest in this class of algebras is due to the following rich homological structure: (1) the existence of Serre duality on the perfect derived category~\cite{Hap91} and (2) the fact that the singularity category is equivalent to the stable category of maximal Cohen--Macaulay modules~\cite{Buc21,Hap91,Orl04}.
The setting of our first main result will follow the terminogolies and framework in~\cite{IK22,GIK23}.


During our investigation in~\cite{CX26}, it appears to us that utilizing combinatorics of poset can be fruitful for studying the distance algebras of the Hasse graphs of finite posets.
While it seems difficult to answer Question~\ref{que:main} completely, it turns out that restricting to the case of finite distributive lattices gives us a partial answer.
By Birkhoff's fundamental theorem~\cite{Bir37}, every finite distributive lattice $L$ is isomorphic to the lattice $\mathcal J(P)$ of order ideals of a finite poset $P$.  In this case, we have shown in~\cite{CX26} that $\Sigma_L$ is Auslander--Gorenstein if and only if every connected component of $P$ is a chain. 
Our first main result shows that the Iwanaga--Gorenstein property holds without any restriction on $P$.

\begin{main}\label{main:gorenstein}
\textnormal{(see Corollary~\ref{cor:distributive-ig} and Proposition~\ref{prop:n-2-inequalities})}
Let $R$ be a commutative noetherian ring, let $P$ be a non-empty finite poset, and let $\Sigma:=\Sigma_{\mathcal J(P)}$ be the distance $R$-algebra of the Hasse graph of $\mathcal J(P)$.
\begin{enumerate}
\item If $R$ is Gorenstein, then $\Sigma$ is a Gorenstein $R$-algebra.
\item If $R$ is Iwanaga--Gorenstein, then $\Sigma$ is Iwanaga--Gorenstein and
\[
 \idim\Sigma_\Sigma=\idim{}_\Sigma\Sigma
 \leq\idim_RR+|P|-1.
\]
If, moreover, $|P|\ge 3$, then the inequality is strict.
\end{enumerate}
\end{main}
Note that for given a finite distributive lattice $L$, the poset $P$ consists of the join-irreducible elements of $L$. Thus the bounds in the theorem can also be stated in terms of the number of join-irreducible elements of $L$.

To prove Theorem~\ref{main:gorenstein}, we fix $X\in\mathcal J(P)$ and consider the (finite) poset
\[
 \mathcal D_X:=\{\underline Z:=X\SDiff Z\mid Z\in\mathcal J(P)\} \text{ ordered by inclusion,}
\]
where $\SDiff$ denotes symmetric difference. 
We use the order complex of $\D_X$ to construct a projective resolution of the relative injective $\Sigma$-module $I_X:=\Hom_R(\Sigma e_X, R)$ associated with $X$. 
Theorem \ref{main:gorenstein} then follows from combining this with duality of $\Sigma$.

In general, this resolution is far from minimal. 
Over a field $\Bbbk$, we can utilize algebraic Morse theory~\cite{CLZ26} to reduce it to a minimal one.
Moreover, it turns out that the projective dimension of $I_X$ can be described explicitly using the simplicial complex on $P$ generated by $\mathcal D_X$, that is,
\[
 \mathcal S_X:=\{U\subseteq P\mid U\subseteq\underline Y\text{ for some }\underline Y\in\mathcal D_X\}.
\]

\begin{main}\label{main:over-field}
\textnormal{(see Theorem~\ref{thm:minimal-resolution} and Theorem~\ref{thm:pdim=reg})}
Let $\Bbbk$ be a field, $P$ be a finite poset, and $\Sigma:=\Sigma_{\mathcal J(P)}$ be the distance $\Bbbk$-algebra of the Hasse graph of $\mathcal J(P)$. 
For $Z\in\mathcal J(P)$, let $P_Z=e_Z\Sigma$ and $I_Z=\Hom_\Bbbk(\Sigma e_Z,\Bbbk)$ be the indecomposable projective and injective $\Sigma$-modules associated with $Z$. 
The following hold for any fixed $X\in\mathcal J(P)$.
\begin{enumerate}
\item The minimal projective resolution $\mathcal P_{\min,\bullet}$ of $I_X$ has terms
\[
 \mathcal P_{\min,k}\cong\bigoplus_{Z\in\mathcal J(P)}P_Z^{\oplus\widetilde\beta_{Z,k}},
 \qquad
 \widetilde\beta_{Z,k}:=\dim_\Bbbk\widetilde H_{k-1}\bigl(\Delta(\mathcal D_X^{>\underline Z})\bigr)
 \qquad(k\geq0),
\]
where $\Delta(-)$ denotes the order complex, $\widetilde{H}$ denotes reduced homology, and $\D_X^{> \underline{Z}}$ is the upset of $\D_X$ consisting of elements strictly larger than $\underline{Z}$.

\item The projective dimension of $I_X$ is precisely the Castelnuovo--Mumford regularity of the Stanley--Reisner ring $\Bbbk[\mathcal S_X]$. Consequently,
\[
 \idim\Sigma_\Sigma=\idim{}_\Sigma\Sigma=\max_{X\in\mathcal J(P)}\{ \text{Castelnuovo--Mumford regularity of }\Bbbk[\mathcal S_X] \}.
\]
\end{enumerate}
\end{main}

Theorem~\ref{main:over-field}\textup{(1)} turns questions about multiplicities into questions about topology. In Section~\ref{sec:low-degree-multiplicity}, we show that $\mathcal P_{\min,0}$ and $\mathcal P_{\min,1}$ are always multiplicity-free, whereas arbitrarily large multiplicities occur in every degree $k\geq2$, already when $P$ is of height two.

\medskip
\noindent
\textbf{Organization of the paper.} Section~\ref{sec:preliminaries} recalls Gorenstein algebras, distance algebras, and the combinatorial and topological tools that we need. Section~\ref{sec:candidate-resolution} constructs the candidate projective resolutions and  proves their exactness, which yields Theorem~\ref{main:gorenstein}. Section~\ref{sec:minimal-morse} obtains the minimal resolutions by algebraic Morse theory and proves Theorem~\ref{main:over-field}. Section~\ref{sec:low-degree-multiplicity} studies the multiplicities. Appendix~\ref{app:thin-modules} explains the projective resolution as a simplicial replacement of a diagram of projective modules over arbitrary finite-dimensional algebras. 
Appendix~\ref{app:lattice-proofs} contains the proofs of two preliminary results from Section~\ref{sec:preliminaries}.

\medskip
\noindent
{\bf Conventions and notation.} Throughout, $R$ denotes a commutative noetherian ring with identity, and $\Bbbk$ denotes a field. All algebras are associative and unital. By an $R$-algebra we mean one that is finitely generated and projective as an $R$-module. Unless otherwise stated, modules are finitely generated right modules that are free over $R$. Unadorned tensors $\otimes$ between modules are tensors $\otimes_R$ over the ground ring $R$. All graphs, posets, lattices, and simplicial complexes are finite. The symbols $\idim$ and $\pdim$ denote injective and projective dimension, respectively. We write $D:=\Hom_R(-,R)$. For a collection $\mathcal{X}$ of objects in an additive category $\mathcal{C}$, we write $\add \mathcal{X}$ for the additive closure of $\mathcal{X}$, that is, the full subcategory consisting of finite direct sum of direct summands of objects in $\mathcal{X}$.

\medskip
\noindent
\textbf{Use of AI.} After we calculated the projective cover of the first syzygy (see subsection~\ref{subsec:multi-free deg 0 1}) and fed it to ChatGPT 5.6 Sol, it suggested the candidate resolution to us; the first half of Appendix \ref{app:thin-modules} records our retrospective explanation of this suggestion in a more general setting, and the second half (material after Example \ref{eg:non exact cpx from simple}) was pointed out by ChatGPT after we fed it this viewpoint.
Lemma \ref{lem:perturbed-retraction}, which is a key ingredient in the proof of the main result (and the major obstacle when we first attempted the problem), was first shown directly by ChatGPT; the (more intuitive) proof we present now is modified by ourselves following the advice of Claude Opus 5.5.
After we fed ChatGPT with the write-up of Theorem \ref{thm:exactness} and 
Theorem \ref{thm:minimal-resolution}, and exploring various phenomena, ChatGPT 5.6 Sol obtained subsection~\ref{subsec:higher-multiplicities}.
A statement similar to Theorem \ref{thm:pdim=reg} was first speculated by ChatGPT 5.5 when we fed it the companion paper \cite{CX26}, without a sound argument; the correct form and proof of Theorem \ref{thm:pdim=reg} was obtained by ChatGPT 5.6 Sol after feeding it our Theorem \ref{thm:minimal-resolution}.
All mathematical arguments were independently verified and finalized by the authors. The authors take full responsibility for the accuracy, originality, and integrity of the paper.
AI was also used for grammar checking and diagram drawing.

\medskip
\noindent
\textbf{Acknowledgments.}
AC is deeply grateful to Yasuhiko Asao and Shun Wakatsuki on teaching him about the connection between distance algebras and magnitude homolgies, and the many discussions on the subject.
AC is supported by JSPS Grant-in-Aid for Scientific Research (C) 24K06666. BX is supported by the China Scholarship Council (No. 202506040127).


\section{Preliminaries}
\label{sec:preliminaries}

\subsection{Duality and Gorenstein algebras}
\label{subsec:duality-gorenstein}

Let $A$ be a finite projective $R$-algebra.

Recall that a ring $B$ is \defn{Iwanaga--Gorenstein} if it is noetherian on both sides and
\[
 \idim B_B<\infty
 \qquad\text{and}\qquad
 \idim {}_B B<\infty.
\]
A finite projective $R$-algebra $A$ is called a \defn{Gorenstein $R$-algebra} if $A_{\mathfrak p}$ is Iwanaga--Gorenstein for every $\mathfrak p\in\operatorname{Spec}R$ such that $A_{\mathfrak p}\neq0$; see \cite[Definition~1.1]{IK22} and \cite[Definition~6.1]{GIK23}. When $R=\Bbbk$, this is precisely the Iwanaga--Gorenstein property.
Following Buchweitz \cite[Section~7.6]{Buc21}, the \defn{dualizing bimodule} of $A$ is
\[
 D(A):=\Hom_R(A,R) \text{ with action } (a f b)(c):=f(bca)
 \quad(a,b,c\in A,\ f\in D(A)).
\]
When $R$ is Gorenstein, a criterion of Iyengar--Krause \cite[Theorem~4.6]{IK22} (see also \cite[Theorem~6.7]{GIK23}) states that $A$ is a Gorenstein $R$-algebra if and only if
\begin{equation}
 \pdim_A D(A)<\infty
 \qquad\text{and}\qquad
 \pdim_{A^{\operatorname{op}}}D(A)<\infty;
 \label{eq:canonical-finite-pdim}
\end{equation}
equivalently, $D(A)$ admits a bounded resolution by finitely generated projective modules on each side.  For the distance algebras of finite distributive lattices, we will prove condition \eqref{eq:canonical-finite-pdim} without assuming that $R$ is Gorenstein.  When $R$ is Gorenstein, the above criterion then shows that these are Gorenstein $R$-algebras.

An \defn{$A$-lattice} is a right $A$-module whose underlying $R$-module is finitely generated and projective. We write $\lat(A)$ for the category of $A$-lattices, equipped with the exact structure inherited from the category of $A$-modules. Every short exact sequence in this category splits over $R$.

For finite projective algebras over a commutative noetherian ring, we recall the following fact from the theory of lattices over orders; see, for example, \cite[Section~25]{CR81}.
The following formulation can be found in \cite[Section~2]{BBIKP25}; we include a proof as it is straightforward.

\begin{lem}\label{lem:lattice-canonical-duality}
The functor $D$ induces an exact duality
\[
 D:\lat(A)^{\operatorname{op}}
 \xrightarrow{\sim}\lat(A^{\operatorname{op}}).
\]
The category $\lat(A)$ has enough projective and injective objects, given by $\add A$ and $\add D(A)$, respectively. Moreover, the projective dimension of an $A$-lattice in $\lat(A)$ equals its ordinary projective dimension.
\end{lem}

\begin{proof}
See Appendix~\ref{app:lattice-proofs}.
\end{proof}

We call the objects of $\add D(A)$ \defn{relative injective} $A$-modules. For $M\in\lat(A)$, write $\idim^{\lat}_A M$ for its injective dimension in $\lat(A)$, namely, the least length of a coresolution by the relative injectives if it exists, or $\infty$ otherwise. When $R=\Bbbk$, relative injectives are the finitely generated ordinary injectives, and relative injective dimension agrees with ordinary injective dimension.
We say that $A$ is \defn{self-injective} if $A_A$ and ${}_AA$ are relative injective; this holds, for instance, when $A\cong D(A)$ as right $A$-modules (which also implies isomorphism as left $A$-modules). For $R=\Bbbk$ this is the usual notion.  A self-injective algebra is Iwanaga--Gorenstein if $\idim_RR<\infty$, by the following Corollary~\ref{cor:relative-ordinary-injective-dimension}.

\begin{cor}\label{cor:relative-ordinary-injective-dimension}
Let $A$ be a finite projective $R$-algebra. Then
\[
 \idim^{\lat}_A A
 =\pdim_{A^{\operatorname{op}}}D(A),
 \qquad
 \idim^{\lat}_{A^{\operatorname{op}}}A
 =\pdim_A D(A).
\]
If $\idim_R R<\infty$, then
\[
\begin{aligned}
 \idim A_A
 &\leq \idim^{\lat}_A A+\idim_R R
 =\pdim_{A^{\operatorname{op}}}D(A)+\idim_R R,\\
 \idim {}_A A
 &\leq \idim^{\lat}_{A^{\operatorname{op}}}A+\idim_R R
 =\pdim_A D(A)+\idim_R R.
\end{aligned}
\]
In particular, if $\idim_R R$, $\pdim_{A^{\operatorname{op}}}D(A)$, and $\pdim_{A}D(A)$ are all finite, then $A$ is Iwanaga--Gorenstein. 
When $R=\Bbbk$ is a field, both inequalities are equalities.
\end{cor}

\begin{rmk}
For related formulas for self-injective dimensions, see \cite[Proposition~1.1(3)]{GN01} and \cite[Lemma~4.17]{IKWY15}.
\end{rmk}

\begin{proof}
See Appendix~\ref{app:lattice-proofs}.
\end{proof}


\subsection{Distance algebra}

Let $G=(V,E)$ be a connected combinatorial graph, where $V$ is the set of vertices and $E$ is the set of edges. In this paper, we assume that $G$ contains neither loops nor multiple edges.

For vertices $x,y\in V$, we write $\{x,y\}\in E$ if $x$ and $y$ are connected by an edge; in this case, they are called \defn{adjacent}. The \defn{distance} between $x$ and $y$, denoted by $d(x,y)$ (or $d_G(x,y)$ if we want to specify the graph), is the length of a shortest path from $x$ to $y$ on $G$, and the \defn{diameter} of $G$ is the maximal distance in $G$, denoted by $\operatorname{diam}(G):=\max_{x,y\in V}d(x,y)$.

\begin{dfn}
Let $G=(V,E)$ be a graph (without multi-edges and loops), and $Q_G$ be its \defn{double quiver}, i.e. the set of vertices of $Q_G$ is $V$ and there is a pair of arrows
\[\begin{tikzcd}
	x \arrow[r, "{a_{xy}}", bend left] & y \arrow[l, "{a_{yx}}", bend left]
\end{tikzcd} \text{ for each edge }\{x,y\}\in E.\]
The \defn{distance $R$-algebra} of $G$ is defined to be $\Sigma_G:=R Q_G/I_G$, where $I_G$ is the ideal generated by the following relations:
\begin{enumerate}
\item $p=0$, for every path $p$ such that $\ell(p)>d(s(p),t(p))$;
\item $p=q$, for every pair of paths $p,q$ in $Q_G$ with
$$s(p)=s(q),\qquad t(p)=t(q),\qquad
\ell(p)=d(s(p),t(p))=\ell(q).$$
\end{enumerate}
\end{dfn}

By definition, in $\Sigma_G$, a path is non-zero if and only if it is represented by a shortest path in $G$, and any two shortest paths with the same starting and ending vertices are identified. Hence, for any two vertices $x,y\in V$, there is a unique non-zero path from $x$ to $y$ modulo $I_G$. In particular, $e_x\Sigma_Ge_y$ is a free $R$-module of rank one for all $x,y\in V$, and $\Sigma_G$ is free of rank $|V|^2$ over $R$. When $R=\Bbbk$, the algebra $\Sigma_G$ is Schurian.

We take an $R$-basis $\{b_{xy}\mid x,y\in V\}$ for $\Sigma_G$, where $b_{xy}$ denotes the equivalence class of any shortest path from $x$ to $y$ on $G$, with the convention that $b_{xx}=e_x$ is the vertex idempotent associated to $x$. Then multiplication in $\Sigma_G$ is given by
$$
b_{xy}b_{yz}=
\begin{cases}
b_{xz}, & \text{if } d(x,z)=d(x,y)+d(y,z),\\
0, & \text{otherwise}.
\end{cases}
$$
Note that $b_{xy}$ is denoted by $(x,y)$ and by $e_{xy}$ in \cite{AI24} and \cite{AW24}, respectively. Thus, each $\Sigma_G$ is an $R$-algebra.
Note that path length induces a natural grading on $\Sigma_G$. 

The graphs that we consider in this paper are the Hasse diagrams of finite distributive lattices $L$.  We will simplify the notation and write $\Sigma_L$.

For each $x\in V(G)$, set
\[
 P_x:=e_x\Sigma_G,
 \qquad
 I_x:=D(\Sigma_G e_x).
\]
These are projective and relative injective right $\Sigma_G$-modules, respectively, and
\[
 (\Sigma_G)_{\Sigma_G}=\bigoplus_{x\in V(G)}P_x,
 \qquad
 D(\Sigma_G)_{\Sigma_G}=\bigoplus_{x\in V(G)}I_x.
\]

When $R=\Bbbk$, $P_x$ and $I_x$ are indecomposable, and $I_x$ is precisely the injective envelope of $S_x:=P_x/\rad P_x$.
Moreover, every finitely generated projective right $\Sigma_G$-module $Q$ admits a decomposition $Q\cong\bigoplus_{x\in V(G)}P_x^{\oplus m_x}$, where the non-negative integers $m_x$ are uniquely determined (by the Krull--Schmidt theorem).  We call $m_x$ the \defn{multiplicity} of $P_x$ in $Q$ and write $[Q:P_x]:=m_x$. We say that $Q$ is \defn{multiplicity-free} if $[Q:P_x]\leq1$ for every $x\in V(G)$.

Path reversal $\tau\colon\Sigma_G\to\Sigma_G$, $b_{xy}\mapsto b_{yx}$, is an $R$-linear anti-involution fixing every $e_x$. It identifies $\Sigma_G^{\operatorname{op}}\cong\Sigma_G$, and so gives an exact equivalence $M\mapsto\tau_*M$ from right to left $\Sigma_G$-modules (same underlying module, with $a\cdot m:=m\tau(a)$) that preserves projective dimension.
Since $\tau(e_x)=e_x$, there is an isomorphism of left $\Sigma_G$-modules
\begin{equation}
 \tau_*I_x\xrightarrow{\sim}D(e_x\Sigma_G),
 \qquad f\longmapsto\bigl(b\longmapsto f(\tau(b))\bigr).
 \label{eq:tau-canonical-summands}
\end{equation}

Applying Corollary~\ref{cor:relative-ordinary-injective-dimension} to $\Sigma_G$ gives the following criterion. Here $D(\Sigma_G)$ is regarded as the dualizing bimodule with its natural left and right actions.

\begin{cor}\label{cor:distance-gorenstein-criterion}
Let $n:=\max_{x\in V(G)}\pdim_{\Sigma_G}I_x
 \in\mathbb Z_{\geq0}\cup\{\infty\}$. Then
\begin{equation}
 \pdim_{\Sigma_G}D(\Sigma_G)
 =\pdim_{\Sigma_G^{\operatorname{op}}}D(\Sigma_G)
 =\idim^{\lat}_{\Sigma_G}\Sigma_G
 =n.
 \label{eq:distance-canonical-dimensions}
\end{equation}
\begin{enumerate}
    \item If $R$ is Gorenstein, then $\Sigma_G$ is a Gorenstein
$R$-algebra if and only if $n<\infty$.
\item If $R$ is Iwanaga--Gorenstein and $n<\infty$, then $\Sigma_G$ is Iwanaga--Gorenstein and
\[
 \idim(\Sigma_G)_{\Sigma_G}
 =\idim{}_{\Sigma_G}\Sigma_G
 \leq n+\idim_R R.
\]
\item When $R=\Bbbk$ is a field, both self-injective dimensions equal $n$.
\end{enumerate}
\end{cor}

\begin{proof}
The decomposition $D(\Sigma_G)_{\Sigma_G}
 =\bigoplus_{x\in V(G)}I_x$ gives $\pdim_{\Sigma_G}D(\Sigma_G)=n$. Summing \eqref{eq:tau-canonical-summands} over the vertices gives an isomorphism of left $\Sigma_G$-modules
\[
 \tau_*\bigl(D(\Sigma_G)_{\Sigma_G}\bigr)
 \cong{}_{\Sigma_G}D(\Sigma_G).
\]
Since $\tau_*$ preserves projective dimensions, we also have
$\pdim_{\Sigma_G^{\operatorname{op}}}D(\Sigma_G)=n$. 
The remaining assertions follow immediately from Corollary~\ref{cor:relative-ordinary-injective-dimension}, the Gorenstein criterion \eqref{eq:canonical-finite-pdim}, and the equivalence induced by $\tau$.
\end{proof}


\subsection{Simplicial complexes and homologies}
\label{subsec:simplicial-topology}

Unless a coefficient ring is specified, all chain groups and homology groups are taken over $R$.  We refer the reader to \cite{Mun84} for any possible missing details.

A finite \defn{simplicial complex} $\Delta$ on a finite vertex set $V$ is a collection of subsets of $V$, called \defn{faces}, such that $\sigma\in\Delta$ and $\tau\subseteq\sigma$ imply $\tau\in\Delta$. A face $\sigma$ has \defn{dimension} $|\sigma|-1$, and $\Delta_r$ denotes the set of its $r$-dimensional faces. Two degenerate complexes must not be confused: the \defn{void complex} $\emptyset$, which has no faces at all, and the \defn{trivial complex} $\{\emptyset\}$, whose only face is the empty face (of dimension $-1$). Every other complex contains the empty face and at least one vertex.

For $r\geq-1$, we define a free $R$-module
\[
C_r(\Delta):= R\Delta_r = \bigoplus_{\sigma\in \Delta_r} R \langle \sigma\rangle; \text{ in particular, }C_{-1}(\Delta)= \begin{cases}
    R \cdot 1, & \Delta\neq \emptyset;\\
    0, & \Delta=\emptyset;
\end{cases}
\]
where $1:=\langle\,\rangle$ is the shorthand for the basis element corresponding to the empty face. 
Fix a total order on $V$.  Then we can write any $r$-face as an \defn{$r$-chain}  $v_\ast := (v_0<v_1<\cdots <v_r)$. 
For $r\geq0$ the \defn{boundary map} $\partial_r\colon C_r(\Delta)\to C_{r-1}(\Delta)$ is
\[
\partial_r\langle v_\ast\rangle =\partial_r\langle v_0<\cdots<v_r\rangle:= \sum_{i=0}^{r}(-1)^i \langle v_*^{\widehat{i}}\rangle=\sum_{i=0}^{r}(-1)^i\langle v_0<\cdots<\widehat{v_i}<\cdots<v_r\rangle,
\]
where $\widehat{v_i}$ means that $v_i$ is removed, and $v_*^{\widehat{i}}:=( v_0<\cdots<\widehat{v_i}<\cdots<v_r)$.
The map $\partial_0$, which sends $\langle v\rangle\mapsto 1$, is also called \defn{augmentation}. 
This makes 
\[
\widetilde C_\bullet(\Delta)\colon\quad\cdots\longrightarrow C_1(\Delta)\xrightarrow{\partial_1}C_0(\Delta)\xrightarrow{\partial_0}C_{-1}(\Delta)\longrightarrow0
\]
a \defn{chain complex}, i.e. $\partial_{r-1}\partial_r=0$ for all $r>0$, and the \defn{reduced homology} $\widetilde H_r(\Delta)$ for $r\geq-1$ is its homology.
A different choice of total order on $V$ changes the boundary by some sign changes, but the reduced homologies will be canonically isomorphic; the majority of simplicial complexes we consider in this paper come from \emph{order complex} (see the next subsection), where there is a canonical choice on the vertex set.

For clarity, let us note that the void complex has vanishing reduced homology everywhere, whereas the trivial complex has a single non-vanishing reduced homology $\widetilde{H}_{-1}(\{\emptyset\})=R$. 
Also, a simplicial complex $\Delta$ with at least one vertex has $\widetilde H_{-1}=0$; in other words, $\widetilde H_{-1}(\Delta)\neq0$ if and only if $\Delta=\{\emptyset\}$. 

For $r\geq1$, $\widetilde H_r(\Delta)$ is the ordinary simplicial homology $H_r(\Delta)$, i.e. the homology of the ``non-augmented complex" $C_\bullet(\Delta):=(C_r(\Delta))_{r\geq0}$.
For later use, we call $C_{-1}(\Delta)$ the \defn{augmentation part} and $C_\bullet(\Delta)$ the \defn{resolution part} of $\widetilde C_\bullet(\Delta)$. 
If clarification of the ground ring is needed, we add ``$;R$'' inside the brackets, for example $C_\bullet(\Delta;R)$.

\subsection{Links and order complexes}
\label{subsec:topological-tools}

A simplicial complex $K$ is a \defn{cone with vertex $p$} if $p\in V(K)$ and $\sigma\cup\{p\}\in K$ for every $\sigma\in K$. A cone has zero reduced homology over $R$; see, for example, \cite[Theorem~8.2]{Mun84}.

\begin{dfn}\label{def:link,face poset,sd}
\begin{enumerate}
    \item Let $K$ be a finite simplicial complex on a vertex set $V$.  The \defn{link} of a face $F$ in $K$ is defined as
\[
\operatorname{lk}_K F:=\{G\subseteq V\setminus F:F\cup G\in K\}.
\]

    \item 
For a finite poset $(P,<)$, its \defn{order complex} $\Delta(P)$ is the simplicial complex on the vertex set $P$ whose $r$-dimensional faces are the chains $x_0<x_1<\cdots<x_r$ of $P$; we will use this choice of ordering for the associated chain complex $\widetilde{C}_\bullet(\Delta(P))$.  Note that $\Delta(P)$ is never a void complex.

\item Define the \defn{face poset} and the \defn{non-empty-face poset} of $K$ as
\[
\mathcal F(K):=\{\text{faces of }K\}, \quad \mathcal{F}^\times(K):=\mathcal{F}(K)\setminus\{\emptyset\},
\]
respectively, ordered by inclusion.

\item The \defn{barycentric subdivision} $\operatorname{sd}K$ of a finite simplicial complex $K$ is the order complex of the non-empty-face poset of $K$.
\end{enumerate}
\end{dfn}

For a non-void finite simplicial complex $K$, write $|K|$ for its geometric realization.  We identify simplicial homology with the singular homology of the realization.  Barycentric subdivision preserves the geometric realization, and hence preserves reduced homology over $R$.

For a finite set $E$ with $|E|\geq2$, the power set $2^E$ is regarded as a simplicial complex on $E$, namely, the $(|E|-1)$-simplex with boundary $\partial(2^E)=2^E\setminus\{E\}$, so
\begin{equation}
\Delta\bigl(2^E\setminus\{\emptyset,E\}\bigr)
=\Delta\bigl(\mathcal F^\times(\partial(2^E))\bigr)
=\operatorname{sd}(\partial(2^E)),\quad \text{and}\quad 
|\operatorname{sd}(\partial(2^E))|\cong S^{|E|-2}.
\label{eq:boolean-boundary-sphere}
\end{equation}
Thus omitting $E$ passes to the boundary of the simplex, and omitting $\emptyset$ passes to its non-empty-face poset.

\begin{lem}\label{lem:munkres}
Let $\mathcal F(K)_{>F}:=\{F'\in K\mid F\subsetneq F'\}$.  Then there is an isomorphism of simplicial complexes
\[
\Delta\bigl(\mathcal F(K)_{>F}\bigr)
\cong\operatorname{sd}(\operatorname{lk}_K F).
\]
In particular, we have $\widetilde{H}_k(\operatorname{lk}_KF)\cong \widetilde{H}_k(\Delta(\mathcal{F}(K)_{>F}))$ for all $k\ge -1$.
\end{lem}
\begin{proof}
The first part is known in \cite[Lemma 15.3]{Mun84}; we can also obtain it using the poset isomorphism
\newcommand{\mapsfrom}{\mathrel{\reflectbox{$\mapsto$}}}
\[
\mathcal F(K)_{>F}:=\{F'\in K\mid F\subsetneq F'\} \rightarrow
\mathcal{F}^\times(\operatorname{lk}_K F),
\qquad F'\mapsto F'\setminus F, \text{ with inverse } F\cup H \mapsfrom  H.
\]
The second part follows from
\[
\widetilde{H}_k(\operatorname{lk}_KF)\cong \widetilde{H}_k\big(\operatorname{sd}(\operatorname{lk}_KF)\big) \cong \widetilde{H}_k\big(\Delta(\mathcal{F}(K)_{>F})\big),
\]
where the first isomorphism uses the fact that barycentric subdivision preserves reduced homology \cite[Theorem 17.2]{Mun84}.
\end{proof}

We also recall an order homotopy theorem for maps between order complexes.

Let $f:Q\to Q'$ be an order-preserving map between two finite posets $Q,Q'$.  Define, for all $k\ge 0$, a map $f_{\#,k}:C_k(\Delta(Q))\to C_k(\Delta(Q'))$ given by
\begin{equation}\label{eq:order-presv-map-induced-chain-map}
\langle x_0<\cdots<x_k\rangle \mapsto
\begin{cases}
\langle f(x_0)<\cdots<f(x_k)\rangle, & f(x_0),\ldots,f(x_k)\text{ are pairwise distinct},\\
0,&\text{otherwise}.
\end{cases}
\end{equation}
We also set $f_{\#,-1}=\id_R$. Then $f_\#$ is the \defn{induced chain map} on the chain complexes and their augmentations.

\begin{lem}\label{lem:order homotopy}
Let $f,g:Q\to Q'$ be order-preserving maps between finite posets.
If $f(q)\leq g(q)$ for every $q\in Q$, then there is a chain homotopy $h$ between $f_\#,g_\#$, i.e. $\partial h+h\partial=g_\#-f_\#$, given by $h_{-1}:=0$ and
\[
h_k\langle x_0<\cdots<x_k\rangle:=\sum_{j=0}^k(-1)^j\langle f(x_0)\leq\cdots\leq f(x_j)\leq g(x_j)\leq\cdots\leq g(x_k)\rangle\quad \text{ for }k\geq0,
\]
where every summand with a repeated entry is regarded as zero.
\end{lem}
\begin{proof}
This is the chain complex version of \cite[Theorem 10.11]{Bjo95}; see also \cite[Theorem 13.3]{Mun84}.  
\end{proof}

The following consequence of Alexander duality will be useful.

\begin{lem}\label{lem:split-boundary-subset-chains}
Assume that $|P|\geq3$, and let $K$ be a proper subcomplex of
$\Delta\bigl(2^P\setminus\{\emptyset,P\}\bigr)$. Then the boundary map
$\partial_{|P|-2}:C_{|P|-2}(K)\longrightarrow C_{|P|-3}(K)$ has an $R$-linear left inverse. 
\end{lem}
\begin{proof}
Put $d:=|P|-2\ge 1$.
If $C_d(K)=0$, then $\partial_d$ is the zero map from the zero module, which has the left inverse $0$, so there is nothing to prove. This covers in particular the void complex $K=\emptyset$ and the trivial complex $K=\{\emptyset\}$, for which $C_q(K)=0$ for all $q\geq0$. Hence we may assume that $K$ has a $d$-dimensional face; in particular $|K|\neq\emptyset$.
Since $P$ is finite, each $C_q(K)=C_q(K;\mathbb Z)\otimes_{\mathbb Z}R$ is a finitely generated free $R$-module. 
Let
\[
  C^q(K;R):=\Hom_R\bigl(C_q(K),R\bigr),\qquad
  \delta^q:=\Hom_R(\partial_{q+1},R)\colon C^q(K;R)\to C^{q+1}(K;R)
\]
be the cochain complex.
Each $C^q(K;R)\cong C^q(K;\mathbb Z)\otimes_{\mathbb Z}R$ is again finitely generated free over $R$.

We first show that $\delta^{d-1}\colon C^{d-1}(K;\mathbb Z)\to C^{d}(K;\mathbb Z)$ is surjective.
Indeed, since $|\Delta(2^P\setminus\{\emptyset,P\})|\cong S^{d}$ by \eqref{eq:boolean-boundary-sphere} and $K$ is a
proper subcomplex, $|K|$ is a compact, locally contractible, non-empty, proper subspace of $S^d$.
In particular, we have $S^d\setminus|K|\neq\emptyset$, and Alexander duality \cite[Cor.~3.45]{Hat02} gives
\[
  H^{d}(K;\mathbb Z)\;\cong\;\widetilde H^{d}(|K|;\mathbb Z)\;\cong\;\widetilde H_{-1}(S^d\setminus|K|;\mathbb Z)\;=\;0,
\]
where the first isomorphism uses $d\ge1$, and the last equality holds because reduced homology in degree $-1$ vanishes for a non-empty space.
As $K$ has no simplices of dimension strictly larger than $d$, we have $C^{d+1}(K;\mathbb Z)=0$, so $H^d(K;\mathbb Z)=\operatorname{coker}\delta^{d-1}$, which means that $\delta^{d-1}$ is surjective, as claimed.

Now we base change back to $R$.
Since $-\otimes_{\mathbb Z}R$ is right exact, $\delta^{d-1}\colon C^{d-1}(K;R)\to C^d(K;R)$ is still surjective with a free (hence, projective) codomain, so there is an $R$-linear right inverse $s: C^d(K;R)\to C^{d-1}(K;R)$ of $\delta^{d-1}$.

Now the claim follows by dualizing back to the chain groups.
Namely, since the chain groups are finitely generated free over $R$, we have $C_q(K;R)\cong \Hom_R(C^q(K;R),R)$ and $\Hom_R(\delta^{d-1},R)=\partial_{d}$.  Applying $\Hom_R(-,R)$ to $\delta^{d-1}\circ s=\mathrm{id}$ thus yields the required 
$\Hom_R(s,R)\circ\partial_d=\mathrm{id}$.
\end{proof}

\subsection{Finite posets and distributive lattices}

Recall that a \defn{lattice} is a poset $L$ in which every pair $x,y\in L$ has a greatest lower bound, called the \defn{meet} and denoted by $x\wedge y$, and a least upper bound, called the \defn{join} and denoted by $x\vee y$. A lattice is \defn{distributive} if
\[
x\wedge(y\vee z)
 =(x\wedge y)\vee(x\wedge z),
\qquad
x\vee(y\wedge z)
 =(x\vee y)\wedge(x\vee z)
\]
for all $x,y,z\in L$. A finite lattice has a least element $0$ and a greatest element $1$.

For a finite poset $P$, we denote by $\max(P)$ the set of its maximal elements. An \defn{order ideal} of $P$ is a subset $X\subseteq P$ such that $y\in X$ and $x\leq y$ implies $x\in X$. We denote by $\mathcal J(P)$ the set of order ideals of $P$, ordered by inclusion. 
$\mathcal J(P)$ is a finite distributive lattice with $X\wedge Y=X\cap Y$ and $X\vee Y=X\cup Y$; conversely, Birkhoff's fundamental theorem of finite distributive lattices states that every finite distributive lattice is isomorphic to $\mathcal J(P)$ for some finite poset $P$; see, for example, \cite{Gra03}.

The \defn{Hasse graph} of a finite poset $Q$ is the undirected graph whose vertex set is $Q$, in which two elements are adjacent if one covers the other, i.e. an edge is given by $\{x,y\}$ with $x<y$ and there is no $z\in Q$ satisfying $x<z<y$. We write $d_Q(x,y)$, or simply $d(x,y)$, for the corresponding graph distance. Moreover, for subsets $X,Y\subseteq Q$, their \defn{symmetric difference} is $X\SDiff Y:=(X\setminus Y)\sqcup(Y\setminus X)$.

\begin{lem}\label{lem:symmetric-difference-identities}
Let $X,Y,Z$ be subsets of a finite set. Then the following hold.
\begin{enumerate}[label=\textup{(\arabic*)},ref=\arabic*]

\item\label{itm:symdiff-associativity} $X\SDiff Y=Y\SDiff X$, $X\SDiff X=\emptyset$, and $(X\SDiff Y)\SDiff Z=X\SDiff(Y\SDiff Z)$.

\item\label{itm:symdiff-cancellation} $(X\SDiff Y)\SDiff(X\SDiff Z)=Y\SDiff Z$, $X\SDiff(X\SDiff Y)=Y.$

\item\label{itm:symdiff-intersection} Let $m(X,Y,Z):=(X\cap Y)\cup (X\cap Z)\cup (Y\cap Z)$.
Then $$(X\SDiff Y)\cap (X\SDiff Z) = \bigl(X\setminus(Y\cup Z)\bigr)\sqcup\bigl((Y\cap Z)\setminus X\bigr) = X\SDiff\, m(X,Y,Z).$$

\item\label{itm:symdiff-difference} $ (X\SDiff Y)\setminus(X\SDiff Z) =\bigl(X\cap(Z\setminus Y)\bigr) \sqcup\bigl((Y\setminus Z)\setminus X\bigr).$

\item\label{itm:symdiff-inclusion} We have 
\[
(X\SDiff Y)\subseteq(X\SDiff Z)\quad \Longleftrightarrow\quad 
 \begin{cases}
     X\cap(Z\setminus Y)=\emptyset,\\
     Y\setminus Z\subseteq X.
 \end{cases}
\]
In particular, if $Y\subseteq Z$, then $ (X\SDiff Y)\subseteq(X\SDiff Z) \;\Longleftrightarrow\; X\cap(Z\setminus Y)=\emptyset$.

\item\label{itm:symdiff-interval} If $(X\SDiff Z)\subseteq(X\SDiff Y)$, then
\[
\begin{aligned}
 X\SDiff Y&=(X\SDiff Z)\sqcup(Z\setminus Y)\sqcup(Y\setminus Z),\\
 X\SDiff(Y\cap Z)&=(X\SDiff Z)\sqcup(Z\setminus Y),\\
 X\SDiff(Y\cup Z)&=(X\SDiff Z)\sqcup(Y\setminus Z).
\end{aligned}
\]

\item\label{itm:symdiff-cardinality} $|Y|=|X|+|X\SDiff Y|$ if and only if  $X\subseteq Y$.
\end{enumerate}

\end{lem}

\begin{proof}
(1), (2), (3) are elementary and left to the reader.

(4)\; Intersecting $(X\SDiff Y)\setminus(X\SDiff Z)$ with $X$ yields the left set in the decomposition, and removing $X$ yields the second.

(5)\; This follows from \textup{(\ref{itm:symdiff-difference})} as the left-hand side is equivalent to $(X\SDiff Y)\setminus (X\SDiff Z)=\emptyset$.

(6)\; The hypothesis $(X\SDiff Z)\subseteq(X\SDiff Y)$ gives $Z\setminus Y\subseteq X$ and $(Y\setminus Z)\cap X=\emptyset$ by \textup{(\ref{itm:symdiff-inclusion})} (with $Y$ and $Z$ exchanged). Thus removing $Z\setminus Y$ from $Z$, or adjoining $Y\setminus Z$ to $Z$, adds precisely the corresponding set to $X\SDiff Z$. The three identities then follow from
\[
\begin{aligned}
 Y&=\bigl(Z\setminus(Z\setminus Y)\bigr)\sqcup(Y\setminus Z),\\
 Y\cap Z&=Z\setminus(Z\setminus Y),\\
 Y\cup Z&=Z\sqcup(Y\setminus Z).
\end{aligned}
\]
Note that the unions are disjoint because $Z\setminus Y$ and $Y\setminus Z$ are disjoint from $X\SDiff Z$ and from each other.

(7)\; We always have $|X|+|X\SDiff Y|-|Y|=2|X\setminus Y|$, and so $|Y|=|X|+|X\SDiff Y|$ holds precisely when $X\setminus Y=\emptyset$, or equivalently, $X\subseteq Y$.
\end{proof}

\begin{lem}\label{lem:d(X,Y)}
For any $X,Y,Z\in\mathcal J(P)$, one has
\[
d(X,Y)=|X\SDiff Y|, \quad\text{ and }\quad 
d(X,Y)=d(X,Z)+d(Z,Y)
\Longleftrightarrow
X\SDiff Z\subseteq X\SDiff Y.
\]
\end{lem}

\begin{proof}
Every edge in the Hasse graph of $\mathcal J(P)$ changes an order ideal by exactly one element. Hence every path from $X$ to $Y$ has length at least $|X\SDiff Y|$.  
Conversely, one may first remove the elements of $X\setminus Y$, choosing a maximal remaining element at each step, and then add the elements of $Y\setminus X$, choosing a minimal missing element at each step.  
Every intermediate subset is an order ideal, and this gives a path of length $|X\SDiff Y|$. By Lemma~\ref{lem:symmetric-difference-identities}\textup{(\ref{itm:symdiff-cancellation})}, $Z\SDiff Y=(X\SDiff Z)\SDiff(X\SDiff Y)$. Hence the second assertion follows from the distance formula and Lemma~\ref{lem:symmetric-difference-identities}\textup{(\ref{itm:symdiff-cardinality})}, applied to $X\SDiff Z$ and $X\SDiff Y$.
\end{proof}

\begin{lem}[{cf.\ \cite[Lemma~6.3 and Corollary~7.24(2)]{CX26}}]
\label{lem:projective-relative-injective}
For $X\in\mathcal J(P)$, the projective $\Sigma_{\mathcal J(P)}$-module $P_X$ is relative injective if and only if $P\setminus X\in\mathcal J(P)$. In this case, $P_X\cong I_{P\setminus X}$ and $P_{P\setminus X}\cong I_X$.
\end{lem}

\begin{proof}
$\underline{\Rightarrow}$: Suppose that $P_X$ is relative injective, that is, $P_X\in\add D(\Sigma_{\mathcal J(P)})$. Choose a maximal ideal $\mathfrak m$ of $R$. Since $\Sigma_{\mathcal J(P)}$ is free over $R$ of finite rank, we have
\[
D(\Sigma_{\mathcal J(P)})\otimes_R R/\mathfrak m
\cong
\Hom_{R/\mathfrak m}
\bigl(\Sigma_{\mathcal J(P)}\otimes_R R/\mathfrak m,
R/\mathfrak m\bigr).
\]
Thus $P_X\otimes_R R/\mathfrak m$ is injective over $\Sigma_{\mathcal J(P)}\otimes_R R/\mathfrak m$. By \cite[Lemma~6.3 and Corollary~7.24(2)]{CX26}, we obtain $P\setminus X\in\mathcal J(P)$.

$\underline{\Leftarrow}:$ Suppose that $P\setminus X\in\mathcal J(P)$. By Lemma \ref{lem:d(X,Y)}, for every $Y\in\mathcal J(P)$, we have $d(X,Y)+d(Y,P\setminus X)=|P|=d(X,P\setminus X)$, which then implies that $b_{XY}b_{Z,\,P\setminus X}=\delta_{Y,Z}b_{X,\,P\setminus X}$ for any $Y,Z\in\mathcal J(P)$. Define
\[
\phi:P_X\longrightarrow I_{P\setminus X},
\qquad
\phi(p)(q):=
\text{the coefficient of $b_{X,\,P\setminus X}$ in $pq$},
\]
where $p\in P_X=e_X\Sigma_{\mathcal{J}(P)}$ and $q\in\Sigma_{\mathcal J(P)}e_{P\setminus X}$. It is routine to verify that this is a homomorphism of modules. The multiplication  $b_{XY}b_{Z,\,P\setminus X}=\delta_{Y,Z}b_{X,\,P\setminus X}$ implies that $\phi$ sends the basis $\{b_{XY}\}_{Y\in\mathcal J(P)}$ to the dual $R$-basis of $\{b_{Y,\,P\setminus X}\}_{Y\in\mathcal J(P)}$. Thus $P_X\cong I_{P\setminus X}$ and is relative injective. Swapping $X$ and $P\setminus X$ gives $P_{P\setminus X}\cong I_X$.
\end{proof}

\begin{eg}\label{eg:posets-at-most-two-elements}
Up to isomorphism, there are exactly four posets $P$ with $|P|\leq2$. In the following, $A_n$ denotes the path graph on $n$ vertices.
\begin{enumerate}
\item If $P=\emptyset$, then $\mathcal J(P)=\{\emptyset\}$ and $\Sigma_{\mathcal J(P)}\cong R$.  Thus, $\Sigma_{\mathcal{J}(P)}$ is self-injective ($\idim^{\lat} \Sigma_{\mathcal{J}(P)}=0$).

\item If $P=\{\ast\}$, then $\mathcal J(P)=\{\emptyset \subset \{\ast\}\}$ and $\Sigma_{\mathcal J(P)}\cong\Sigma_{A_2}$. The algebra $\Sigma_{A_2}$ is given by the quiver
\[
\begin{tikzcd}[column sep=large]
{1} & {2}
\arrow["{b_{12}}", from=1-1, to=1-2, bend left]
\arrow["{b_{21}}", from=1-2, to=1-1, bend left]
\end{tikzcd} \;\; \text{with relations } b_{12}b_{21}=0=b_{21}b_{12}.
\]
Then, by Lemma~\ref{lem:projective-relative-injective}, $\Sigma_{A_2}$ is self-injective.

\item If $P=\{a,b\}$ with $a,b$ incomparable, then $\mathcal J(P)$ is the Boolean lattice of rank two, which is the direct product of two copies of $\{1<2\}$.  The property $P\setminus X\in\mathcal J(P)$ for every $X\in\mathcal J(P)$ still holds, and so it follows from Lemma \ref{lem:projective-relative-injective} that $\Sigma_{\mathcal{J}(P)}$ is self-injective.
Note that the Hasse graph of $\mathcal{J}(P)$ is a square, which is the Cartesian product of two copies of $A_2$, and so it follows from \cite[Proposition~3.2]{CX26} that
\[
\Sigma_{\mathcal J(P)}\cong\Sigma_{A_2}\otimes_R\Sigma_{A_2},
\]
and one can then deduce that $\Sigma_{\mathcal{J}(P)}$ is self-injective by the fact that the tensor product over the base ring $R$ preserves self-injectivity.

\item If $P=\{a<b\}$, then $\mathcal J(P)=\{\emptyset \subset (a)\subset (a<b)\}\cong \{1<2<3\}$ and $\Sigma_{\mathcal J(P)}\cong\Sigma_{A_3}$. The algebra $\Sigma_{A_3}$ is given by the quiver
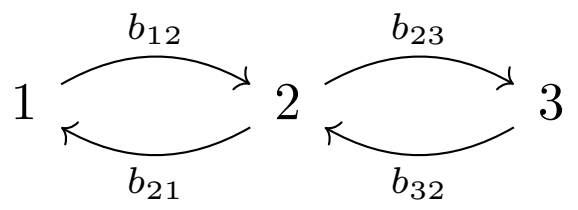
\[
\begin{tikzcd}[column sep=large]
{1} & {2} & {3}
\arrow["{b_{12}}", from=1-1, to=1-2, bend left]
\arrow["{b_{21}}", from=1-2, to=1-1, bend left]
\arrow["{b_{23}}", from=1-2, to=1-3, bend left]
\arrow["{b_{32}}", from=1-3, to=1-2, bend left]
\end{tikzcd}
\]
with every path of length two that returns to its starting vertex equal to zero.

By Lemma~\ref{lem:projective-relative-injective}, $P_1\cong I_3$ and $P_3\cong I_1$, whereas $P_2$ is not relative injective. Let $\{f_i\}_{i=1}^3$ be the $R$-basis of $I_2=D(\Sigma_{A_3}e_2)$ dual to $\{b_{i2}\}_{i=1}^3$. There is an exact sequence
\[
0\longrightarrow P_2
\xrightarrow{\ \iota\ }P_1\oplus P_3
\xrightarrow{\ \pi\ }I_2\longrightarrow0,
\]
where $\iota(p)=(b_{12}p,-b_{32}p)$ and $\pi(p,q)=f_1p+f_3q$. This sequence is a relative injective coresolution of $P_2$ of length one, and $\idim^{\lat}\Sigma_{A_3}=1$. The same equality holds on the left by Corollaries~\ref{cor:relative-ordinary-injective-dimension} and~\ref{cor:distance-gorenstein-criterion}.
For the special case where $R$ is a field, see \cite[Proposition~7.6]{CX26}.
\end{enumerate}
\end{eg}


\subsection{Order complexes of strict upsets}
\label{subsec:strict-upsets}

Let $P$ be a finite poset and fix $X\in\mathcal J(P)$.
We write $\underline{Z}:=X\SDiff Z$ (and $\underline{Z}_X$ if clarification is needed). Consider the poset
\[\mathcal D_X := \{X\SDiff Z \mid Z\in\mathcal J(P)\},\] 
ordered by inclusion.  We write $\le$ for this partial order on $\D_X$, and keep $\subseteq$ as the symbol for the inclusion partial order on $\mathcal{J}(P)$.

For $Z\in\mathcal J(P)$, let $\D_X^{>\underline{Z}}$ be the upset of $\D_X$ containing all elements strictly larger than $\underline{Z}$; we call this the \defn{strict upset above $\underline{Z}$} for convenience. Denote by
\[
\Delta_X:=\Delta(\D_X), \qquad \text{and} \qquad \Delta^{>\underline{Z}}:=\Delta_X^{>\underline{Z}}:=\Delta(\D_X^{>\underline{Z}})
\]
the associated order complexes.  Note that the vertices of $\Delta^{>\underline{Z}}$ are the elements $\underline{T}>\underline{Z}$, and its $r$-dimensional faces are the strict chains $\underline{T_0}<\underline{T_1}<\cdots<\underline{T_r}$ (with $\underline{Z}<\underline{T_0}$), and it is a subcomplex of $\Delta_X$.

\begin{eg}\label{subsec:eg1}
Let $P=\{1<2\}\sqcup\{3\}$.  For a more compact display, we use, for instance, $[ijk]$ to denote $\{i,j,k\}\subset P$, and write $\underline{xy}$ when $\underline{Z}=\{x,y\}$. We consider the case when $X=\{1\}=[1]$.
\[
\begin{array}{r|cccccc}
 Z&[1]&[\emptyset]&[12]&[13]&[3]&[123]\\ \hline
 \underline{Z} = X\SDiff Z&\underline{\emptyset}&\underline{1}&\underline{2}&\underline{3}&\underline{13}&\underline{23}
\end{array}
\]

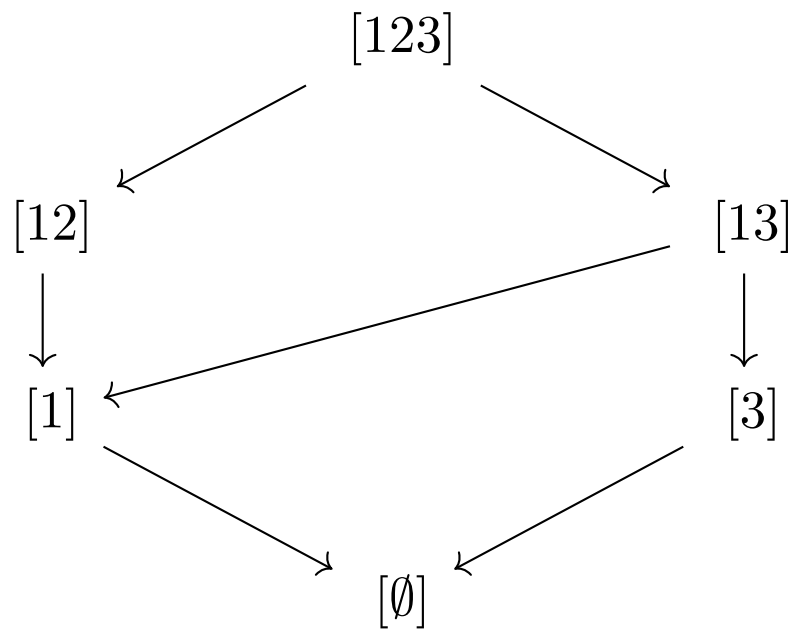
\[
L: \quad \begin{tikzcd}[column sep=large]
 & \ [123]\ar[ld]\ar[rd] & \\
\ [12] \ar[d] & & \ [13] \ar[d]\ar[lld] \\
\ [1] \ar[rd] & & \ [3] \ar[ld] \\
& \ [\emptyset] &
\end{tikzcd} \qquad
\D_X: \quad 
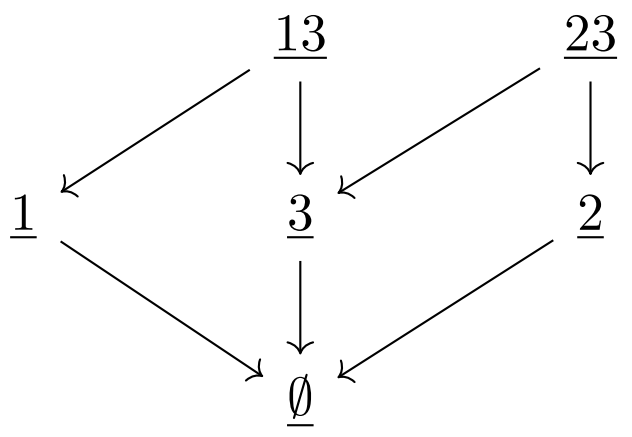
\begin{tikzcd}[column sep=large]
 & \underline{13}\ar[ld]\ar[d] & \underline{23}\ar[d]\ar[dl] \\
 \underline{1}\ar[rd] & \underline{3}\ar[d] & \underline{2}\ar[ld] \\
 & \underline{\emptyset} &
\end{tikzcd}
\]
Note that $\D_X$ is really just a re-orientation of $L$.

The six strict upsets $\D_X^{>\underline{Z}}$ and their order complexes $\Delta_X^{>\underline{Z}}$ are displayed below. In this example, every $\D_X^{>\underline{Z}}$ has height at most one, so each $\Delta_X^{>\underline{Z}}$ is at most $1$-dimensional and may be visualized as the underlying graph of the Hasse diagram of $\Delta_X^{> \underline{Z}}$.

\[
\begin{array}{c|c|c|c}
 Z& \underline{Z} & \D_X^{>\underline{Z}}&\Delta_X^{>\underline{Z}}\\[1mm] \hline
 \ [1]&\underline{\emptyset} &\{\underline{1},\underline{2},\underline{3},\underline{13},\underline{23}\}&
 \underline{1}-\underline{13}-\underline{3}-\underline{23}-\underline{2}\\[1mm]
 \ [\emptyset]& \underline{1}&\{\underline{13}\}& \{\underline{13}\} \\[1mm]
 \ [12]& \underline{2}&\{\underline{23}\}& \{\underline{23}\} \\[1mm]
 \ [13]& \underline{3}&\{\underline{13}, \underline{23}\}&  \{\underline{13}\}\sqcup\{\underline{23}\}\\[1mm]
 \ [3]& \underline{13}&{\emptyset}&\{\emptyset\}\\[1mm]
 \ [123]& \underline{23}&{\emptyset}&\{\emptyset\}
\end{array}
\]

Consequently,
\[
 \widetilde H_r(\Delta^{>\underline{Z}}_{[1]})=0
 \quad\text{except for}\quad
 \widetilde H_0(\Delta^{>\underline{3}}_{[1]})\cong R(\langle \underline{13}\rangle -  \langle\underline{23}\rangle),
 \quad
 \widetilde H_{-1}(\Delta^{>\underline{13}}_{[1]})\cong R,
 \quad
 \widetilde H_{-1}(\Delta^{>\underline{23}}_{[1]})\cong R.
\]
\end{eg}


\section{A projective resolution from strict chains}
\label{sec:candidate-resolution}

We retain the notation and assumptions of subsection~\ref{subsec:strict-upsets}, and let $A$ be the distance algebra $\Sigma_{\mathcal{J}(P)}$.
In this section we construct a complex $\mathcal P_\bullet$ of projective $A$-modules indexed by the strict chains in $\D_X$, equip it with an augmentation to $I_X$, and prove that the augmented complex is exact.

\subsection{The candidate complex}

To distinguish between different copies of the same module in a direct sum, we will use $P_Z^{(i)}$ to denote the copy of $P_Z$ labeled by $i$ in some indexing set.  Similarly, for $p\in P_Z$, we write $p^{(i)}$ for the element $p$ in the $i$-labeled copy of $P_Z$; this is nothing but a compact way to write $\iota_i(p)$, where $\iota_i$ is the canonical inclusion of $P_Z^{(i)}$ into the direct sum that contains it.

For $k\geq0$, we write
\[
 \underline{Z}_\ast=(\underline{Z}_0<\cdots<\underline{Z}_k)
 \in\Delta_{X,k}=\Delta(\mathcal D_X)_k
\]
for a strict chain, where $Z_i$ is the order ideal such that $\underline{Z}_i = Z_i \SDiff X$. For $0\leq i\leq k$, define
\[
 \underline{Z}_\ast^{\widehat{i}}
 :=(\underline{Z}_0<\cdots<\widehat{\underline{Z}_i}<\cdots<\underline{Z}_k)
\]
to be the strict chain obtained by deleting $\underline{Z}_i$. When $k=0$, the chain $\underline{Z}_\ast^{\widehat{0}}$ is empty.

For $k\geq0$, define $\mathcal{P}_k$ to be the direct sum of projective $A$-modules $P_Z$ over all $\underline{Z}_\ast\in\Delta_{X,k}$ with $Z_0 = Z$. Observe that, when $k=0$ (so $\Delta_{X,0}=\mathcal{D}_X=\mathcal{J}(P)$), we have
\[
\mathcal{P}_0 = \bigoplus_{Z\in\mathcal{J}(P)} P_Z^{(\underline{Z})} \cong \bigoplus_{Z\in\mathcal{J}(P)} P_Z\otimes_{R} C_{-1}(\Delta^{>\underline{Z}})\]
via $p^{(\underline{Z})}\mapsto p\otimes 1\in P_Z\otimes_{R}C_{-1}(\Delta^{>\underline{Z}_0})$. On the other hand, when $k>0$, a chain $\underline{Z}_\ast \in \Delta_{X,k}$ can be identified with the $(k-1)$-simplex $\underline{Z}_\ast^{\widehat{0}}$ of $\Delta^{>\underline{Z}}$. This yields an isomorphism
\begin{equation}
\mathcal P_k := \bigoplus_{\underline{Z}_\ast\in \Delta_{X,k}} P_{Z_0}^{(\underline{Z}_\ast)} \cong
\bigoplus_{Z\in\mathcal J(P)}
P_Z\otimes_{R}C_{k-1}(\Delta^{>\underline{Z}}), \qquad  p^{(\underline{Z}_\ast)} \mapsto \begin{cases}
    p\otimes \langle \underline{Z}_\ast^{\widehat{0}} \rangle, & k>0;\\ p\otimes 1, & k=0.
\end{cases}
\label{eq:candidate-terms}
\end{equation}

Note that, for $Z,T\in\mathcal J(P)$ with $\underline{Z}<\underline{T}$ (equivalently, $d(X,T)=d(X,Z)+d(Z,T)$), we have a right $A$-module homomorphism
$$\mu_{T,Z}\colon P_Z\longrightarrow P_T, \qquad p\longmapsto b_{TZ}p,$$ given by left-multiplying $b_{TZ}\in e_TA e_Z$.
For $k>0$, define an $A$-module homomorphism $d_k\colon\mathcal P_k\longrightarrow\mathcal P_{k-1}$ by
\begin{equation}
 d_k\bigl(p^{(\underline{Z}_\ast)}\bigr)
:=
\bigl(\mu_{Z_1,Z_0}(p)\bigr)^{(\underline{Z}_\ast^{\widehat{0}})}
+ \sum_{i=1}^k
(-1)^i  p^{(\underline{Z}_\ast^{\widehat{i}})}.
 \label{eq:candidate-differential}
\end{equation}

\begin{prop}\label{prop:resolution-differential}
The maps $d_k$ satisfy $d_{k-1}d_k=0$ for every $k\geq2$.
\end{prop}
\begin{proof}
This is the special case of Proposition \ref{prop:cpx from module} applied to $\D_X$ viewed as a category.
\end{proof}

\begin{eg}\label{subsec:eg2}
We continue with the example presented in Example \ref{subsec:eg1}. For a strict chain $\underline{Z}_\ast$ in $\D_X$, let $P_{Z_0}^{(\underline{Z}_\ast)}$ be a labeled copy of $P_{Z_0}$.  
In this example,
 \begin{align*}
 \mathcal P_0 & = P_{[1]}^{(\underline{\emptyset})}\oplus P_{\emptyset}^{(\underline{1})}\oplus P_{[12]}^{(\underline{2})}\oplus P_{[13]}^{(\underline{3})}\oplus P_{[3]}^{(\underline{13})}\oplus P_{[123]}^{(\underline{23})},\\
 \mathcal P_1 & = P_{[1]}^{(\underline{\emptyset}<\underline{1})}\oplus P_{[1]}^{(\emptyset<\underline{2})}\oplus P_{[1]}^{(\underline{\emptyset}<\underline{3})}  \oplus P_{[1]}^{(\underline{\emptyset}<\underline{13})}\oplus P_{[1]}^{(\underline{\emptyset}<\underline{23})}\\
 & \quad\oplus P_{\emptyset}^{(\underline{1}<\underline{13})}\oplus P_{[12]}^{(\underline{2}<\underline{23})}\oplus P_{[13]}^{(\underline{3}<\underline{13})}\oplus P_{[13]}^{(\underline{3}<\underline{23})},\\
 \mathcal P_2 & = P_{[1]}^{(\underline{\emptyset}<
 \underline{1}<
 \underline{13})}\oplus P_{[1]}^{(\underline{\emptyset}<
 \underline{3}<
 \underline{13})}\oplus P_{[1]}^{(\underline{\emptyset}<
 \underline{2}<
 \underline{23})}\oplus P_{[1]}^{(\underline{\emptyset}<
 \underline{3}<
 \underline{23})}.
\end{align*}
Note that $\mathcal P_{\ge 3}=0$.  For $k>0$, the differential is
\[
 d_k\bigl(p^{(\underline{Z}_\ast)}\bigr)=\bigl(b_{Z_1Z_0}p\bigr)^{(\underline{Z}_\ast^{\widehat{0}})}+\sum_{i=1}^{k}(-1)^ip^{(\underline{Z}_\ast^{\widehat{i}})}.
\]
We display the resolution more graphically in Figure \ref{fig:complex-picture}, where each component of the maps involved is shown more clearly.  The labeling on each arrow represents the map, with the unlabeled ones representing the identity map.  We group (highlight) the terms coming from the same subcomplex $\Delta^{>\underline{Z}}$ by the same color. 

\begin{figure}
    \centering
    \begin{tikzpicture}[
  >=stealth,
  every node/.style={inner sep=2pt},
  labl/.style={fill=white, inner sep=1pt,
font=\small} ,
  hilit/.style={line width=10pt, opacity=0.5, line cap=round, line join=round}
]
\colorlet{pinkc}{magenta!60!white}
\colorlet{greenc}{green!70!black}

\node[fill=yellow] (P1-0-1-13)      at (2.5,-3.2)  {$P_{[1]}^{(\emptyset<1<13)}$};
\node[fill=yellow] (P1-0-3-13)      at (2.5,-6.4)  {$P_{[1]}^{(\emptyset<3<13)}$};
\node[fill=yellow] (P1-0-3-23)      at (2.5,-9.6)  {$P_{[1]}^{(\emptyset<3<23)}$};
\node[fill=yellow] (P1-0-2-23)      at (2.5,-12.8) {$P_{[1]}^{(\emptyset<2<23)}$};

\node[fill=yellow] (P1-0-1)        at (7.5,-3.2)  {$P_{[1]}^{(\emptyset<1)}$};
\node[fill=cyan] (P0-1-13)        at (7.5,-1.6)  {$P_{\emptyset}^{(1<13)}$};
\node[fill=yellow] (P1-0-13)        at (7.5,-4.8)  {$P_{[1]}^{(\emptyset<13)}$};
\node[fill=green] (P13-3-13)       at (7.5,-6.4)  {$P_{[13]}^{(3<13)}$};
\node[fill=yellow] (P1-0-3)         at (7.5,-8.0)  {$P_{[1]}^{(\emptyset<3)}$};
\node[fill=green] (P13-3-23)       at (7.5,-9.6)  {$P_{[13]}^{(3<23)}$};
\node[fill=yellow] (P1-0-23)        at (7.5,-11.2) {$P_{[1]}^{(\emptyset<23)}$};
\node[fill=yellow] (P1-0-2)         at (7.5,-12.8) {$P_{[1]}^{(\emptyset<2)}$};
\node[fill=orange] (P12-2-23)       at (7.5,-14.4) {$P_{[12]}^{(2<23)}$};

\node[fill=cyan] (P0-1)           at (15,-1.6)   {$P_{\emptyset}^{(1)}$};
\node (P3-13)          at (15,-3.2)   {$P_{[3]}^{(13)}$};
\node[fill=yellow] (P1-0)           at (15,-6.4)   {$P_{[1]}^{(\emptyset)}$};
\node[fill=green] (P13-3)          at (15,-9.6)   {$P_{[13]}^{(3)}$};
\node (P123-23)        at (15,-12.8)  {$P_{[123]}^{(23)}$};
\node[fill=orange] (P12-2)          at (15,-14.4)  {$P_{[12]}^{(2)}$};

\draw[->] (P1-0-1) -- node[labl, pos=0.2]{$b_{\emptyset[1]}$} (P0-1);
\draw[hilit, color=yellow] (P1-0-1) --  (P1-0);
\draw[->] (P1-0-1) -- node[labl, pos=0.1]{$-1$} (P1-0);
\draw[hilit, color=yellow] (P1-0-1-13) -- (P1-0-1);
\draw[->] (P1-0-1-13) -- (P1-0-1);

\draw[->] (P1-0-1-13) -- node[labl]{$b_{\emptyset[1]}$} (P0-1-13);
\draw[hilit, color=yellow] (P1-0-1-13) --  (P1-0-13);
\draw[->] (P1-0-1-13) -- node[labl]{$-1$} (P1-0-13);
\draw[hilit, color=cyan] (P0-1-13) -- (P0-1);
\draw[->] (P0-1-13) -- node[labl]{$-1$} (P0-1);
\draw[->] (P0-1-13) -- node[labl,pos=0.3]{$b_{[3]\emptyset}$} (P3-13);
\draw[->] (P1-0-13) -- node[labl]{$b_{[3][1]}$} (P3-13);
\draw[hilit, color=yellow] (P1-0-13) -- (P1-0);
\draw[->] (P1-0-13) -- node[labl]{$-1$} (P1-0);
\draw[hilit, color=yellow] (P1-0-3-13) --  (P1-0-13);
\draw[->] (P1-0-3-13) -- node[labl]{$-1$} (P1-0-13);
\draw[->] (P1-0-3-13) -- node[labl]{$b_{[13][1]}$} (P13-3-13);
\draw[hilit, color=yellow] (P1-0-3-13) -- (P1-0-3);
\draw[->] (P1-0-3-13) -- (P1-0-3);
\draw[->] (P13-3-13) -- node[labl,pos=0.15]{$b_{[3][13]}$} (P3-13);
\draw[hilit, color=green] (P13-3-13) --  (P13-3);
\draw[->] (P13-3-13) -- node[labl]{$-1$} (P13-3);
\draw[hilit, color=yellow] (P1-0-3) --  (P1-0);
\draw[->] (P1-0-3) -- node[labl]{$-1$} (P1-0);
\draw[->] (P1-0-3) -- node[labl, pos=0.2]{$b_{[13][1]}$} (P13-3);
\draw[hilit, color=yellow] (P1-0-3-23) -- (P1-0-3);
\draw[->] (P1-0-3-23) -- (P1-0-3);

\draw[->] (P1-0-3-23) -- node[labl]{$b_{[13][1]}$} (P13-3-23);
\draw[hilit, color=yellow] (P1-0-3-23) --  (P1-0-23);
\draw[->] (P1-0-3-23) -- node[labl]{$-1$} (P1-0-23);
\draw[hilit, color=green] (P13-3-23) -- (P13-3);
\draw[->] (P13-3-23) -- node[labl,pos=0.4]{$-1$} (P13-3);
\draw[->] (P13-3-23) -- node[labl]{$b_{[123][13]}$} (P123-23);
\draw[hilit, color=yellow] (P1-0-23) --  (P1-0);
\draw[->] (P1-0-23) -- node[labl,pos=0.4]{$-1$} (P1-0);
\draw[hilit, color=yellow] (P1-0-2-23) --  (P1-0-23);
\draw[->] (P1-0-23) -- node[labl]{$b_{[123][1]}$} (P123-23);
\draw[->] (P1-0-2-23) -- node[labl]{$-1$} (P1-0-23);
\draw[hilit, color=yellow] (P1-0-2-23) -- (P1-0-2);
\draw[->] (P1-0-2-23) -- (P1-0-2);
\draw[->] (P1-0-2-23) -- node[labl]{$b_{[12][1]}$} (P12-2-23);
\draw[hilit, color=yellow] (P1-0-2) --  (P1-0);
\draw[->] (P1-0-2) -- node[labl, pos=0.4]{$-1$} (P1-0);
\draw[->] (P1-0-2) -- node[labl, pos=0.2]{$b_{[12][1]}$} (P12-2);
\draw[->] (P12-2-23) -- node[labl, pos=0.3]{$b_{[123][12]}$} (P123-23);
\draw[hilit, color=orange] (P12-2-23) -- (P12-2);
\draw[->] (P12-2-23) -- node[labl]{$-1$} (P12-2);

\end{tikzpicture}

    \caption{$\mathcal{P}_\bullet$ and the contribution from each strict upset subcomplex $\Delta^{>\underline{Z}}$}
    \label{fig:complex-picture}
\end{figure}
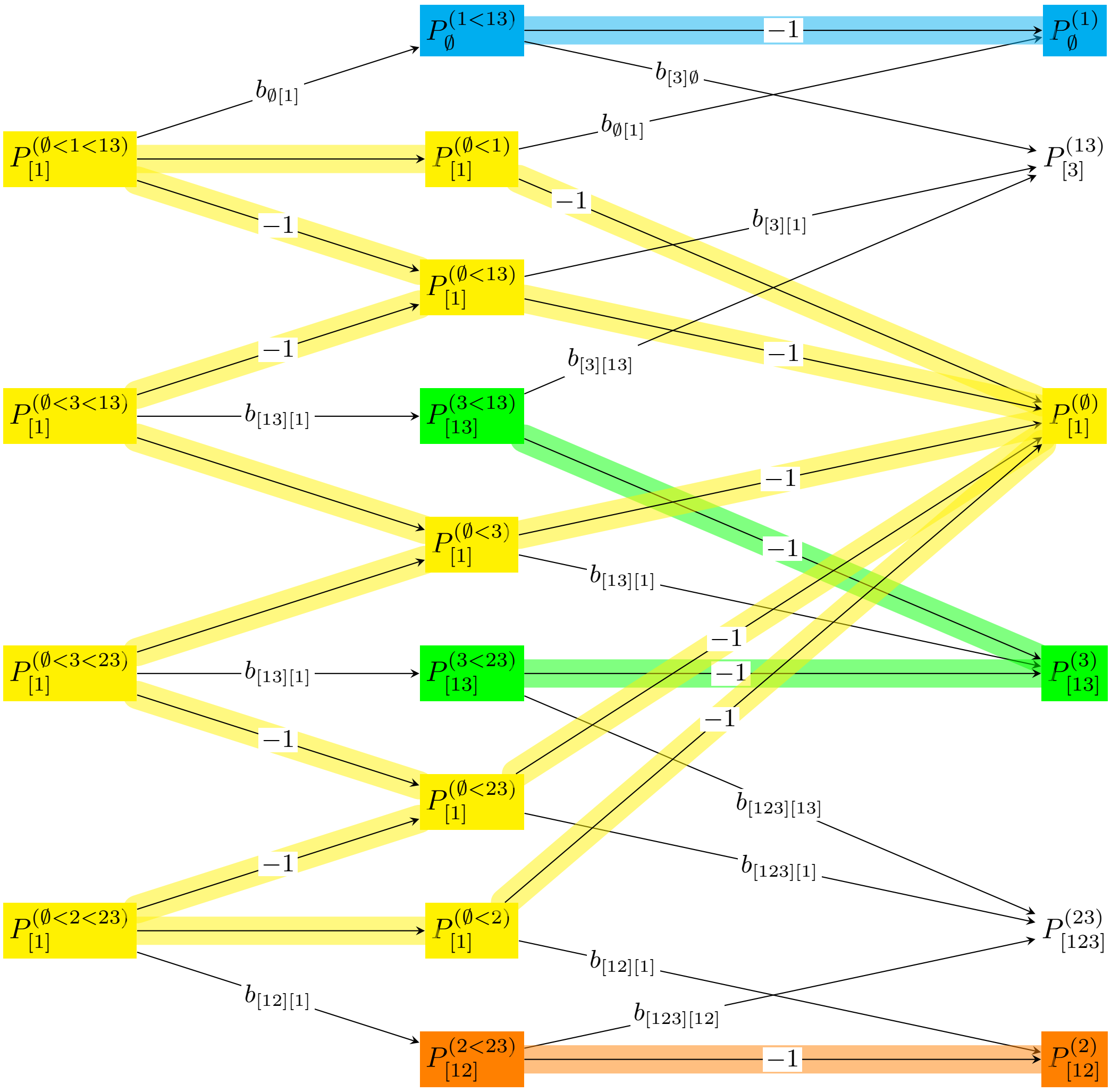
\end{eg}

\medskip 

For each $Y\in\mathcal J(P)$, we have a free $R$-module $I_Xe_Y=D(e_YAe_X)$ of rank $1$.  Let $\phi_Y\colon P_Y\to I_X$ be the homomorphism determined by sending $e_Y$ to the canonical $R$-basis element dual to $b_{YX}\in e_YAe_X$. If $\underline{Z}_0<\underline{Z}_1$, then the corresponding distances are additive and $b_{Z_1Z_0}b_{Z_0X}=b_{Z_1X}$.  Consequently, $\phi_{Z_1}\mu_{Z_1,Z_0}=\phi_{Z_0}$ whenever $\underline{Z}_0<\underline{Z}_1$. Define
\[
 d_0:=\bigl(\phi_Y\bigr)_{Y\in\mathcal J(P)}\colon\mathcal P_0=\bigoplus_{Y\in\mathcal J(P)}P_Y\longrightarrow I_X.
\]
For \(\underline{Z}_0<\underline{Z}_1\), one has $d_0d_1\bigl(p^{(\underline{Z}_0<\underline{Z}_1)}\bigr)=\phi_{Z_1}(b_{Z_1Z_0}p)-\phi_{Z_0}(p)=0$. Hence, by Proposition \ref{prop:resolution-differential}, we obtain a chain complex
\begin{equation}
 \widetilde{\mathcal{P}}_\bullet: \quad \cdots\longrightarrow\mathcal P_2\xrightarrow{d_2}\mathcal P_1\xrightarrow{d_1}\mathcal P_0\xrightarrow{d_0}I_X\longrightarrow0.
 \label{eq:candidate-complex}
\end{equation}
We call $\widetilde{\mathcal P}_\bullet$ the \defn{augmented candidate complex} and $\mathcal{P}_\bullet$ its resolution part.

The key result we want to show in this section is the following.
\begin{thm}\label{thm:exactness}
The complex $\widetilde{\mathcal{P}}_\bullet$ is exact, i.e. $\mathcal{P}_\bullet$ defines a projective resolution of $I_X$.
\end{thm}

Note that the projective resolution constructed above gives an upper bound on $\pdim_A I_X$ through the lengths of strict chains in $\D_X$. Together with Corollary~\ref{cor:distance-gorenstein-criterion}, this yields the Iwanaga--Gorenstein property of $A$.

\begin{cor}\label{cor:distributive-ig}
Assume that $P\neq\emptyset$. Then, for every $X\in\mathcal J(P)$,
\[
\pdim_{\Sigma_{\mathcal J(P)}}I_X\leq |P|-1.
\]
In particular, if $R$ is Gorenstein, then $\Sigma_{\mathcal J(P)}$ is a Gorenstein $R$-algebra. If $R$ is Iwanaga--Gorenstein, then $\Sigma_{\mathcal J(P)}$ is Iwanaga--Gorenstein and
\[
\idim(\Sigma_{\mathcal J(P)})_{\Sigma_{\mathcal J(P)}}
=\idim{}_{\Sigma_{\mathcal J(P)}}\Sigma_{\mathcal J(P)}
\leq |P|-1+\idim_R R.
\]
\end{cor}

\begin{proof}
Fix $X\in\mathcal J(P)$. If $P\setminus X\in\mathcal J(P)$, then $I_X\cong P_{P\setminus X}$ by Lemma~\ref{lem:projective-relative-injective}, so $I_X$ is projective. Otherwise $P\notin\mathcal D_X$, so every strict chain in $\mathcal D_X$ has at most $|P|$ members and $\mathcal P_k=0$ for $k\geq|P|$. Theorem~\ref{thm:exactness} therefore gives $\pdim_{\Sigma_{\mathcal J(P)}}I_X\leq|P|-1$. The remaining assertions follow from Corollary~\ref{cor:distance-gorenstein-criterion}.
\end{proof}

\subsection{Exactness of the candidate complex}
To show that the candidate complex $\widetilde{\mathcal{P}}_\bullet$ is exact, we do this locally, i.e. show the exactness of the complex 
\begin{equation*}
 \widetilde{\mathcal P}_\bullet e_W:\quad
 \cdots\longrightarrow\mathcal P_1e_W
 \longrightarrow\mathcal P_0e_W
 \longrightarrow I_Xe_W\longrightarrow 0
\end{equation*}
of $R$-modules obtained by right-multiplying by the idempotent $e_W$ everywhere.

The first lemma records the multiplication criterion used in both steps of the argument.

\begin{lem}\label{lem:multiplication-criterion}
Let $U,V,W\in\mathcal J(P)$ with $\underline{U}\leq \underline{V}$.  Then $b_{VU}b_{UW}\neq0 \; \Leftrightarrow\; (\underline{V}\setminus \underline{U})\cap \underline{W}=\emptyset$.
\end{lem}

\begin{proof}We have
\[
 b_{VU}b_{UW}\neq0
 \quad\Longleftrightarrow\quad W\,\SDiff\, U\subseteq W\,\SDiff\, V
 \quad\Longleftrightarrow\quad
 \underline{W}\,\SDiff\,\underline{U}\subseteq\underline{W}\,\SDiff\,\underline{V},
\]
where the first equivalence is Lemma~\ref{lem:d(X,Y)}, and the second follows from Lemma~\ref{lem:symmetric-difference-identities}\textup{(\ref{itm:symdiff-cancellation})}. The claim then follows from Lemma~\ref{lem:symmetric-difference-identities}\textup{(\ref{itm:symdiff-difference})} by taking $(X,Y,Z)=(W,V,U)$.
\end{proof}

When restricting to the local part $\widetilde{\mathcal P}_\bullet e_W$, given how the elements in each indecomposable projective are defined, we obtain $R$-modules whose bases are controlled by chains of $\Delta(\D_X)$ bounded above by $\underline{W}$.
It is natural to expect that the resulting complex of $R$-modules is related to the chain complex associated to $$\D_X^{\le \underline{W}} := \{\underline{Z}\in\D_X\mid \underline{Z}\le \underline{W}\}.$$  
In particular, we may use tools from the study of order complexes, namely, Lemma~\ref{lem:order homotopy}, to obtain the desired exactness.

Our strategy is as follows.  The first step (Lemma \ref{lem:evaluated-differential}) is to show that $\widetilde{\mathcal P}_\bullet e_W$ is actually isomorphic to the chain complex $C_\bullet(\Delta(\D_X))$, not with the usual differential $\partial$ but with a `perturbed'  differential $\partial-\omega$.  Then we show in Lemma~\ref{lem:perturbed-retraction} that this perturbed chain complex is homotopy equivalent to the perturbed chain complex associated to $\D_X^{\le \underline{W}}$.  Finally, we show that this smaller perturbed complex is acyclic (hence, an exact sequence) in Lemma \ref{lem:perturbed-lower-exact}.

With a view to using Lemma \ref{lem:order homotopy}, we consider
\[
r_W\colon\D_X\longrightarrow\D_X^{\leq\underline{W}},\qquad
 \underline{Z}\longmapsto\underline{Z}\cap\underline{W}=X\SDiff\,m(X,W,Z)
\]
where the equality is Lemma~\ref{lem:symmetric-difference-identities}\textup{(\ref{itm:symdiff-intersection})}; since $m(X,W,Z)$ is an order ideal, $r_W$ is well-defined. Clearly $r_W$ is order-preserving, $r_W(\underline{Z})\leq\underline{Z}$, and $r_W$ is the identity on $\D_X^{\leq\underline{W}}$.

\begin{lem}\label{lem:evaluated-differential}
Consider the $R$-linear maps
\begin{align*}
\theta_{-1}\colon &I_Xe_W\to C_{-1}(\Delta_X)=R, & b_{WX}^* &\mapsto 1, \\
\theta_k\colon& \mathcal P_ke_W\to C_k(\Delta_X), & (b_{Z_0W})^{(\underline{Z}_\ast)}&\mapsto\langle\underline{Z}_\ast\rangle \qquad (k\ge 0),
\end{align*} 
where $b_{WX}^*$ is the basis of $D(e_W\Sigma e_X)$ dual to $b_{WX}$.
Define $\omega_k\colon C_k(\Delta_X)\to C_{k-1}(\Delta_X)$ by
\[
 \omega_k\langle\underline{Z}_0<\cdots<\underline{Z}_k\rangle:=
 \begin{cases}
 1, & k=0 \text{ and } \underline{W}\not\subseteq \underline{Z}_0,\\
 \langle\underline{Z}_1<\cdots<\underline{Z}_k\rangle, & k>0 \text{ and } r_W(\underline{Z}_0)\neq r_W(\underline{Z}_1),\\
 0, & \text{otherwise}.
 \end{cases}
\]
Then $\theta\colon\widetilde{\mathcal P}_\bullet e_W\to\bigl(\widetilde C_\bullet(\Delta_X),\partial-\omega\bigr)$ is an isomorphism of complexes of $R$-modules.
\end{lem}

\begin{proof}
Each $\theta_k$ is an isomorphism, since $(P_{Z_0}e_W)^{(\underline{Z}_\ast)}=R\,(b_{Z_0W})^{(\underline{Z}_\ast)}$ and $I_Xe_W=D(e_WAe_X)$ are free of rank one. 
Hence, we only need to show that $\theta_{k-1}d_k=(\partial_k-\omega_k)\theta_k$ for $k\geq0$.

Let $k\geq1$. Writing out \eqref{eq:candidate-differential} explicitly yields
\begin{equation}
 d_k\bigl((b_{Z_0W})^{(\underline{Z}_\ast)}\bigr)
 =\bigl(b_{Z_1Z_0}b_{Z_0W}\bigr)^{(\underline{Z}_\ast^{\widehat{0}})}
 +\sum_{i=1}^k(-1)^i
 (b_{Z_0W})^{(\underline{Z}_\ast^{\widehat{i}})}.
 \label{eq:evaluated-differential}
\end{equation}
We compare this with $\partial_k\langle\underline{Z}_\ast\rangle=\sum_{i=0}^k(-1)^i\langle\underline{Z}_\ast^{\widehat{i}}\rangle$ term by term.

For $i\geq1$, the chain $\underline{Z}_\ast^{\widehat{i}}$ starts with $\underline{Z}_0$, so $\theta_{k-1}$ maps $(b_{Z_0W})^{(\underline{Z}_\ast^{\widehat{i}})}$ to $\langle\underline{Z}_\ast^{\widehat{i}}\rangle$. These terms agree.

For $i=0$, the chain $\underline{Z}_\ast^{\widehat{0}}$ starts with $\underline{Z}_1$, so $\theta_{k-1}$ maps $(b_{Z_1W})^{(\underline{Z}_\ast^{\widehat{0}})}$ to $\langle\underline{Z}_\ast^{\widehat{0}}\rangle$. 
We have
\[
 b_{Z_1Z_0}b_{Z_0W}=b_{Z_1W}
 \iff(\underline{Z}_1\setminus\underline{Z}_0)\cap\underline{W}=\emptyset
 \iff r_W(\underline{Z}_0)=r_W(\underline{Z}_1),
\]
where the first equivalence is Lemma~\ref{lem:multiplication-criterion} with $(U,V)=(Z_0,Z_1)$, and the second holds as $\underline{Z}_0\subseteq\underline{Z}_1$. Hence $\theta_{k-1}$ sends the term $(b_{Z_1Z_0}b_{Z_0W})^{(\underline{Z}_\ast^{\widehat{0}})}$ of \eqref{eq:evaluated-differential} to
\begin{align*}
\theta_{k-1}\big( (b_{Z_1Z_0}b_{Z_0W})^{(\underline{Z}_\ast^{\widehat{0}})} \big) &= \begin{cases}
 \langle\underline{Z}_\ast^{\widehat{0}}\rangle, & \text{if }r_W(\underline{Z}_0)=r_W(\underline{Z}_1),\\
 0, & \text{otherwise}
 \end{cases}
 \\
 &=\quad
 \langle\underline{Z}_\ast^{\widehat{0}}\rangle-\omega_k\langle\underline{Z}_\ast\rangle.
\end{align*}
This shows that $\theta_{k-1}d_k=(\partial_k-\omega_k)\theta_k$ holds for $k\ge 1$.

For the case when $k=0$, recall that $\phi_Z:P_Z\to I_X$ is the homomorphism that sends $e_Z$ to the element $b_{ZX}^*\in D(e_Z \Sigma e_X )\subset I_X$ dual to $b_{ZX}\in e_Z \Sigma e_X$, and observe that
\[
 \phi_Z(b_{ZW})\neq0
 \iff b_{ZW}b_{WX}=b_{ZX}
 \iff d(X,Z)=d(X,W)+d(W,Z)
 \iff \underline{W}\subseteq\underline{Z},
\]
where the last equivalence is Lemma~\ref{lem:d(X,Y)}.
Then we obtain 
\[
\theta_{-1}d_0\bigl((b_{ZW})^{(\underline{Z})}\bigr) = \theta_{-1}\big(\phi_Z(b_{ZW})\big) = \begin{cases}
    1, &\text{if }\underline{W}\subseteq\underline{Z},\\
    0, &\text{otherwise}.
\end{cases}
\]
This is exactly $(\partial_0-\omega_0)\langle\underline{Z}\rangle$ as required.
\end{proof}

\begin{lem}\label{lem:perturbed-retraction}
Consider the chain complex $\widetilde{C}_\bullet:=(\widetilde C_\bullet(\Delta_X), \delta:=\partial-\omega)$ of Lemma~\ref{lem:evaluated-differential}.
\begin{enumerate}[(1)]
\item $\widetilde{C}^{\leq\underline{W}}_\bullet:=\widetilde{C}_\bullet(\Delta(\D_X^{\leq\underline{W}}))$ is a subcomplex of $\widetilde{C}_\bullet$, and its differential $\delta|_{\widetilde{C}^{\leq\underline{W}}_\bullet}$ is given by
\[
 \delta_k\langle\underline{Z}_\ast\rangle=(\partial_k-\omega_k)\langle\underline{Z}_\ast\rangle=
 \begin{cases}
 \sum_{i=1}^k(-1)^i\langle\underline{Z}_\ast^{\widehat{i}}\rangle, & k\geq1,\\
 1, & k=0 \text{ and }\underline{Z}=\underline{W},\\
 0, & k=0 \text{ and }\underline{Z}\neq \underline{W}.
 \end{cases}
\]
\item The inclusion $\mathrm{inc}\colon \widetilde{C}^{\leq\underline{W}}_\bullet\to \widetilde{C}_\bullet$ is a homotopy equivalence, with homotopy inverse given by the induced map $(r_W)_\#$.
\end{enumerate}
\end{lem}

\begin{proof}
(1) This amounts to showing that $\delta_k(\widetilde{C}_k^{\le\underline{W}}) \subseteq \widetilde{C}_{k-1}^{\le\underline{W}}$.  First consider the case $k\geq1$, and let $\underline{Z}_\ast=(\underline{Z}_0<\cdots<\underline{Z}_k)$ be a chain in $\D_X^{\leq\underline{W}}$. 
Since $r_W$ is the identity on $\D_X^{\leq\underline{W}}$, we have $r_W(\underline{Z}_0)=\underline{Z}_0\neq\underline{Z}_1=r_W(\underline{Z}_1)$, so $\omega_k\langle\underline{Z}_\ast\rangle=\langle\underline{Z}_\ast^{\widehat{0}}\rangle$. This gives the displayed formula for $k\geq1$. In particular, as each $\underline{Z}_\ast^{\widehat{i}}$ is still bounded by $\underline{W}$, $\delta_k(\widetilde{C}_k^{\le\underline{W}})\subseteq \widetilde{C}_{k-1}^{\leq\underline{W}}$ holds.

For $k=0$ and $\underline{Z}\in\D_X^{\leq\underline{W}}$, we have $\underline{W}\subseteq\underline{Z}$ if and only if $\underline{Z}=\underline{W}$, which gives the formula for $\delta_0=\partial_0-\omega_0$; it takes values in $\widetilde{C}_{-1}^{\leq\underline{W}}=R$ as required.

\medskip

(2) Since $r_W$ is the identity on $\D_X^{\leq\underline{W}}$, the composition $(r_W)_\#\circ\mathrm{inc}$ is the identity chain map on the unperturbed complex $(\widetilde{C}_\bullet(\Delta(\D_X^{\le \underline{W}})),\partial)$.  Apply Lemma~\ref{lem:order homotopy} to $f:=r_W$, regarded as a map $\D_X\to\D_X$, and $g:=\id_{\D_X}$. Then $f\leq g$, $f_\#=\mathrm{inc}\circ(r_W)_\#$, $g_\#=\id$, and there exists $h$ with $\partial h+h\partial=\id-\mathrm{inc}\circ(r_W)_\#$, which is explicitly given by $h_{-1}:=0$ and
\[
 h_k\langle\underline{Z}_0<\cdots<\underline{Z}_k\rangle
 =\sum_{j=0}^k(-1)^{j}
 \langle r_W(\underline{Z}_0)\leq\cdots\leq r_W(\underline{Z}_{j})\leq\underline{Z}_{j}<\cdots<\underline{Z}_k\rangle
 \qquad(k\geq0),
\]
where summands with a repeated entry are zero. Therefore
\[
 \delta h+h\delta=(\partial h+h\partial)-(\omega h+h\omega)=\id-\mathrm{inc}\circ(r_W)_\#-(\omega h+h\omega).
\]

We claim that $\omega_{k+1}h_k\langle\underline{Z}_\ast\rangle+h_{k-1}\omega_k\langle\underline{Z}_\ast\rangle=0$ for all $k$-chains $\underline{Z}_\ast\in \Delta_{X,k}$.  
Indeed, by the definition of $h_k$, the first two entries of the $j$-th summand of $h_k\langle\underline{Z}_\ast\rangle$ are
\begin{equation}\label{eq:first two entries}
 \begin{cases}
 r_W(\underline{Z}_0)\leq\underline{Z}_0, & \text{if }j=0,\\
 r_W(\underline{Z}_0)\leq r_W(\underline{Z}_1), & \text{if }j\geq1.
 \end{cases}
\end{equation}
For $j=0$, applying $r_W$ to both entries yields $r_W(\underline{Z})$, so $\omega_{k+1}$ kills this summand.  In particular, when $k=0$, by combining with $h_{-1}=0$ we obtain
\[
\omega_{1}h_0\langle \underline{Z}_0 \rangle + h_{-1}\omega_0\langle \underline{Z}_0\rangle = \omega_1 \bigl( \langle r_W(\underline{Z}_0) < \underline{Z}_0 \rangle \bigr) = 0,
\]
as required.

Consider the case when $k\geq1$. 
For $j\geq1$ in \eqref{eq:first two entries}, both entries are fixed by $r_W$, so the summand is zero if $r_W(\underline{Z}_0)=r_W(\underline{Z}_1)$, and otherwise $\omega_{k+1}$ deletes $r_W(\underline{Z}_0)$. Hence, if $r_W(\underline{Z}_0)=r_W(\underline{Z}_1)$, then $\omega_{k+1}h_k\langle\underline{Z}_\ast\rangle=0=h_{k-1}\omega_k\langle\underline{Z}_\ast\rangle$. Otherwise
\begin{align*}
 \omega_{k+1}h_k\langle\underline{Z}_\ast\rangle
 &=\sum_{j=1}^k(-1)^{j}
 \langle r_W(\underline{Z}_1)\leq\cdots\leq r_W(\underline{Z}_{j})\leq\underline{Z}_{j}<\cdots<\underline{Z}_k\rangle\\
 &=-h_{k-1}\langle\underline{Z}_1<\cdots<\underline{Z}_k\rangle
 =-h_{k-1}\omega_k\langle\underline{Z}_\ast\rangle,
\end{align*}
where the second equality follows from the definition of $h_{k-1}$, and the third is the definition of $\omega_k$.  This proves our claim.

Finally, by the above claim, we get that $\delta h+h\delta=\id-\mathrm{inc}\circ(r_W)_\#$, which means that
\[\delta(\id-\mathrm{inc}\circ(r_W)_\#)=\delta h\delta=(\id-\mathrm{inc}\circ(r_W)_\#)\delta,\] so the composition $\mathrm{inc}\circ (r_W)_\#$ commutes with the perturbed differential $\delta$. 
Since $\mathrm{inc}$ is an injective chain map by (1), $(r_W)_\# : \widetilde{C}_\bullet \to \widetilde{C}_\bullet^{\le \underline{W}}$ is also a chain map on the perturbed complexes. 
In particular, as chain maps on the perturbed complexes, $(r_W)_\#\circ\mathrm{inc}=\id_{\widetilde{C}_\bullet^{\le\underline{W}}}$ and $\mathrm{inc}\circ(r_W)_\#$ is homotopic to $\id$.
\end{proof}

\begin{lem}\label{lem:perturbed-lower-exact}
The complex $\widetilde{C}^{\leq\underline{W}}_\bullet$ of Lemma~\ref{lem:perturbed-retraction}\textup{(1)} is exact.
\end{lem}

\begin{proof}
Write $\Delta^{\le \underline{W}}:= \Delta(\D_X^{\le \underline{W}})$.
Then we have a decomposition of $R$-modules for all $k\ge 0$:
\[
\widetilde{C}_k(\Delta^{\leq\underline{W}})=\bigoplus_{\underline{Y}\in\D_X^{\leq\underline{W}}}K_{Y,k}, \;\;\text{ where }\;\; K_{Y,k}:= \bigoplus_{\substack{\underline{Z}_\ast\in \Delta^{\le \underline{W}}_k\\\underline{Z}_0=\underline{Y}}}  R\langle \underline{Z}_\ast \rangle.
\]

First observe that $K_{W,k}=0$ for $k\ge 1$, and that $\delta_0|_{K_{W,0}}:K_{W,0}\xrightarrow{\sim} \widetilde{C}_{-1}(\Delta^{\le\underline{W}}) = R$ is an isomorphism.  Indeed, there cannot be a chain in $\D_X^{\leq\underline{W}}$ with at least two members and starts with the maximum $\underline{W}$, which gives us $K_{W,>0}=0$.  On the other hand, $K_{W,0}=R\langle W\rangle$ by definition, and $\delta_0\langle W\rangle = 1\in R$ by Lemma \ref{lem:perturbed-retraction}(1).

Next, we claim that
\[
 \Im\delta_k|_{K_{Y,k}} \,\text{ is }\,
 \begin{cases}
 \text{in } K_{Y,k-1}, & \text{if }k\geq1,\\
 0, & \text{if }k=0 \text{ and } \underline{Y}<\underline{W}.
 \end{cases}
\]
Consequently, combining with the observation from the previous paragraph (which handles the case $k=0$ with $\underline{Y}=\underline{W}$), we have the following decomposition of complexes of $R$-modules 
\begin{equation}
 \widetilde{C}^{\leq\underline{W}}_\bullet=\Bigl(\bigoplus_{\underline{Y}<\underline{W}}K_{Y,\bullet}\Bigr)\oplus \underbrace{\Bigl(K_{W,0}\xrightarrow[\delta_0]{\ \sim\ }R\Bigr)}_{K_{W,\bullet}}.
 \label{eq:perturbed-lower-decomposition}
\end{equation}
Indeed, the $k=0$ part follows by applying Lemma~\ref{lem:perturbed-retraction}\textup{(1)} directly, whereas for the $k\geq1$ part, we observe from Lemma~\ref{lem:perturbed-retraction}\textup{(1)} that $\delta_k\langle\underline{Z}_\ast\rangle$ is a linear combination of the faces $\langle\underline{Z}_\ast^{\widehat{i}}\rangle$ with $i\geq1$, all of which start with the same member as $\underline{Z}_\ast$.  Hence $\delta_k(K_{Y,k})\subseteq K_{Y,k-1}$ as claimed.

With the decomposition \eqref{eq:perturbed-lower-decomposition} in hand, it remains to show that $K_{Y,\bullet}$ is acyclic for every $\underline{Y}<\underline{W}$.  Fix $\underline{Y}<\underline{W}$, and consider the subposet $(\underline{Y},\underline{W}]:=\{\underline{T}\in\D_X\mid\underline{Y}<\underline{T}\leq\underline{W}\}$.
For $k\ge 1$, deleting the zeroth face yields a bijection:
\[
\{ \underline{Z}_\ast \in \Delta_k^{\le \underline{W}} \mid \underline{Z}_0=\underline{Y} \} \xrightarrow{1:1} \Delta\bigl((\underline{Y},\underline{W}]\bigr)_{k-1}, \qquad \langle \underline{Z}_\ast\rangle \longmapsto \langle \underline{Z}_\ast^{\widehat{0}} \rangle.
\]
This extends to $R$-linear isomorphisms $\psi_{Y,k}\colon K_{Y,k}\xrightarrow{\ \sim\ \ }C_{k-1}\bigl(\Delta((\underline{Y},\underline{W}])\bigr)$ for all $k\ge 1$.
On the other hand, we also have an isomorphism $\psi_{Y,0}:K_{Y,0} = R\langle \underline{Y}\rangle \to C_{-1}\bigl(\Delta((\underline{Y},\underline{W}])\bigr)=R$ that
sends $\langle\underline{Y}\rangle$ to $1\in R$.

We claim that $\psi_{Y,k-1}\circ\delta_k=-\partial_{k-1}\circ\psi_{Y,k}$ for $k\geq1$, i.e. they define an isomorphism of complexes
\[
 \psi_Y\colon K_{Y,\bullet}\xrightarrow{\ \sim\ }\widetilde C_\bullet\bigl(\Delta((\underline{Y},\underline{W}])\bigr)[1],
\]
where $[1]$ is the suspension (which lowers the homological degrees of the chain groups by one and negates the differential).
Indeed, for $k\geq1$ we have
\begin{align*}
 \psi_{Y,k-1}\,\delta_k\langle\underline{Y}<\underline{Z}_1<\cdots<\underline{Z}_k\rangle
 &=\psi_{Y,k-1}\Bigl(\sum_{i=1}^k(-1)^i\langle\underline{Y}<\underline{Z}_1<\cdots<\widehat{\underline{Z}_i}<\cdots<\underline{Z}_k\rangle\Bigr)\\
 &=\sum_{i=1}^k(-1)^i\langle\underline{Z}_1<\cdots<\widehat{\underline{Z}_i}<\cdots<\underline{Z}_k\rangle\\
 &=-\partial_{k-1}\langle\underline{Z}_1<\cdots<\underline{Z}_k\rangle\\
 &=-\partial_{k-1}\,\psi_{Y,k}\langle\underline{Y}<\underline{Z}_1<\cdots<\underline{Z}_k\rangle,
\end{align*}
where the first equality is Lemma~\ref{lem:perturbed-retraction}\textup{(1)}, and the third holds because $\underline{Z}_i$ is the $(i-1)$-st member of $\underline{Z}_1<\cdots<\underline{Z}_k$ (so its deletion has sign $(-1)^{i-1}$ in $\partial_{k-1}$).

Finally, note that $\widetilde{C}_\bullet\bigl(\Delta((\underline{Y},\underline{W}])\bigr)$ is an exact sequence, because $\underline{W}$ is the maximum of $(\underline{Y},\underline{W}]$, which makes its order complex a cone.  Thus, the isomorphism $K_{Y,\bullet}\cong \widetilde{C}_\bullet\bigl(\Delta((\underline{Y},\underline{W}])\bigr)[1]$ implies that every $K_{Y,\bullet}$ is acyclic, as required.
\end{proof}

\begin{proof}[Proof of Theorem~\ref{thm:exactness}]
Since $\widetilde{\mathcal P}_\bullet=\bigoplus_{W\in\mathcal J(P)}\widetilde{\mathcal P}_\bullet e_W$, it suffices to show that each $\widetilde{\mathcal P}_\bullet e_W$ is exact. By Lemma~\ref{lem:evaluated-differential}, $\widetilde{\mathcal P}_\bullet e_W$ is isomorphic to the complex $\widetilde{C}_\bullet$ of Lemma~\ref{lem:perturbed-retraction}, which is homotopy equivalent to $\widetilde{C}^{\leq\underline{W}}_\bullet$ by Lemma~\ref{lem:perturbed-retraction}\textup{(2)}, and $\widetilde{C}^{\leq\underline{W}}_\bullet$ is exact by Lemma~\ref{lem:perturbed-lower-exact}.
\end{proof}

\begin{figure}[htbp]
\centering
\begin{tikzpicture}[
  >=stealth,
  every node/.style={inner sep=2pt},
  toE/.style={fill=cyan!25},
  toT/.style={fill=green!30},
  toET/.style={fill=orange!35},
  toZ/.style={draw=black!30},
  sub/.style={draw=red, very thick},
  kept/.style={->},
  gone/.style={->, dashed, black!40},
  hmap/.style={->, magenta, thick},
  labl/.style={fill=white, inner sep=1pt, font=\scriptsize, text=blue!70!black}
]
\fill[fill=yellow!50!white, rounded corners=5pt] (3.6,-6.4) -- (8.2,-7.6) -- (9,-7.6) -- (11.6,-6.2) -- (11.6,-5.8) -- (9,-4.4) -- (8.2,-4.4) -- (3.6,-5.6) -- cycle;
\node[toZ] (a2) at (0,-2.4)  {$\langle\underline{\emptyset}<\underline{1}<\underline{13}\rangle$};
\node[toZ] (b2) at (0,-4.8)  {$\langle\underline{\emptyset}<\underline{3}<\underline{13}\rangle$};
\node[toZ] (c2) at (0,-7.2)  {$\langle\underline{\emptyset}<\underline{3}<\underline{23}\rangle$};
\node[toZ] (d2) at (0,-9.6)  {$\langle\underline{\emptyset}<\underline{2}<\underline{23}\rangle$};
\node[toET] (one-13)  at (4.3,-1.2)  {$\langle\underline{1}<\underline{13}\rangle$};
\node[toZ] (e-1)      at (4.3,-2.4)  {$\langle\underline{\emptyset}<\underline{1}\rangle$};
\node[toET] (e-13)    at (4.3,-3.6)  {$\langle\underline{\emptyset}<\underline{13}\rangle$};
\node[toZ] (three-13) at (4.3,-4.8)  {$\langle\underline{3}<\underline{13}\rangle$};
\node[toET] (e-3) at (4.3,-6.0)  {$\langle\underline{\emptyset}<\underline{3}\rangle$};
\node[toZ] (three-23) at (4.3,-7.2)  {$\langle\underline{3}<\underline{23}\rangle$};
\node[toET] (e-23)    at (4.3,-8.4)  {$\langle\underline{\emptyset}<\underline{23}\rangle$};
\node[toZ] (e-2)      at (4.3,-9.6)  {$\langle\underline{\emptyset}<\underline{2}\rangle$};
\node[toET] (two-23)  at (4.3,-10.8) {$\langle\underline{2}<\underline{23}\rangle$};
\node[toE] (one)         at (8.6,-1.2)  {$\langle\underline{1}\rangle$};
\node[toT] (thirteen)    at (8.6,-2.4)  {$\langle\underline{13}\rangle$};
\node[toE] (e)       at (8.6,-4.8)  {$\langle\underline{\emptyset}\rangle$};
\node[toT] (three)   at (8.6,-7.2)  {$\langle\underline{3}\rangle$};
\node[toT] (twentythree) at (8.6,-9.6)  {$\langle\underline{23}\rangle$};
\node[toE] (two)         at (8.6,-10.8) {$\langle\underline{2}\rangle$};
\node (R) at (11.4,-6.0) {$R$};

\draw[kept] (a2) -- (one-13);
\draw[kept] (a2) -- (e-13);
\draw[kept] (a2) -- (e-1);
\draw[gone] (b2) -- (three-13);
\draw[kept] (b2) -- (e-13);
\draw[kept] (b2) -- (e-3);
\draw[gone] (c2) -- (three-23);
\draw[kept] (c2) -- (e-23);
\draw[kept] (c2) -- (e-3);
\draw[kept] (d2) -- (two-23);
\draw[kept] (d2) -- (e-23);
\draw[kept] (d2) -- (e-2);
\draw[gone] (one-13) -- (thirteen);
\draw[kept] (one-13) -- (one);
\draw[kept] (e-1) -- (one);
\draw[kept] (e-1) -- (e);
\draw[gone] (e-13) -- (thirteen);
\draw[kept] (e-13) -- (e);
\draw[kept] (three-13) -- (thirteen);
\draw[kept] (three-13) -- (three);
\draw[gone] (e-3) -- (three);
\draw[kept, very thick] (e-3) -- (e);
\draw[kept] (three-23) -- (twentythree);
\draw[kept] (three-23) -- (three);
\draw[gone] (e-23) -- (twentythree);
\draw[kept] (e-23) -- (e);
\draw[kept] (e-2) -- (two);
\draw[kept] (e-2) -- (e);
\draw[gone] (two-23) -- (twentythree);
\draw[kept] (two-23) -- (two);
\draw[gone] (one) -- (R);
\draw[kept] (thirteen) -- (R);
\draw[gone] (e) -- (R);
\draw[kept, very thick] (three) -- (R);
\draw[kept] (twentythree) -- (R);
\draw[gone] (two) -- (R);

\draw[hmap] (one) to[bend right=30] (e-1);
\draw[hmap] (thirteen) to[bend right=20] (three-13);
\draw[hmap] (twentythree) to[bend left=20] (three-23);
\draw[hmap] (two) to[bend left=30] (e-2);
\draw[hmap] (one-13) to[bend right=30] (a2);
\draw[hmap] (one-13) to[bend left=10] node[labl,pos=0.6]{$-1$} (b2);
\draw[hmap] (e-13) to[bend left=18] node[labl,pos=0.3]{$-1$} (b2.8);
\draw[hmap] (e-23) to[bend right=18] node[labl,pos=0.3]{$-1$} (c2.-8);
\draw[hmap] (two-23) to[bend right=10] node[labl,pos=0.6]{$-1$} (c2);
\draw[hmap] (two-23) to[bend left=15] (d2);
\end{tikzpicture}
\caption{The perturbed complex $C_\bullet$ and its subcomplex  $C^{\le \underline{W}}_\bullet$ (yellow-highlighted part), along with the non-zero components (in magenta) of the homotopy $h$, and the effect of applying $(r_W)_\#$ (color-coded nodes).}
\label{fig:perturbed-complex}
\end{figure}
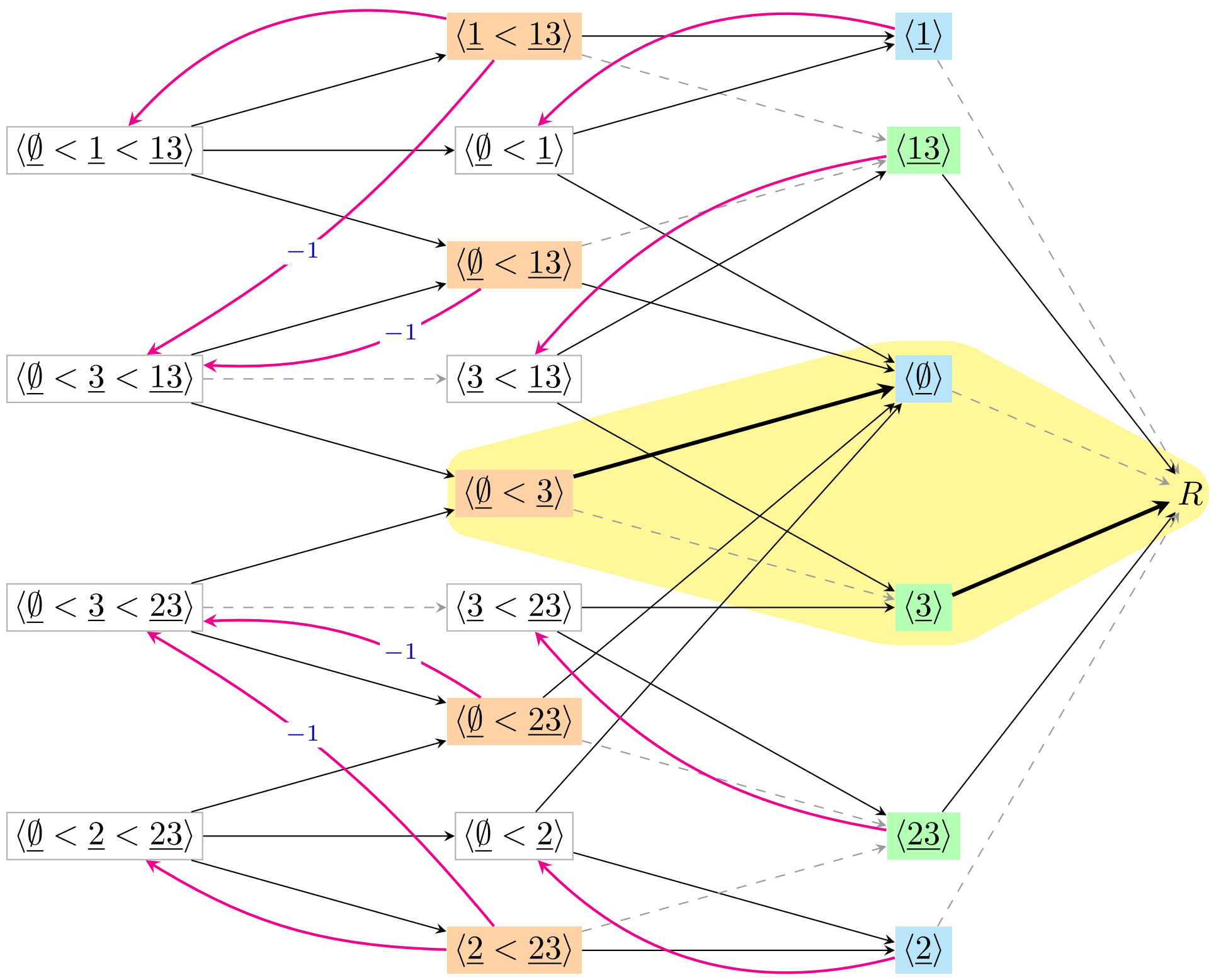

\begin{eg}\label{eg:perturbed-lemmas} We continue with Examples~\ref{subsec:eg1} and~\ref{subsec:eg2}, where $X=[1]$, and take $W=[13]$, so that $\underline{W}=\underline{3}$. Then 
$$
r_W: \begin{array}{ccc} \underline{\emptyset},\underline{1},\underline{2} &\mapsto&  \underline{\emptyset} \\
\underline{3},\underline{13},\underline{23} &\mapsto& \underline{3} \end{array} \quad \text{ and }\quad  \D_X^{\le \underline{W}} = \{\underline{\emptyset} < \underline{3}\}.
$$
Figure~\ref{fig:perturbed-complex} shows the complex $\widetilde{C}_\bullet$ of Lemma~\ref{lem:perturbed-retraction}; it is drawn in the same layout as Figure~\ref{fig:complex-picture} (with the augmentation part $R$ inserted) in order to demonstrate Lemma~\ref{lem:evaluated-differential}.
\begin{itemize}
 \item The forward (left-to-right) arrows correspond to the components of the differentials in $\widetilde{P}_\bullet$.  The solid arrows are the non-zero contributions of the perturbed differential $\partial-\omega$. The dashed arrows are the components of $\partial$ that are killed by the perturbation $\omega$.
 
 \item The yellow-highlighted part is $C_\bullet^{\le \underline{W}}$, which is explicitly of the form:
\[
\begin{tikzcd}[row sep=tiny, column sep=large]
  0 \arrow[r] & R\langle\underline{\emptyset}<\underline{3}\rangle \arrow[r, "\delta_1"] & R\langle\underline{\emptyset}\rangle\oplus R\langle\underline{3}\rangle \arrow[r, "\delta_0"] & R \arrow[r] & 0 \\
  & \langle\underline{\emptyset}<\underline{3}\rangle \arrow[r, mapsto] & -\langle\underline{\emptyset}\rangle & & \\
  & & \langle\underline{\emptyset}\rangle \arrow[r, mapsto] & 0 & \\
  & & \langle\underline{3}\rangle \arrow[r, mapsto] & 1 &
\end{tikzcd}
\]
This is an acyclic complex which, as shown in the proof of Lemma~\ref{lem:perturbed-lower-exact}, decomposes into $K_{\underline{\emptyset},\bullet}\oplus\bigl(K_{\underline{3},0}\xrightarrow{\sim}R\bigr)$ with $K_{\underline{\emptyset}} \cong \widetilde{C}_\bullet(\Delta\{\underline{3}\})[1]$.

  \item The node color encodes its image under $(r_W)_\#$, apart from the degree $-1$ term $R$, which is sent to $0$ by $(r_W)_\#$.  Namely, blue nodes are sent to $R\langle \underline{\emptyset}\rangle$, green nodes to $R\langle \underline{3}\rangle$, orange nodes to $R\langle \underline{\emptyset}< \underline{3})\rangle$, and white nodes to $0$.

  \item The backward arrows correspond to the components of the homotopy $h$; the label is the sign of the corresponding map.
\end{itemize}
\end{eg}

\subsection{Refining projective dimension}
\label{subsec:sharper-dimension-bound}

The cases $|P|\leq2$ are described in Example~\ref{eg:posets-at-most-two-elements}. For $|P|\geq3$, we will improve the bound in Corollary~\ref{cor:distributive-ig} by splitting the last differential of $\mathcal P_\bullet$. We first record the highest possible projective term, then use the topology of the boundary of a simplex to split its differential.

\begin{lem}\label{lem:chain-degree-bound}
Let $X\in\mathcal J(P)$ and suppose that $P\setminus X\notin\mathcal J(P)$. Then
\[
\Delta^{>\underline X}_r=\emptyset
\quad\text{for }r\geq|P|-1,
\qquad
\mathcal P_k=0
\quad\text{for }k\geq|P|.
\]
Moreover, the identification in \eqref{eq:candidate-terms} restricts to an isomorphism
\[
\mathcal P_{|P|-1}
\cong P_X\otimes_R C_{|P|-2}(\Delta^{>\underline X}).
\]
\end{lem}

\begin{proof}
Since $P\setminus X\notin\mathcal J(P)$, we have $P\notin\mathcal D_X$. Every strict chain of non-empty members of $\mathcal D_X$ therefore has at most $|P|-1$ members, i.e. $\Delta^{>\underline X}_{|P|-1}=\emptyset$.   This also implies that every strict chain in $\mathcal D_X$ has at most $|P|$ members, so $\mathcal P_k=0$ for $k\geq|P|$. For the final isomorphism, note that every chain in $\mathcal D_X$ with $|P|$ members begins at $\emptyset=\underline X$, so the claim follows from \eqref{eq:candidate-terms}.
\end{proof}

\begin{eg}\label{eg:chain-degree-bound}
Continue with $P=\{1<2\}\sqcup\{3\}$ and $X=\{1\}$ from Example~\ref{subsec:eg1}. The contribution coming from the complex $\Delta^{>\underline X}$ is the yellow-highlighted part in Figure \ref{fig:complex-picture}.  We have $C_1(\Delta^{>\underline X})\cong R^{\oplus 4}$, and so $\mathcal{P}_2\cong P_X^{\oplus 4}$ and $\mathcal{P}_k=0$ for $k\ge 3$.
\end{eg}

\begin{prop}\label{prop:n-2-inequalities}
Assume that $|P|\geq3$. Then, for every $X\in\mathcal J(P)$,
\[
\pdim_{\Sigma_{\mathcal J(P)}}I_X\leq |P|-2.
\]
If $R$ is Iwanaga--Gorenstein, then $\Sigma_{\mathcal J(P)}$ is Iwanaga--Gorenstein and
\[
\idim(\Sigma_{\mathcal J(P)})_{\Sigma_{\mathcal J(P)}}
=\idim{}_{\Sigma_{\mathcal J(P)}}\Sigma_{\mathcal J(P)}
\leq |P|-2+\idim_R R.
\]
\end{prop}

\begin{proof}
Fix $X\in\mathcal J(P)$. If $P\setminus X\in\mathcal J(P)$, then $I_X$ is projective by Lemma~\ref{lem:projective-relative-injective}, and the projective dimension bound follows. Assume that $P\setminus X\notin\mathcal J(P)$.

\smallskip
\noindent
We first claim that there is an $R$-linear map
\[
\rho:
C_{|P|-3}(\Delta^{>\underline X})
\longrightarrow C_{|P|-2}(\Delta^{>\underline X})
\]
such that $\rho\partial_{|P|-2}=\id$.
By Lemma~\ref{lem:split-boundary-subset-chains}, it suffices to show that $\Delta^{>\underline X}$ is a proper subcomplex of $\Delta(2^P\setminus\{\emptyset,P\})$.  Indeed, since $P\setminus X\notin\mathcal J(P)$, we have $P\notin\mathcal D_X$. 
By the assumption $P\setminus X\notin \mathcal J(P)$ and the fact that $\underline{X}=\emptyset$, every vertex of $\Delta^{>\underline X}$ is a non-empty proper subset of $P$. 
Another consequence of $P\setminus X\notin\mathcal J(P)$ is that there exist $a<b$ in $P$ with $a\in X$ and $b\in P\setminus X$. 
Hence for any $T\in\mathcal J(P)$, the set $\underline{T}=X\SDiff T$ cannot contain both $a$ and $b$. 
In particular,  $\{a,b\}\notin\mathcal D_X$.  Since $|P|\geq3$, the set $\{a,b\}$ is a vertex of $\Delta(2^P\setminus\{\emptyset,P\})$. Thus $\Delta^{>\underline X}$ is a proper subcomplex of $\Delta(2^P\setminus\{\emptyset,P\})$, as required.

Now, we will construct the left inverse of $d_{|P|-1}:\mathcal{P}_{|P|-1}\to \mathcal{P}_{|P|-2}$. 
Consider the projection of $\mathcal{P}_{|P|-2}$ onto the component $P_X\otimes_R C_{|P|-3}(\Delta^{> \underline{X}})$ via  \eqref{eq:candidate-terms}, that is,
\[
\pi:\mathcal P_{|P|-2}
\longrightarrow P_X\otimes_R C_{|P|-3}(\Delta^{>\underline X}), \quad p^{(\underline Z_\ast)}\mapsto \begin{cases}
p\otimes\langle\underline Z_*^{\widehat{0}}\rangle,
& Z_0=X \text{ and } p\in P_{Z_0},\\
0,& \text{otherwise}.
\end{cases}
\]

Consider a chain $\underline{Z}_\ast=(\underline Z_0<\underline Z_1<\cdots<\underline Z_{|P|-1})\in \D_X$ with $Z_0=X$ and let $p\in P_X$. 
Then by definition of the differential \eqref{eq:candidate-differential} we have
\begin{align}
\pi d_{|P|-1}
\bigl(p^{(\underline{Z}_\ast)}\bigr)
&=
\sum_{i=1}^{|P|-1}(-1)^i p\otimes
\langle\underline Z_1<\cdots<\widehat{\underline Z_i}
<\cdots<\underline Z_{|P|-1}\rangle\notag
\\
&=-p\otimes\partial_{|P|-2}
\langle\underline Z_1<\cdots<\underline Z_{|P|-1}\rangle. \label{eq: pi dP-1}
\end{align}
Combining with the isomorphism in Lemma \ref{lem:chain-degree-bound} we obtain the following composition
\[ \mathcal P_{|P|-1} \xrightarrow{d_{|P|-1}} \mathcal P_{|P|-2}
\xrightarrow{\pi}
P_X\otimes_R C_{|P|-3}(\Delta^{>\underline X})\xrightarrow{\ -\id_{P_X}\otimes\rho\ }
P_X\otimes_R C_{|P|-2}(\Delta^{>\underline X})\cong \mathcal P_{|P|-1}.
\]
Since $\rho$ is the left inverse of $\partial_{|P|-2}$, we get from \eqref{eq: pi dP-1} that
\[
(-\id_{P_X}\otimes\rho)\pi d_{|P|-1} = 
(-\id_{P_X}\otimes\rho) (-\id_{P_X}\otimes\partial_{|P|-2}) = \id_{P_X}\otimes \id_{C_{|P|-2}(\Delta^{>\underline X})},
\]
which means that $d_{|P|-1}$ is a split monomorphism. Thus we can truncate the resolution by one more term and get that
$\pdim_{\Sigma_{\mathcal J(P)}}I_X\leq|P|-2$.
The Iwanaga--Gorenstein assertion and the self-injective dimension bound now follow from Corollary~\ref{cor:distance-gorenstein-criterion}.
\end{proof}

\begin{rmk}
The bound in Proposition~\ref{prop:n-2-inequalities} is optimal at least for $|P|=3$: see Example~\ref{eg:morse-matching-change-of-basis} where the bound is attained.
\end{rmk}

\section{Morse reduction and the minimal projective resolution}
\label{sec:minimal-morse}

While we now know that ${\mathcal P}_\bullet$ is a projective resolution of $I_X$, it is far from being minimal.  To obtain a minimal resolution, we apply algebraic Morse theory \cite{CLZ26} to remove its invertible differential components. 
This allows us to describe the term using homologies of $\Delta(\D_X^{>Z})$, which can be used to show (in subsection~\ref{subsec:CM regularity}) a reformulation of $\pdim I_X$ in terms of Castelnuovo--Mumford regularity of certain Stanley--Reisner rings.

Throughout, we require the ground ring $R=\Bbbk$ to be a field.

\subsection{Algebraic Morse theory}
\label{subsec:algebraic-morse}

Algebraic Morse theory~\cite{Sko06,JW09} reduces a complex to a homotopy equivalent complex by matching summands along invertible differential components. We use the formulation of Chen--Liu--Zhou~\cite[Section~3]{CLZ26}.

\begin{dfn}\label{dfn:morse-matching}
Let $(X_\bullet,d)$ be a complex of right $A$-modules with fixed decompositions $X_n=\bigoplus_{\alpha\in I_n}X_{n,\alpha}$. Write $d_{\beta,\alpha}\colon X_{n,\alpha}\to X_{n-1,\beta}$ for the matrix components of $d_n$.
\begin{enumerate}
\item The \defn{weighted quiver} $Q(X_\bullet)$ has the summands $X_{n,\alpha}$ as vertices and an arrow $X_{n,\alpha}\xrightarrow{d_{\beta,\alpha}}X_{n-1,\beta}$ weighted by $d_{\beta,\alpha}$ whenever $d_{\beta,\alpha}\neq0$.

\item A \defn{partial matching} $\mathcal M$ is a collection of arrows, called \defn{matching}, such that
\begin{enumerate}
\item no two selected arrows share a vertex;
\item the weight of every selected arrow is an isomorphism.
\end{enumerate}
Arrows not in $\mathcal M$ are said to be \defn{unmatched}.

\item Form a new weighted quiver $Q^{\mathcal M}$ by
replacing each matched arrow $W\xrightarrow{w}U$ by a \defn{reversed arrow}, denoted by $U\morsedottedarrow W$, with weight $-w^{-1}$.
A vertex not incident with a matched arrow is called \defn{critical}.

\item A directed path in $Q^{\mathcal M}$ is a \defn{zigzag path} if unmatched arrows and reversed arrows alternate; a path consisting of one arrow is also allowed.  We shall use the following single path notation.  For a vertex $U$ in degree $n$ and a vertex $V$ in degree $n-1$, let $\Gamma_n^{\mathcal M}(U,V)$ be the set of zigzag paths from $U$ to $V$ of the form
\[
 U=U_0\morsethickarrow U_1
 \morsedottedarrow U_2
 \morsethickarrow\cdots
 \morsedottedarrow U_{2r}
 \morsethickarrow U_{2r+1}=V,
 \qquad r\geq0.
\]
The weight $\weight(\gamma)$ of a zigzag path $\gamma$ is the composition, along the order of the zigzag path $\gamma$, of the weights of its arrows.

\item The partial matching $\mathcal M$ is a \defn{Morse matching} if it satisfies the following local finiteness hypothesis: for every matched arrow $W\xrightarrow{w}U$, every vertex $V$ in the same degree, say $n-1$, as $U$, and every $x\in U$, the sum
\[
 \sum_{\gamma\in\Gamma_n^{\mathcal M}(W,V)}
 \weight(\gamma)\bigl(-w^{-1}(x)\bigr)
\]
is well-defined, and it is non-zero for only finitely many vertices $V$.
\end{enumerate}
\end{dfn}

One special case where the local finiteness hypothesis  can be conveniently verified is the following.

\begin{prop}\textnormal{(\cite[Proposition~3.2]{CLZ26})}
\label{prop:morse-sufficient-condition}
Let $\mathcal M$ be a partial matching of $Q(X_\bullet)$.  If $Q^{\mathcal M}$ contains no infinite zigzag path, then $\mathcal M$ is a Morse matching.
\end{prop}

We next record the main result of algebraic Morse theory.

\begin{thm}\textnormal{(\cite[Theorem~3.3]{CLZ26})}
\label{thm:algebraic-morse}
Let $\mathcal M$ be a Morse matching of $Q(X_\bullet)$, and write $\Crit_n(\mathcal M)$ for its critical vertices in degree $n$.  Then we have a complex $(X_\bullet^{\mathcal{M}},d^{\mathcal{M}})$, called the \defn{Morse complex}, given by
\[
 X_n^{\mathcal M}
 :=\bigoplus_{U\in\Crit_n(\mathcal M)}U,  \qquad d_{V,U}^{\mathcal M}
 :=\sum_{\gamma\in\Gamma_n^{\mathcal M}(U,V)} \weight(\gamma),
\]
where $d^{\mathcal{M}}_{V,U}$ is the $(V,U)$-component of $d_n^{\mathcal{M}}:X_n^{\mathcal{M}}\to X_{n-1}^{\mathcal{M}}$ for $U\in\Crit_n(\mathcal M)$ and $V\in\Crit_{n-1}(\mathcal M)$.
Moreover, it is homotopy equivalent to $(X_\bullet,d)$ -- more precisely, one has a decomposition of chain complexes
\[
 (X_\bullet,d)\cong
 (X_\bullet^{\mathcal M},d^{\mathcal M})
 \oplus(Y_\bullet,d^Y),
\]
where $Y_\bullet$ is null-homotopic.
\end{thm}

In particular, if $X_\bullet$ is a projective resolution of a module, then its Morse complex is also a projective resolution of the same module.

\subsection{Construction of the Morse matching}
\label{subsec:construction-morse-matching}

Let us first explain the idea of how to find a Morse matching. The differential formula~\eqref{eq:candidate-differential} separates the faces of a strict chain into two types: the zeroth face deletes its initial member, whereas the remaining faces preserve the initial member.  Under the identification~\eqref{eq:candidate-terms}, these two contributions give a decomposition $d=\eta+\partial$, where the non-zero components of $\partial$ are signed identity maps between copies of the same projective, while those of $\eta$ are induced by the left-multiplication map (a more precise formulation will come after the next paragraph); a related separation is used in Lemma~\ref{lem:evaluated-differential}.

We want to use a Morse matching to eliminate the identity components in the differential, but several of them may share a source or target (see Example \ref{subsec:eg2}), which obstructs any naive choice of components to form a Morse matching. 
The `solution' is to use homological perturbation theory~\cite{GL89}: decompose using a choice of kernel--image splitting for each chain complex $C_\bullet(\Delta^{>\underline{Z}})$. 
With respect to this new choice of decomposition, $\partial$ will naturally define a collection of isomorphisms with pairwise disjoint endpoints, whereas $\eta$ only contributes radical maps.  This allows us to obtain a Morse matching that reduces $\mathcal{P}_\bullet$ to a minimal resolution.

\medskip

We now construct the matching. All tensor products below are over $\Bbbk$. For $Z\in\mathcal J(P)$, put
\[
 C_{Z,\bullet}:=\widetilde C_\bullet(\Delta^{>\underline{Z}})[1],
 \qquad C_{Z,k}=\widetilde C_{k-1}(\Delta^{>\underline{Z}}),
 \qquad \partial^Z_k:C_{Z,k}\longrightarrow C_{Z,k-1},
\]
where $\partial^Z_k = -\partial_{k-1}$ is the negative of the boundary map on $\widetilde{C}_\bullet(\Delta^{>\underline{Z}})$. Thus $C_{Z,0}=\Bbbk 1_Z$, $C_{Z,k}=0$ for $k<0$, and $\partial^Z_0=0$.

For $\underline{Z}<\underline{T}$ and $k\geq1$, define $\eta_{T,Z,k}:C_{Z,k}\to C_{T,k-1}$ on the standard basis by
\[
 \eta_{T,Z,k}\langle \underline{Z}_1<\cdots<\underline{Z}_k\rangle
 =
 \begin{cases}
  1_T,&T=Z_1\text{ and }k=1,\\
  \langle \underline{Z}_2<\cdots<\underline{Z}_k\rangle,
       &T=Z_1\text{ and }k>1,\\
  0,&T\neq Z_1.
 \end{cases}
\]
Note that here we have $\underline{Z}<\underline{Z}_1<\cdots<\underline{Z}_k$; when $T=Z_1$, the output is a basis element of $C_{T,k-1}$. 
Under \eqref{eq:candidate-terms},  we have
\[
 \begin{aligned}
 d_k&=\partial_k+\eta_k:
       \mathcal P_k\longrightarrow\mathcal P_{k-1}, \quad \text{ for all $k\ge 1$,}\\
 \text{ with } \partial_k&:=\bigoplus_Z\id_{P_Z}\otimes\partial^Z_k,
 \qquad
 \eta_k:=\sum_{\underline{Z}<\underline{T}}\mu_{T,Z}\otimes\eta_{T,Z,k}.
 \end{aligned}
\]
The $(T,Z)$-component of $\eta_k$ has source $P_Z\otimes C_{Z,k}$ and target $P_T\otimes C_{T,k-1}$. Recall that $\mu_{T,Z}:P_Z\to P_T$ is given by $p\mapsto b_{TZ}p$.


Since $\Bbbk$ is a field, choose, for every $Z$ and $k\geq0$, a decomposition
\begin{equation}
 C_{Z,k}=B_{Z,k}\oplus H_{Z,k}\oplus L_{Z,k}\;
\text{ for all $k\ge 0$, where } 
\begin{cases}
    B_{Z,k}:=\Image\partial^Z_{k+1},\\
    B_{Z,k}\oplus H_{Z,k}=\Ker\partial^Z_k,\\
    \left.\partial^Z_k\right|_{L_{Z,k}}:
 L_{Z,k}\xrightarrow{\sim}B_{Z,k-1}.
 \end{cases}
 \label{eq:homology-splitting}
\end{equation}
In particular, $L_{Z,0}=0$, and the canonical projection from $\Ker \partial_k^Z$ restricts to an isomorphism
\[
 H_{Z,k}\xrightarrow{\sim}H_k(C_{Z,\bullet})
 \cong\widetilde H_{k-1}(\Delta^{>\underline{Z}}),\qquad h\longmapsto h+B_{Z,k}.
\]
Use the resulting decomposition
\[
 \mathcal P_k=\bigoplus_{Z\in\mathcal J(P)}
 \bigl((P_Z\otimes B_{Z,k})\oplus
       (P_Z\otimes H_{Z,k})\oplus
       (P_Z\otimes L_{Z,k})\bigr)
\]
to form the weighted quiver of the resolution part $\mathcal P_\bullet$, that is, the vertices of the weighted quiver at degree $k$ are given by $(Z,E)$ for $Z\in\mathcal{J}(P)$ and $E\in\{B,H,L\}$ such that the summand $P_Z\otimes E_{Z,k}$ is non-zero.  Note that the map $\partial_k$ sends $P_Z\otimes C_{Z,k}$ into $P_Z\otimes C_{Z,k-1}$, whereas $\eta_k$ sends it into $\bigoplus_{\underline{Z}<\underline{T}}P_T\otimes C_{T,k-1}$.

\begin{thm}\label{thm:minimal-resolution}
Suppose the ground ring is a field $R=\Bbbk$. With respect to the decomposition \eqref{eq:homology-splitting}, define $\mathcal M$ to be the collection of all non-zero components of 
\begin{equation}
 \id_{P_Z}\otimes\left.\partial^Z_k\right|_{L_{Z,k}}:
 P_Z\otimes L_{Z,k}\xrightarrow{\sim}P_Z\otimes B_{Z,k-1}
 \label{eq:morse-matching}
\end{equation}
for all $Z\in \mathcal{J}(P)$ and all $k\ge 1$.
Then $\mathcal M$ is a Morse matching, and its associated Morse complex defines a minimal projective resolution
\[
 \cdots\longrightarrow\mathcal P_{\min,2}
 \longrightarrow\mathcal P_{\min,1}
 \longrightarrow\mathcal P_{\min,0}
 \longrightarrow I_X\longrightarrow0
\]
whose terms, for $k\geq0$, are
\begin{equation}
 \begin{aligned}
 \mathcal P_{\min,k}
 &=\bigoplus_{Z\in\mathcal J(P)}P_Z\otimes H_{Z,k} \cong\bigoplus_{Z\in\mathcal J(P)}
 P_Z\otimes\widetilde H_{k-1}(\Delta^{>\underline{Z}})
 \cong\bigoplus_{Z\in\mathcal J(P)}
 P_Z^{\oplus\widetilde\beta_{Z,k}},
 \end{aligned}
 \label{eq:minimal-terms}
\end{equation}
where $\widetilde\beta_{Z,k}:= \dim_{\Bbbk}\widetilde H_{k-1}(\Delta^{>\underline{Z}})$.
\end{thm}

\begin{proof} 
To see that $\mathcal M$ is a partial matching, notice first that both the source $P_Z\otimes L_{Z,k}$ and the target $P_Z\otimes B_{Z,k-1}$ of each displayed map have label $Z$. Since every non-zero component of $\eta_k$ changes the label from $Z$ to some $T$ with $\underline{Z}<\underline{T}$, the component of $\eta_k$ between these two summands is zero. Hence the component of $d_k=\partial_k+\eta_k$ from $P_Z\otimes L_{Z,k}$ to $P_Z\otimes B_{Z,k-1}$ is exactly $\id_{P_Z}\otimes\left.\partial^Z_k\right|_{L_{Z,k}}$. Their endpoints are pairwise disjoint and their weights are isomorphisms, so $\mathcal M$ is a partial matching.

To show that $\mathcal M$ is a Morse matching, by Proposition~\ref{prop:morse-sufficient-condition}, it suffices to show that $\mathcal M$ has no infinite zigzag path.
Indeed, since every non-zero component of $\partial$ belongs to $\mathcal M$, every unmatched arrow in $Q^{\mathcal M}$ is a component of $\eta$. Thus, for $k\geq1$, every unmatched arrow has the form
\[
 P_Z\otimes E_{Z,k}
 \xrightarrow{\mu_{T,Z}\otimes\eta'_{T,Z}}
 P_T\otimes F_{T,k-1}
 \quad \text{ with }\underline{Z}<\underline{T}
\]
for appropriate $E,F\in\{B,H,L\}$, and the map $\eta'_{T,Z}$ is given by the composition
\begin{equation*}\label{eq:delta'}
 \eta'_{T,Z}:
 E_{Z,k}\hookrightarrow C_{Z,k}
 \xrightarrow{\eta_{T,Z,k}}C_{T,k-1}
 \twoheadrightarrow F_{T,k-1}.
\end{equation*}
On the other hand, every reversed arrow goes from $P_Z\otimes B_{Z,k-1}$ to $P_Z\otimes L_{Z,k}$. Hence, along a zigzag path, reversed arrows preserve the label $\underline{Z}\in\mathcal D_X$, whereas unmatched arrows strictly increase it. 
Since $\mathcal D_X$ is finite and the two types of arrows alternate, zigzag paths have bounded length. In particular, there is no infinite zigzag path, as required.

Theorem~\ref{thm:algebraic-morse} yields a Morse complex $(\mathcal P_{\min,\bullet},d_{\min})$ homotopy equivalent to $\mathcal P_\bullet$.
Its terms are projective, and $\mathcal P_\bullet$ resolves $I_X$ by Theorem~\ref{thm:exactness}.  Hence $\mathcal P_{\min,\bullet}$ is a projective resolution of $I_X$.  To see that the individual terms are as described in \eqref{eq:minimal-terms}, recall from our choice of decomposition \eqref{eq:homology-splitting} and $\mathcal{M}$ that every non-zero summand $P_Z\otimes B_{Z,k}$ is matched with $P_Z\otimes L_{Z,k+1}$.
Thus the critical summands in degree $k$ are precisely the non-zero $P_Z\otimes H_{Z,k}$. This, together with $H_{Z,k}\cong\widetilde H_{k-1}(\Delta^{>\underline{Z}})$, yields the required description \eqref{eq:minimal-terms}. 

It remains to prove minimality. Note that every unmatched arrow has weight $\mu_{T,Z}\otimes\eta'_{T,Z}:P_Z\otimes E_{Z,k}\to P_T\otimes F_{T,k-1}$, where $k\geq1$ and $\underline{Z}<\underline{T}$, with appropriate $E,F\in\{B,H,L\}$.
Since $Z\neq T$, we have $b_{TZ}\in e_T\rad\Sigma_{\mathcal J(P)}$, which implies that the map $\mu_{T,Z}$ given by left-multiplying $b_{TZ}$ is radical.  Hence, we have
\begin{equation}\label{eq:image muTZ otimes delta'}    
 \Image(\mu_{T,Z}\otimes\eta'_{T,Z})
 \subseteq (P_T\rad\Sigma_{\mathcal J(P)})\otimes F_{T,k-1}
 = \rad P_T\otimes F_{T,k-1}.
\end{equation}
Fix $k\geq1$ and a zigzag path $\gamma$ from $P_Z\otimes H_{Z,k}$ to $P_T\otimes H_{T,k-1}$ contributing to $d_{\min,k}$. Since the target of $\gamma$ is critical, the last arrow of $\gamma$ has to be unmatched. 
The image of $\weight(\gamma)$ is contained in the image of this last arrow, so the inclusion \eqref{eq:image muTZ otimes delta'} applied to $F=H$ gives $\Image\weight(\gamma)\subseteq\rad P_T\otimes H_{T,k-1}$. By Theorem~\ref{thm:algebraic-morse}, each matrix component of $d_{\min,k}$ is a sum of these path weights; the sum is finite because $Q^{\mathcal M}$ is finite and its zigzag paths have bounded length. Hence $\Image d_{\min,k}\subseteq\rad\mathcal P_{\min,k-1}$ for every $k\geq1$, which proves minimality.
\end{proof}

The choices in \eqref{eq:homology-splitting} may change the individual differentials in the Morse complex.  Nevertheless, the resulting complexes are minimal projective resolutions of the same module and are therefore isomorphic.  In particular, $P_Z$ occurs in degree $k$ precisely when $\widetilde H_{k-1}(\Delta^{>\underline{Z}})\neq0$.

\begin{rmk}
Note that the multiplicity of $P_Z$ appearing in the minimal projective resolution of $I_X$ is given by $\dim_\Bbbk\Ext_A^k(I_X,S_Z)$, and so we have 
\[\widetilde{\beta}_{Z,k}=\dim_{\Bbbk}\widetilde H_{k-1}(\Delta^{>\underline{Z}})=\dim_{\Bbbk}\Ext_A^k(I_X,S_Z).
 \]
\end{rmk}

We can now refine Corollary \ref{cor:distributive-ig} when working over a field.
\begin{cor}\label{cor:SigmaJ(P) is IG} Assume that the ground ring is a field $R=\Bbbk$ and $P\neq\emptyset$.
For every $X\in\mathcal J(P)$, the indecomposable injective module $I_X$ over the distance $\Bbbk$-algebra $\Sigma_L$ satisfies
\begin{align}
 \pdim_{\Sigma_L} I_X
 &=\max\left\{k\geq0\ \middle|\
 \widetilde H_{k-1}(\Delta^{>\underline{Z}})\neq0
 \text{ for some }Z\in\mathcal J(P)\right\}.
 \label{eq:projective-dimension}
\end{align}
In particular, for every $Z\in\mathcal J(P)$, $ \widetilde\beta_{Z,k} =\dim_{\Bbbk}\widetilde H_{k-1}(\Delta^{>\underline{Z}})=0$ if any of the following is satisfied
\begin{enumerate}
    \item $k\ge |P|$.
    \item $k=|P|-1$ and $|P|\ge 3$. 
\end{enumerate}
\end{cor}
\begin{proof}
The projective dimension formula follows from the minimal terms in \eqref{eq:minimal-terms}. 
The Betti number vanishing follows from combining Corollary~\ref{cor:distributive-ig} and Proposition~\ref{prop:n-2-inequalities} with \eqref{eq:projective-dimension}.
\end{proof}

\begin{eg}\label{eg:morse-matching-change-of-basis}
Continue with Examples~\ref{subsec:eg1} and~\ref{subsec:eg2}. We can regard Figure~\ref{fig:complex-picture} as the weighted quiver with respect to the standard decomposition of $\mathcal{P}_k$ into indecomposable direct summands, which are indexed by chains in $\Delta(\D_X)$. Let us demonstrate the alternative choice used in Theorem \ref{thm:minimal-resolution}.

For $Z=[1]$, take
\begin{align*}
 L_{[1],2}=C_{[1],2}
  &=\Bbbk\{
  \langle\underline{1}<\underline{13}\rangle,
  \langle\underline{3}<\underline{13}\rangle,
  \langle\underline{3}<\underline{23}\rangle,
  \langle\underline{2}<\underline{23}\rangle\},\\
 B_{[1],1}&=\Bbbk\{
  \langle\underline{1}\rangle-\langle\underline{13}\rangle,
  \langle\underline{3}\rangle-\langle\underline{13}\rangle,
  \langle\underline{3}\rangle-\langle\underline{23}\rangle,
  \langle\underline{2}\rangle-\langle\underline{23}\rangle\},\\
 L_{[1],1}&=\Bbbk\{\langle\underline{13}\rangle\}.
\end{align*}
For the remaining sets with non-zero degree-one terms, take
\begin{align*}
 L_{\emptyset,1}&=\Bbbk\{\langle\underline{13}\rangle\}, \qquad L_{[12],1}=\Bbbk\{\langle\underline{23}\rangle\}, \qquad 
 L_{[13],1}=\Bbbk\{\langle\underline{13}\rangle\}, \qquad 
 H_{[13],1}=
 \Bbbk\{\langle\underline{13}\rangle-
              \langle\underline{23}\rangle\},\\
 B_{Z,0}&=C_{Z,0} \;\;\text{ for } Z\in\{[1],\emptyset,[12],[13]\},&&\\
 H_{Z,0}&=C_{Z,0} \;\;\text{ for } Z\in\{[3],[123]\}.
 &&
\end{align*}
All unlisted spaces in the decomposition \eqref{eq:homology-splitting} are zero.

The map $\partial^{Z=[1]}_2$ is the negative of the simplicial boundary, so on $L_{[1],2}=C_{[1],2}$ it sends each basis element $\langle\underline a<\underline c\rangle$ to $\langle\underline a\rangle-\langle\underline c\rangle$.  In particular it maps the four generators of $L_{[1],2}$ displayed above, in the given order, onto the four generators of $B_{[1],1}$.
Moreover, $\partial^Z_1$ sends the displayed generator of each $L_{Z,1}$ to $-1\in B_{Z,0}$.  These eight isomorphisms have pairwise disjoint endpoints and are selected for the Morse matching.  In terms of the colors in Figure~\ref{fig:complex-picture}, the yellow part gives four matching components from degree two to degree one and one from degree one to degree zero, while the blue, green, and orange parts each give one matching component from degree one to degree zero. In particular, the two green arrows have the same target and both have weight $-1$. After a change of basis, the green part therefore contributes one matching component from $P_{[13]}\otimes L_{[13],1}$ to $P_{[13]}\otimes B_{[13],0}$,  while the difference $\langle\underline{13}\rangle-\langle\underline{23}\rangle$ spans $H_{[13],1}$ and determines the critical summand $P_{[13]}\otimes H_{[13],1}$.

The matrix components of $\eta$ also change.  For example, $ \eta_2\bigl( p\otimes\langle\underline{3}<\underline{23}\rangle \bigr) =b_{[13][1]}p\otimes\langle\underline{23}\rangle$,  and, in the chosen decomposition of $C_{[13],1}$, we get that $\langle\underline{23}\rangle =\langle\underline{13}\rangle -\bigl( \langle\underline{13}\rangle- \langle\underline{23}\rangle \bigr)$. It thus produces two matrix components, of weights $b_{[13][1]}$ and $-b_{[13][1]}$, respectively.  Similarly,
\[
 \eta_1\bigl(
 p\otimes(
 \langle\underline{3}\rangle-\langle\underline{13}\rangle)
 \bigr)
 =b_{[13][1]}p\otimes1-b_{[3][1]}p\otimes1.
\]
These are precisely the corresponding arrows displayed in Figure~\ref{fig:morse-matching-change-of-basis}, where we write $P_Z\langle c\rangle:=P_Z\otimes_{\Bbbk}\Bbbk c$ for compactness.

\begin{figure}[htbp]
\centering
\begin{tikzpicture}[
  >=stealth,
  every node/.style={inner sep=2pt},
  summand/.style={draw,rounded corners,fill=white},
  critical/.style={draw,rounded corners,fill=yellow!40},
  match/.style={->,very thick,draw=blue!75!black},
  other/.style={->,draw=black!45},
  labl/.style={fill=white,inner sep=1pt,font=\small,text=black}
]

\node[draw=none,font=\normalsize\bfseries] at (0,0.9) {$\mathcal P_2$};
\node[draw=none,font=\normalsize\bfseries] at (5,0.9) {$\mathcal P_1$};
\node[draw=none,font=\normalsize\bfseries] at (12.5,0.9) {$\mathcal P_0$};

\node[summand] (E1) at (0,-1.2)  {$P_{[1]}\langle\underline{1}<\underline{13}\rangle$};
\node[summand] (E2) at (0,-5.6)  {$P_{[1]}\langle\underline{3}<\underline{13}\rangle$};
\node[summand] (E3) at (0,-10.0) {$P_{[1]}\langle\underline{3}<\underline{23}\rangle$};
\node[summand] (E4) at (0,-14.4) {$P_{[1]}\langle\underline{2}<\underline{23}\rangle$};

\node[summand]  (Le)  at (5, 0.0)  {$P_{\emptyset}\langle\underline{13}\rangle$};
\node[summand]  (B1)  at (5,-1.8)  {$P_{[1]}\langle\underline{1}-\underline{13}\rangle$};
\node[summand]  (L1)  at (5,-3.6)  {$P_{[1]}\langle\underline{13}\rangle$};
\node[summand]  (B2)  at (5,-5.4)  {$P_{[1]}\langle\underline{3}-\underline{13}\rangle$};
\node[summand]  (L13) at (5,-7.2)  {$P_{[13]}\langle\underline{13}\rangle$};
\node[critical] (H13) at (5,-9.0)  {$P_{[13]}\langle\underline{13}-\underline{23}\rangle$};
\node[summand]  (B3)  at (5,-10.8) {$P_{[1]}\langle\underline{3}-\underline{23}\rangle$};
\node[summand]  (L12) at (5,-12.6) {$P_{[12]}\langle\underline{23}\rangle$};
\node[summand]  (B4)  at (5,-14.4) {$P_{[1]}\langle\underline{2}-\underline{23}\rangle$};

\node[summand]  (Ue)   at (12.5,-0.5)  {$P_{\emptyset}$};
\node[critical] (U3)   at (12.5,-4.0)  {$P_{[3]}$};
\node[summand]  (U1)   at (12.5,-6.0)  {$P_{[1]}$};
\node[summand]  (U13)  at (12.5,-8.0)  {$P_{[13]}$};
\node[critical] (U123) at (12.5,-11.5) {$P_{[123]}$};
\node[summand]  (U12)  at (12.5,-14.5) {$P_{[12]}$};

\draw[match] (E1) -- node[labl] {$1$} (B1);
\draw[match] (E2) -- node[labl] {$1$} (B2);
\draw[match] (E3) -- node[labl] {$1$} (B3);
\draw[match] (E4) -- node[labl] {$1$} (B4);

\draw[other] (E1) -- node[labl,pos=0.48] {$b_{\emptyset[1]}$} (Le);
\draw[other] (E2) -- node[labl,pos=0.48] {$b_{[13][1]}$} (L13);
\draw[other] (E3) to[bend left=10]
  node[labl,pos=0.35] {$b_{[13][1]}$} (L13);
\draw[other] (E3) -- node[labl,pos=0.65] {$-b_{[13][1]}$} (H13);
\draw[other] (E4) -- node[labl,pos=0.48] {$b_{[12][1]}$} (L12);

\draw[match] (Le)  -- node[labl] {$-1$} (Ue);
\draw[match] (L1)  -- node[labl,pos=0.72,above] {$-1$} (U1);
\draw[match] (L13) -- node[labl,pos=0.72,above] {$-1$} (U13);
\draw[match] (L12) -- node[labl,pos=0.72,above] {$-1$} (U12);

\draw[other] (B1) to[out=0,in=-165]
  node[labl,pos=0.48,below] {$b_{\emptyset[1]}$} (Ue);
\draw[other] (B1) to[out=0,in=155]
  node[labl,pos=0.28,below] {$-b_{[3][1]}$} (U3);

\draw[other] (B2) to[out=0,in=165]
  node[labl,pos=0.35,above] {$b_{[13][1]}$} (U13);
\draw[other] (B2) to[out=0,in=-175]
  node[labl,pos=0.30,below] {$-b_{[3][1]}$} (U3);

\draw[other] (B3) to[out=0,in=-165]
  node[labl,pos=0.65,below] {$b_{[13][1]}$} (U13);
\draw[other] (B3) to[out=0,in=170]
  node[labl,pos=0.40,below] {$-b_{[123][1]}$} (U123);

\draw[other] (B4) to[out=0,in=-165]
  node[labl,pos=0.45,below] {$b_{[12][1]}$} (U12);
\draw[other] (B4) to[out=0,in=-155]
  node[labl,pos=0.74,below] {$-b_{[123][1]}$} (U123);

\draw[other] (Le)  to[out=0,in=145]
  node[labl,pos=0.18,above] {$b_{[3]\emptyset}$} (U3);
\draw[other] (L1)  to[out=0,in=175]
  node[labl,pos=0.38,above] {$b_{[3][1]}$} (U3);
\draw[other] (L13) to[out=0,in=-165]
  node[labl,pos=0.55,above] {$b_{[3][13]}$} (U3);
\draw[other] (L12) to[out=0,in=-170]
  node[labl,pos=0.57,above] {$b_{[123][12]}$} (U123);

\draw[other] (H13) to[out=0,in=-145]
  node[labl,pos=0.36,below] {$b_{[3][13]}$} (U3);
\draw[other] (H13) to[out=0,in=155]
  node[labl,pos=0.25,above] {$-b_{[123][13]}$} (U123);

\end{tikzpicture}
\caption{The weighted quiver after the $B\oplus H\oplus L$ change
of basis.}
\label{fig:morse-matching-change-of-basis}
\end{figure}
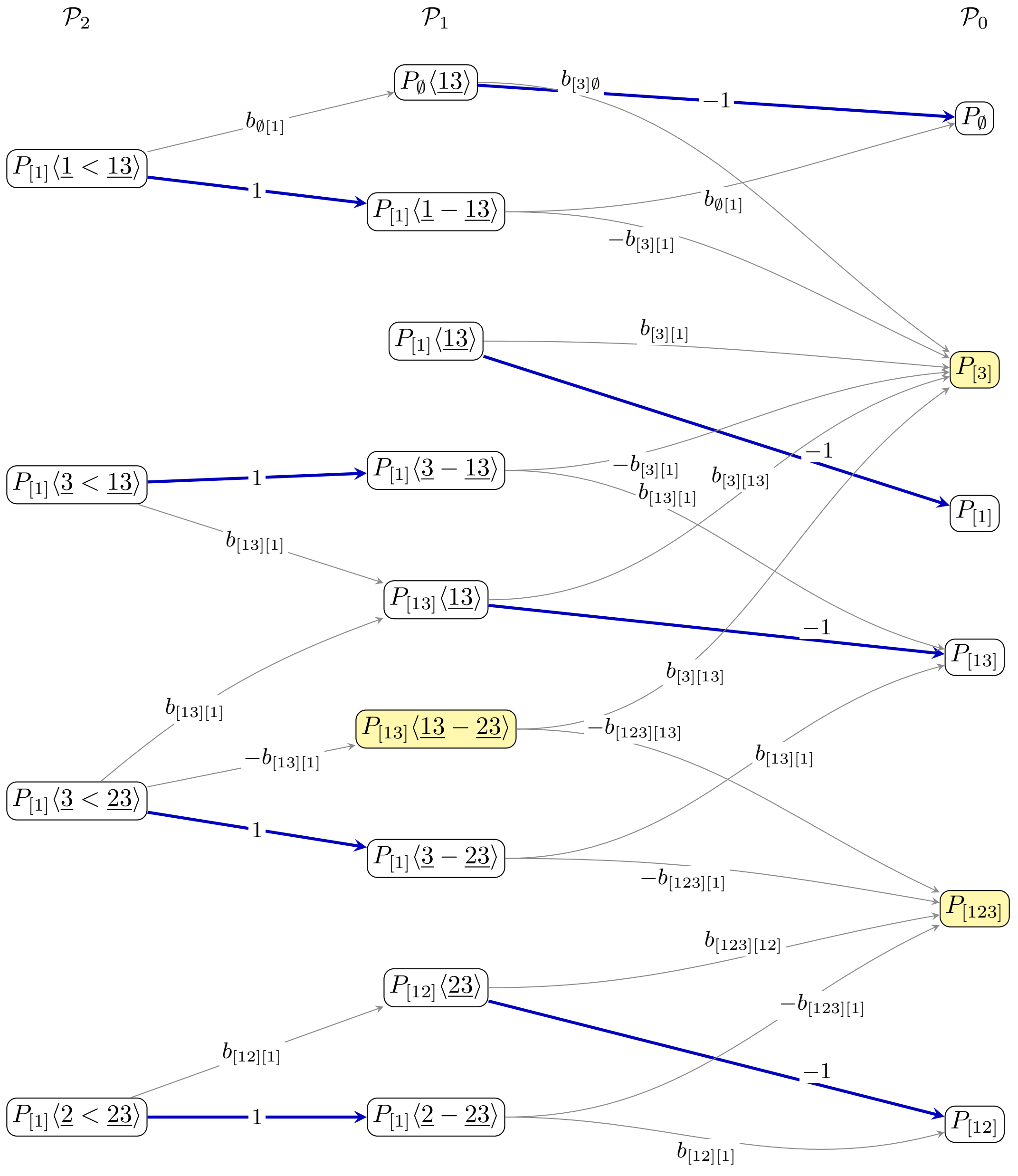

In Figure~\ref{fig:morse-matching-change-of-basis}, the thick blue arrows form the Morse matching, the thin gray arrows are the remaining non-zero matrix components of the differential, and the yellow vertices are critical points.

The two summands $P_{[3]}$ and $P_{[123]}$ without any highlight in Figure~\ref{fig:complex-picture} are precisely the degree-zero terms of the resulting Morse complex.  Its degree-one term is the copy of $P_{[13]}$ represented in the original basis by $p^{(\underline{3}<\underline{13})}-p^{(\underline{3}<\underline{23})}$, for $p\in P_{[13]}$. For the displayed choice of splitting, the remaining differential is
\[
 d_{\min,1}\colon P_{[13]}\longrightarrow P_{[3]}\oplus P_{[123]},
 \qquad
 p\longmapsto
 \bigl(b_{[3][13]}p,-b_{[123][13]}p\bigr).
\]
Note that this example shows that the bound $|P|-2=1$ in Proposition~\ref{prop:n-2-inequalities} is attained.
\end{eg}

\subsection{Projective dimension and Castelnuovo--Mumford regularity}\label{subsec:CM regularity}
The formula for $\pdim I_X$ looks very close to the following characterization of Castelnuovo--Mumford regularity via reduced homology:
\begin{align}
\pdim I_X
&=\max\{k\geq 0:\widetilde H_{k-1}(\Delta^{>\underline{Z}})\neq 0 \text{ for some }Z\in \mathcal{J}(P)\}, \\
\operatorname{reg}\Bbbk[\Delta]
&=\max\{k\geq 0:\widetilde H_{k-1}(\operatorname{lk}_\Delta F)\neq 0 \text{ for some face }F\in \Delta\}. \label{eq:reg=non-vanishing homology}
\end{align}
Here $\Bbbk[\Delta]$ is the Stanley--Reisner ring of a finite simplicial complex $\Delta$, and $\operatorname{reg}$ denotes the Castelnuovo--Mumford regularity (see Definition \ref{def:SR ring} and Lemma \ref{lem:Hochshster easier}). This motivates us to seek a complex $\Delta$ for which the two equations coincide; it turns out that such a $\Delta$ is given by the simplicial complex $\mathcal{S}_X$ `generated by $\mathcal{D}_X$', namely,\[
\mathcal S_X:=\{U\subseteq P \mid U\subseteq \underline{Y} \text{ for some }\underline{Y}\in\mathcal D_X\} = \bigcup_{\underline{Y}\in \max \D_X} 2^{\underline{Y}}.
\]

\begin{thm}\label{thm:pdim=reg}
Let $\Bbbk$ be a field and let $\Sigma_{\mathcal{J}(P)}$ be the distance $\Bbbk$-algebra associated to the lattice $\mathcal{J}(P)$. Then, for any $X\in\mathcal{J}(P)$, we have
\[
\operatorname{pdim}_{\Sigma_{\mathcal{J}(P)}} I_X=\operatorname{reg}\Bbbk[\mathcal S_X].
\]
\end{thm}

In fact, the first step of the proof is to show that $\widetilde H_r(\Delta^{>\underline{Z}})\cong\widetilde H_r(\operatorname{lk}_{\mathcal S_X}\underline{Z})$ for all $r\ge -1$, which gives us $\pdim I_X\le \operatorname{reg} \Bbbk[\mathcal{S}_X]$. The second step is to show that $\operatorname{lk}_{\mathcal S_X}F$ has zero reduced homology whenever $F\in\mathcal S_X\setminus\mathcal D_X$. Both steps use the fact that $\mathcal D_X$ is closed under intersection.

We use the topological preliminaries from subsection~\ref{subsec:topological-tools} to compare $\D_X$ and $\mathcal S_X$.

\begin{lem}\label{lem:face closure map}
The poset $\mathcal D_X$ is closed under intersections of non-empty finite families. The map
\[
c:\mathcal F(\mathcal S_X)\longrightarrow\mathcal D_X,\qquad
F\longmapsto\bigcap_{\substack{\underline{Y}\in\max\mathcal D_X\\F\subseteq \underline{Y}}}\underline{Y}
\]
is well-defined, preserves inclusion, and satisfies $F\subseteq c(F)$.
\end{lem}

\begin{proof}
For $Y,Z\in\mathcal{J}(P)$, we have an order ideal $m(X,Y,Z)$ so that $\underline{m(X,Y,Z)} = \underline{Y}\cap \underline{Z}$ by Lemma~\ref{lem:symmetric-difference-identities}\textup{(\ref{itm:symdiff-intersection})}. Hence, closure under intersections of non-empty finite families follows by induction.

Every $F\in\mathcal S_X$ is contained in a facet, so the family defining $c(F)$ is non-empty. Closure under intersections gives $c(F)\in\mathcal D_X$, and every set in the intersection contains $F$. If $F\subseteq G$, the family of facets containing $G$ is contained in the family containing $F$. Therefore $c(F)\subseteq c(G)$.
\end{proof}

\begin{lem}\label{lem:H(Delta Z)=H(lk DZ)}
For every $Z\in\mathcal J(P)$ and every $k\geq-1$,
\[
\widetilde H_k(\Delta^{>\underline{Z}})\cong
\widetilde H_k\bigl(\operatorname{lk}_{\mathcal S_X}(\underline{Z})\bigr).
\]
\end{lem}

\begin{proof}
The set $\underline{Z}$ is a face of $\mathcal S_X$. Consider the order-preserving maps
\begin{align*}
j_Z:\quad \D_X^{>\underline{Z}} &\longrightarrow\mathcal F(\mathcal S_X)_{>\underline{Z}},
&\underline{T}&\longmapsto \underline{T},\\
c_Z:\mathcal F(\mathcal S_X)_{>\underline{Z}}&\longrightarrow\D_X^{>\underline{Z}},
&F&\longmapsto c(F).
\end{align*}
The map $c_Z$ is well-defined because $\underline{Z}<F\subseteq c(F)$ by Lemma \ref{lem:face closure map}. Moreover, we have
\[
\underline{T}\subset c_Zj_Z(\underline{T}), \text{ for all $\underline{T}\in \D_X^{>\underline{Z}}$} \quad \text{and} \quad 
F\subseteq  j_Zc_Z \text{ for all $F\in\mathcal{F}(\mathcal{S}_X)_{>\underline{Z}}$.}
\]
Hence, Lemma \ref{lem:order homotopy} says that $j_{Z,\#}$ and $c_{Z,\#}$ are chain homotopy inverses. Thus, the order complex of $\D_X^{>\underline{Z}})$ and the order complex of $\mathcal F(\mathcal S_X)_{>\underline{Z}}$ have isomorphic reduced homology. Now we apply Lemma \ref{lem:munkres} and get that
\[
\widetilde{H}_k(\Delta^{>\underline{Z}}) 
\cong \widetilde{H}_k(\Delta\big(\mathcal{F}(\mathcal{S}_X)_{>\underline{Z}}\big)) 
\cong \widetilde{H}_k(\operatorname{lk}_{\mathcal{S}_X}{\underline{Z}}) \text{ for all }k\ge -1,
\]
which proves the claim.
\end{proof}

Now we need to look at the reduced homologies of links of faces $F\in \mathcal{S}_X\setminus \D_X$.

\begin{lem}\label{lem:bad face has 0 homology link}
If $F\in\mathcal S_X$ and $F\notin\mathcal D_X$, then $\operatorname{lk}_{\mathcal S_X}(F)$ is a cone. In particular, $\widetilde H_i\bigl(\operatorname{lk}_{\mathcal S_X}(F)\bigr)=0$ for every $i$.
\end{lem}

\begin{proof}
Since $c(F)\in\mathcal D_X$ and $F\notin\mathcal D_X$, the inclusion $F\subseteq c(F)$ is strict. Choose $p\in c(F)\setminus F$.

For $G\in\operatorname{lk}_{\mathcal S_X}(F)$, choose a facet $\underline{Y}$ containing $F\cup G$. The definition of $c(F)$ gives $p\in \underline{Y}$, so $F\cup G\cup\{p\}\subseteq \underline{Y}$. Since $p\notin F$, it follows that $G\cup\{p\}\in\operatorname{lk}_{\mathcal S_X}(F)$. Thus the link is a cone with vertex $p$.
\end{proof}

Let us review the necessary background we need about Stanley--Reisner rings before proving Theorem \ref{thm:pdim=reg}.

\begin{dfn}\label{def:SR ring}
Let $\Delta$ be a finite simplicial complex and put $S(\Delta):=\Bbbk[x_v\mid v\in V(\Delta)]$, with $\deg x_v=1$. Its \emph{Stanley--Reisner ideal} and \emph{Stanley--Reisner ring} are
\[
I(\Delta):=\left(\prod_{v\in F}x_v\ \middle|\ F\subseteq V(\Delta),\ F\notin \Delta\right), \quad \text{and}
\quad \Bbbk[\Delta]:=S(\Delta)/I(\Delta) \text{ respectively}.
\]
Thus a product of distinct variables is non-zero in $\Bbbk[\Delta]$ precisely when its support is a face of $\Delta$.

Consider a minimal graded free resolution of $\Bbbk[\Delta]$ over $S(\Delta)$, and let $\beta_{i,j}$ be the multiplicity of $S(\Delta)(-j)$ in the $i$-th term, where $(-j)$ denotes the grading shift so that $1\in S(\Delta)(-j)$ is in degree $j$. The \emph{Castelnuovo--Mumford regularity} is
\[
\operatorname{reg}\Bbbk[\Delta]:=\max\{j-i\mid\beta_{i,j}\neq0\}.
\]
For the trivial complex $\Delta=\{\emptyset\}$, set $S(\Delta)=\Bbbk$, so $\operatorname{reg}\Bbbk[\Delta]=0$.
\end{dfn}

\begin{lem}\label{lem:Hochshster easier}{\rm\cite{HH11}}
For every finite simplicial complex $\Delta$,
\[
\operatorname{reg}\Bbbk[\Delta]
=\max\{k\geq0\mid\widetilde H_{k-1}(\operatorname{lk}_\Delta(F))\neq0
\text{ for some }F\in \Delta\}.
\]
\end{lem}

\begin{proof}
This is a reformulation of Hochster's formula (see \cite{Hoc77}) and an interpretation of Castelnuovo--Mumford regularity using local cohomology that seems to be known to experts but is not explicitly written in the literature in the setting when $\Bbbk$ is an arbitrary field; we clarify the reason for completeness and refer to \cite[Section A.7]{HH11} for unexplained terminology and references.  

Let $S:=S(\Delta)$ and $M:=\Bbbk[\Delta]=S(\Delta)/I(\Delta)$.
Consider the local cohomology $H^i_{\mathfrak{m}}(M)$, which is a $\ZZ^{V(\Delta)}$-graded module.  The Castelnuovo--Mumford regularity can be rephrased in terms of the totalized $\ZZ$-grading, that is,
\begin{align*}
\operatorname{reg}(M) 
&=\max\{i+j\mid H^i_{\mathfrak m}(M)_j\neq 0, \ j\in \ZZ\}
=\max \{ i+\sum_v \alpha_v \mid H^i_{\mathfrak m}(M)_\alpha \neq 0, \ \alpha\in \ZZ^{V(\Delta)}\}.
\end{align*}
We claim that
\[
\max \{ i+\sum_v \alpha_v \mid H^i_{\mathfrak m}(M)_\alpha \neq 0, \ \alpha\in \ZZ^{V(\Delta)}\}
=\max\{i-|F| \,:\, \widetilde H_{i-|F|-1}(\operatorname{lk}_\Delta F)\neq0,\ F\in \Delta\}.
\]
Indeed, Hochster's formula \cite[Theorem A.7.3]{HH11} states that
\begin{equation}\label{eq:hochster-local-cohomology}
\dim_{\Bbbk}H^i_{\mathfrak m}(M)_\alpha=
\begin{cases}
\dim_{\Bbbk}\widetilde H_{i-|F_\alpha|-1}(\operatorname{lk}_\Delta(F_\alpha)), & \alpha\leq0\text{ and }F_\alpha\in \Delta,\\
0, & \text{otherwise}.
\end{cases}
\end{equation}
Here $\alpha =(\alpha_v)_{v\in V(\Delta)} \le 0$ means that $\alpha_v\le 0$ for all $v\in V(\Delta)$, and $F_\alpha:=\{v\in V(\Delta)\mid\alpha_v<0\}$ is the face supported at $\alpha$.  Hence
\begin{align*}
\max\bigl\{i+{\textstyle\sum_v}\alpha_v \bigm| H^i_{\mathfrak m}(M)_\alpha\neq0\bigr\}
&=\max\bigl\{i+{\textstyle\sum_v}\alpha_v \bigm| \alpha\leq0,\ F_\alpha\in \Delta,\ \widetilde H_{i-|F_\alpha|-1}(\operatorname{lk}_\Delta (F_\alpha))\neq0\bigr\}\\
&=\max\bigl\{i-|F| \bigm| F\in \Delta,\ \widetilde H_{i-|F|-1}(\operatorname{lk}_\Delta(F))\neq0\bigr\},
\end{align*}
where the second equality follows from the observation that every $\alpha\leq0$ with $F_\alpha=F$ satisfies
\[
\sum_v\alpha_v=-\sum_{v\in F}|\alpha_v|\leq-|F|.
\]

Finally, $\widetilde{H}_{i-|F|-1}(\operatorname{lk}_\Delta (F))\neq 0$ implies that $i-|F|-1\ge -1$, and so we get that
\[
\max\bigl\{i-|F| \bigm| F\in \Delta,\ \widetilde H_{i-|F|-1}(\operatorname{lk}_\Delta (F))\neq0\bigr\} = \max \{ k\ge 0 \mid F\in \Delta\ \widetilde H_{k-1}(\operatorname{lk}_\Delta (F))\neq0 \}, 
\]
as required.
\end{proof}

\begin{proof}[Proof of Theorem \ref{thm:pdim=reg}]
This follows from
\begin{align*}
\operatorname{pdim}_{\Sigma_{\mathcal{J}(P)}} I_X
&=\max\{k\geq0\mid\widetilde H_{k-1}(\Delta^{>\underline{Z}})\neq0
             \text{ for some }Z\in\mathcal J(P)\}\\
&=\max\{k\geq0\mid\widetilde H_{k-1}(\operatorname{lk}_{\mathcal S_X}(F))\neq0
             \text{ for some }F\in\mathcal S_X\}\\
&=\operatorname{reg}\Bbbk[\mathcal S_X],
\end{align*}
where the first equality follows from Corollary \ref{cor:SigmaJ(P) is IG} and the last equality follows from Lemma \ref{lem:Hochshster easier}. For the middle equality, we get $\le$ from Lemma \ref{lem:H(Delta Z)=H(lk DZ)}, whereas $\ge$ follows from Lemma \ref{lem:bad face has 0 homology link}.
\end{proof}


\section{Observations about multiplicities}
\label{sec:low-degree-multiplicity}

We now investigate further properties of the minimal projective resolution constructed above.

For the minimal resolution of $I_X$ constructed in Theorem~\ref{thm:minimal-resolution}, equation~\eqref{eq:minimal-terms} gives
\begin{equation*}
 \bigl[\mathcal P_{\min,k}:P_Z\bigr]
 =\widetilde\beta_{Z,k}
 =\dim_{\Bbbk}\widetilde H_{k-1}(\Delta^{>\underline{Z}}).
\end{equation*}
In the following subsections, we exhibit some possible values of these multiplicities.

\subsection{Multiplicity-free in degrees zero and one}
\label{subsec:multi-free deg 0 1}
In degrees zero and one, the topology of $\Delta^{>\underline{Z}}$ is sufficiently restricted to give the following multiplicity-free results.  We recall that there is a difference between the void complex $\emptyset$ and the trivial complex $\{\emptyset\}$; see subsection~\ref{subsec:simplicial-topology}.

\begin{prop}\label{prop:low-degree-multiplicity}
The projective modules $\mathcal P_{\min,0}$ and $\mathcal P_{\min,1}$ are multiplicity-free.
\end{prop}

\begin{proof}
For $k=0$, the multiplicity of $P_Z$ is $\dim_{\Bbbk}\widetilde H_{-1}(\Delta^{>\underline{Z}})$, which is $1$ when $\Delta^{>\underline{Z}}=\{\emptyset\}$ is the trivial complex (i.e.\ $\D_X^{>\underline{Z}}=\emptyset$) and is $0$ otherwise.

Now consider $k=1$. If $\Delta^{>\underline{Z}}$ is trivial, then $\widetilde\beta_{Z,1}=0$; otherwise, we have
\begin{equation*}
 \widetilde\beta_{Z,1}
 =\dim_{\Bbbk}\widetilde H_0(\Delta^{>\underline{Z}})
 =\#\{\text{connected components of the $1$-skeleton of }\Delta^{>\underline{Z}}\}-1.
\end{equation*}
The $1$-skeleton of $\Delta^{>\underline{Z}}$ is a graph whose vertices are $\underline{Y}\in\D_X^{>\underline{Z}}$ where two vertices $\underline{Y},\underline{Y'}$ are adjacent precisely when they are comparable. It then suffices to show that the $1$-skeleton of $\Delta^{>\underline{Z}}$ has at most two connected components. Note that if $\D_X^{>\underline{Z}}=\emptyset$, then the $1$-skeleton of $\Delta^{>\underline{Z}}$ is empty and $\widetilde H_0(\Delta^{>\underline{Z}})=0$.

Every vertex $\underline{Y}\in \D_X^{>\underline{Z}}$ falls into one of the three cases (a) $Y\subsetneq Z$, (b) $Z\subsetneq Y$, or (c) $Y$ and $Z$ are incomparable.  In case (a), if $Y,Y'\subsetneq Z$, then $Y\cap Y'$ is an order ideal. Lemma~\ref{lem:symmetric-difference-identities}\textup{(\ref{itm:symdiff-inclusion})} gives $Z\setminus Y,Z\setminus Y'\subseteq X$. Similarly applying Lemma~\ref{lem:symmetric-difference-identities}\textup{(\ref{itm:symdiff-inclusion})} to $Y,Y\cap Y'$ and to $Y',Y\cap Y'$ yields
\[
 \underline{Y}\le \underline{Y\cap Y'},\qquad
 \underline{Y'}\le \underline{Y\cap Y'},\qquad
 \underline{Z}\lneqq \underline{Y\cap Y'}.
\]
Thus all vertices in case (a) lie in one connected component.

Case (b) follows from case (a) applied to $P^{\mathrm{op}}$ and $P\setminus X$: the complement map $Y\mapsto P\setminus Y$ identifies $\mathcal J(P)$ with $\mathcal J(P^{\mathrm{op}})$, reverses inclusions, preserves symmetric differences, and therefore leaves $\mathcal D_X$ unchanged but interchanges the cases (a) and (b).

Finally, suppose that $Y$ and $Z$ are incomparable. Since $\underline{Z}\subsetneq\underline{Y}$, Lemma~\ref{lem:symmetric-difference-identities}\textup{(\ref{itm:symdiff-interval})} gives
\[
 \underline{Y\cap Z}=\underline{Z}\sqcup(Z\setminus Y),\qquad
 \underline{Y\cup Z}=\underline{Z}\sqcup(Y\setminus Z).
\]
Both are vertices of the $1$-skeleton and are comparable with $\underline{Y}$; the first corresponds to an order ideal properly contained in $Z$, whereas the second corresponds to one properly containing $Z$.  Consequently, any $Y$ incomparable with $Z$ belongs to one of the two components described.
\end{proof}

\subsection{Higher multiplicities via independence complexes}
\label{subsec:higher-multiplicities}

In contrast to the behavior in degrees $0$ and $1$, we demonstrate that multiplicities can be arbitrarily large in higher degrees already when $P$ is a height-two poset.

Let $B=(X\sqcup Y,E)$ be a (finite) bipartite graph so we can regard $E\subseteq X\times Y$.  This defines a height-two poset $(P_B=X\sqcup Y, \le)$ where
\[
 x<y
 \quad\Leftrightarrow\quad
 \{x,y\}\in E
 \qquad \text{ for any } x\in X,\ y\in Y,
\]
with no other strict relations.  Conversely, the Hasse diagram of a height-two poset $(P=X\sqcup Y, \le)$, where $X$ is of height $0$ and $Y$ is of height $1$, is precisely the bipartite graph $B$ with $P_B=P$.

A subset $\sigma\subseteq V(B)$ is \defn{independent} if no two vertices of $\sigma$ are joined by an edge of $B$.  The \defn{independence complex} $\Ind(B)$ is the simplicial complex on $V(B)$ whose faces are the independent subsets of $B$.

\begin{prop}\label{prop:height-two-independence}
For a bipartite graph $B=(X\sqcup Y,E)$, the poset $\mathcal F\bigl(\Ind(B)\bigr)$ is the same as $\D_X$ associated to the height-two poset $P_B$. In particular, we have $\sd\Ind(B)=\Delta(\D_X\setminus\{\emptyset\}) = \Delta^{>\underline{X}}$, and hence, over $\Sigma_{\mathcal{J}(P_B)}$, we have
\[
 [\mathcal P_{\min,k}:P_X] = \beta_{X,k} = \dim_{\Bbbk}\widetilde H_{k-1}(\Delta^{>\underline{X}})
 \cong \dim_{\Bbbk}\widetilde H_{k-1}(\Ind(B))
 \;\;\text{ for all } k\ge 0.
\]
\end{prop}

\begin{proof}
Fix $T\in \mathcal{J}(P_B)$.  We show that $\underline{T}$ is an independent set of $B$.  Indeed, note that $\underline{T}=(X\setminus T)\sqcup(T\cap Y)$, so if, on the contrary, $x\in X\setminus T$ and $y\in T\cap Y$ were adjacent in $B$, then $x<y$ in $P_B$.  But $T$ being an order ideal implies that $x\in T$, a contradiction.

Conversely, let $U$ be an independent set of $B$; we claim that $T_U:=\bigl(X\setminus(U\cap X)\bigr)\cup(U\cap Y)\in \mathcal{J}(P_B)$.  Indeed, if $y\in U\cap Y$ and $x<y$ in $P_B$, then independence gives $x\notin U$, so $x\in T_U$. Since $T_U=X\SDiff U$, Lemma~\ref{lem:symmetric-difference-identities}\textup{(\ref{itm:symdiff-cancellation})} gives
$\underline{T_U}=X\SDiff(X\SDiff U)=U$.
Hence we get that $\mathcal F(\Ind(B))$ is the same collection of subsets as $\D_X$.  Since both posets use inclusion order, they are the same poset.

The last assertion follows from $\underline{X}=\emptyset$, which holds by Lemma~\ref{lem:symmetric-difference-identities}\textup{(\ref{itm:symdiff-associativity})}.
\end{proof}

By Proposition~\ref{prop:height-two-independence}, a bipartite graph $B$ with $\Ind(B)$ that has a reduced homology group of dimension at least $2$ gives a non-multiplicity-free term we desire.

\begin{eg}
Let $B$ be the bipartition graph $C_6\sqcup K_2$ (the disjoint union of a $6$-point cycle and a $2$-point complete graph) with parts $X=\{1,2,3,4\}$ and $Y=\{5,6,7,8\}$
so that $P_B$ is of the form:
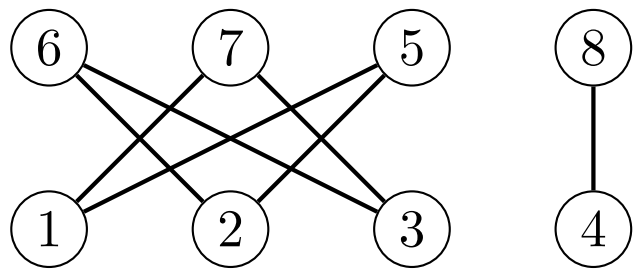
\begin{center}
\begin{tikzpicture}[scale=0.8,
  vertex/.style={
    circle,
    draw,
    fill=white,
    inner sep=0pt,
    minimum size=5mm
  }
]
 \node at (-2,1.5) {$P_B$:};\node[vertex] (1) at (0,0) {1};
  \node[vertex] (2) at (1.5,0) {2};
  \node[vertex] (3) at (3,0) {3};
  \node[vertex] (4) at (4.5,0) {4};
  \node[vertex] (5) at (3,1.5) {5};
  \node[vertex] (6) at (0,1.5) {6};
  \node[vertex] (7) at (1.5,1.5) {7};
  \node[vertex] (8) at (4.5,1.5) {8};
  \draw[thick] (5)--(1);
  \draw[thick] (5)--(2);
  \draw[thick] (6)--(2);
  \draw[thick] (6)--(3);
  \draw[thick] (7)--(3);
  \draw[thick] (7)--(1);
  \draw[thick] (8)--(4);
\end{tikzpicture}
\end{center}
We can draw $\Ind(C_6)$ as follows, where the 2-simplices are filled triangles $\{1,2,3\}$ and $\{5,6,7\}$, and there are three additional maximal faces $\{1,6\}, \{2,7\}, \{3,5\}$.
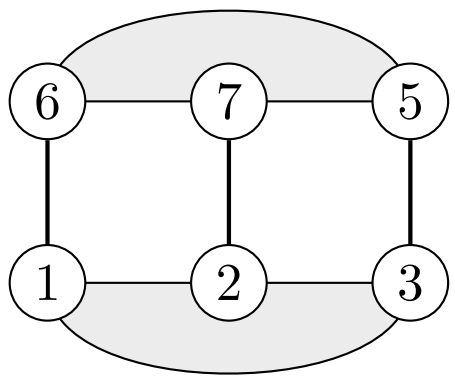
\begin{center}
\begin{tikzpicture}[scale=0.8,
  vertex/.style={
    circle,
    draw,
    fill=white,
    inner sep=0pt,
    minimum size=5mm
  }
]  
\node at (-2,1.5) {$\Ind(C_6)$:};
  \path[draw,fill=gray!15]
    (0,0) -- (1.5,0) -- (3,0)
    .. controls +(0,-1) and +(0,-1) .. (0,0)
    -- cycle;
  \path[draw,fill=gray!15]
    (0,1.5) -- (1.5,1.5) -- (3,1.5)
    .. controls +(0,1) and +(0,1) .. (0,1.5)
    -- cycle;
\node[vertex] (1) at (0,0) {1};
  \node[vertex] (2) at (1.5,0) {2};
  \node[vertex] (3) at (3,0) {3};
  \node[vertex] (5) at (3,1.5) {5};
  \node[vertex] (6) at (0,1.5) {6};
  \node[vertex] (7) at (1.5,1.5) {7};
  \draw[thick] (6)--(1);
  \draw[thick] (7)--(2);
  \draw[thick] (5)--(3);
\end{tikzpicture}
\end{center}
Contracting the two $2$-simplices yields a graph with two vertices joined by three edges, which means that we have homotopy equivalent space $|\Ind(C_6)|\simeq S^1\vee S^1$, where $\vee$ here denotes the wedge sum.  Also, we have $|\Ind(K_2)|=S^0$.

Note that an independent set of a disjoint union $G\sqcup H$ is the union of an independent set of $G$ and one of $H$, i.e. $\Ind(G\sqcup H)=\Ind(G)\ast \Ind(H)$, where $\ast$ denotes the simplicial join (faces are the union of faces in the two components).  Since joining with $S^0$ is homotopic to taking the suspension $\Sigma(-)$, and suspension preserves wedge sums, we get a homotopy equivalence
\begin{align*}
 |\Ind(B)|
 &=|\Ind(C_6)*\Ind(K_2)| \simeq(S^1\vee S^1)*S^0
 \simeq\Sigma(S^1\vee S^1)\simeq S^2\vee S^2.
\end{align*}
Hence, for $k\ge 0$, we have
\[
\widetilde{H}_k(\Ind(B)) \cong  \widetilde{H}_k(S^2\vee S^2) \cong \widetilde{H}_k(S^2)^{\oplus 2} = \begin{cases}
    \Bbbk^{\oplus 2}, & k=2, \\ 0, & \text{otherwise.}
\end{cases}
\]
In particular, we have $[\mathcal{P}_{\min,3}:P_X]=2$, i.e. $\mathcal P_{\min,3}$ is not multiplicity-free.
\end{eg}

In general, we can wedge more spheres to get arbitrary multiplicity in reduced homology.

\begin{eg}\label{eg:complexes-with-arbitrary-homology}
For $m\geq1$ and $k\ge 2$, take a finite triangulation $\Delta$ of $\bigvee_{j=1}^{m}S^{k-2}$ (for $k=2$ this is a set of $m+1$ points). Then $\widetilde H_{k-2}(\Delta)\cong\Bbbk^{\oplus m}$ and $\widetilde H_j(\Delta)=0$ for $j\neq k-2$.
\end{eg}

The final ingredient we need is the following realization theorem.

\begin{prop}\textnormal{(\cite[Proposition~6.2]{NR09})}
\label{prop:nagel-reiner}
For every finite simplicial complex $\Delta$, there exists a finite bipartite graph $B_\Delta$ such that $\Ind(B_\Delta)\simeq\Sigma \Delta$.
\end{prop}

\begin{cor}\label{cor:arbitrary-higher-multiplicities}
For every $k\geq2$ and every $m\geq1$, there exist a finite height-two poset $P$ and an order ideal $X\in\mathcal J(P)$ such that the multiplicity $[\mathcal{P}_{\min,k}:P_X]=m$ in the degree $k$ term of any minimal projective resolution $\mathcal{P}_{\min,\bullet}$ of the $\Sigma_{\mathcal{J}(P)}$-module $I_X$.
\end{cor}

\begin{proof}
As shown in Example~\ref{eg:complexes-with-arbitrary-homology}, we have $\dim_{\Bbbk}\widetilde H_{k-2}(\Delta)=m$ for any finite triangulation $\Delta$ of $\bigvee_{j=1}^m S^{k-2}$.  By Proposition~\ref{prop:nagel-reiner}, we have a bipartite graph $B_\Delta=X\sqcup Y$ so that $\Ind(B_\Delta)\simeq \Sigma \Delta$.  
Hence we have $\dim_{\Bbbk}\widetilde H_{k-1}(\Ind(B_\Delta))=m$ as suspension raises the degree of reduced homology by one. 
By Proposition~\ref{prop:height-two-independence}, over the distance algebra $\Sigma_{P_{\mathcal{J}(B_\Delta)}}$, we have $[\mathcal P_{\min,k}:P_X]=\dim_{\Bbbk}\widetilde H_{k-1}(\Ind(B_\Delta))=m$.
\end{proof}

\appendix

\section{A simplicial construction for thin modules}
\label{app:thin-modules}

We show that the projective resolution in Section~\ref{sec:candidate-resolution} is an instance of the simplicial replacement (bar construction of \cite[Definition 4.2.1]{Rie14}) of a diagram of projective modules.
This explains why it looks very much like a bar construction. 

Let $A=\Bbbk Q/I$ be a basic finite-dimensional bound-quiver algebra, and let $M$ be a \defn{thin} right $A$-module, i.e. $\dim_{\Bbbk}Me_x\leq1$ for every $x\in Q_0$.  Let $Q_M$ be the subquiver of $Q$ with vertex set $(Q_M)_0:=\{x\in Q_0\mid Me_x\neq0\}$ and whose arrows are those $a:x\to y$ for which right-multiplying by $a$ induces a non-zero map $Me_x\to Me_y$.  The key assumption we need is that non-zero vectors $m_x\in Me_x$ can be chosen so that
\begin{equation}
m_x a=m_y
\qquad
\text{for every arrow $a:x\to y$ in $Q_M$}.
\label{eq:thin-normalized-basis}
\end{equation}

The quiver $Q_M$ is acyclic.  Indeed, if it contained an oriented cycle, then \eqref{eq:thin-normalized-basis} would imply that every positive power of this cycle acts non-trivially on one of the vectors $m_x$.  This contradicts the admissibility of $I$, since every sufficiently long path belongs to $I$.

For every $x\in(Q_M)_0$, let $\pi_x\colon P_x=e_xA\to M$ be the homomorphism determined by $\pi_x(e_x)=m_x$.  Since each non-zero space $Me_x$ is spanned by $m_x$, these maps give a surjection
\begin{equation}\label{eq: pi to M}
d_0:=\bigl(\pi_x\bigr)_{x\in(Q_M)_0}\colon
\bigoplus_{x\in(Q_M)_0}P_x\to M.
\end{equation}

Consider the small category $\mathcal C_M$ associated to $M$ given as follows.  We take $(Q_M)_0$ as the objects of $\mathcal{C}_M$. For $x,y\in(Q_M)_0$, define
\[
\mathcal C_M(y,x)
:= \left\{ \mu_p:P_y\to P_x \mid p:x\to y\text{ is a path on }Q_M \right\},
\]
where $\mu_p:P_y\to P_x$ is the homomorphism of projective $A$-modules given by left-multiplying $p$. In particular, two paths $p,p':x\to y$ determine the same morphism of $\mathcal C_M$ precisely when $\mu_p=\mu_{p'}$, and the identity of $x$ is $\mu_{e_x} = \id_{P_x}$. If there is no path between $x,y$ on $Q_M$, then we take $\mathcal{C}_M(y,x)=\emptyset$. Composition is given by $\mu_p\mu_q=\mu_{pq}$ whenever $p:x\to y$ and $q:y\to z$.

Define a functor
\begin{equation*}
F_M\colon\mathcal C_M\longrightarrow\proj A,
\qquad
x\longmapsto P_x,
\qquad
\mu_p\longmapsto\mu_p.
\end{equation*}
By \eqref{eq:thin-normalized-basis}, it follows that $\pi_x\mu_p=\pi_y$, so the homomorphisms $\pi_x$ form a cocone over $F_M$ with vertex $M$, which means that there is a surjection induced by $d_0$ of \eqref{eq: pi to M}
\[
\colim_{\mathcal{C}_M} F_M
\cong \Big(\bigoplus_{x\in \mathcal{C}_M} P_x\Big)
\Big/ \sum_{\substack{p:x\to y\\ m\in P_y}} (m-pm)A
\twoheadrightarrow M.
\]

For $k\geq0$, let $N_k(\mathcal C_M)$ be the set of composable strings of $k$ non-identity morphisms
\begin{equation}\label{eq:x ast string}
x_\ast=
\bigl(
x_0\xrightarrow{f_1}x_1\xrightarrow{f_2}\cdots
\xrightarrow{f_k}x_k
\bigr).
\end{equation}
In particular, $N_0(\mathcal C_M)$ is the object set of $\mathcal C_M$.  For $0\leq r\leq k$, write $x_\ast^{\widehat r}$ for the string obtained by deleting $x_r$ and (only in the case of $r\neq 0,k$) composing the adjacent morphisms.  Define
\begin{equation*}
\mathcal P_k(M):=
\bigoplus_{x_\ast\in N_k(\mathcal C_M)}
F_M(x_0)^{(x_\ast)} = \bigoplus_{x_\ast\in N_k(\mathcal C_M)}
P_{x_0}^{(x_\ast)}.
\end{equation*}
Here \(F_M(x_0)^{(x_\ast)} = P_{x_0}^{(x_\ast)}\) denotes a copy of \(F_M(x_0)=P_{x_0}\) indexed by $x_\ast$.  For $k>0$, define
\[
d_k\colon\mathcal P_k(M)\longrightarrow\mathcal P_{k-1}(M),
\qquad u^{(x_\ast)}
\mapsto
\bigl(F_M(f_1)(u)\bigr)^{(x_\ast^{\widehat0})}
+\sum_{r=1}^k(-1)^r
u^{(x_\ast^{\widehat r})}.
\]
Thus the zeroth face changes the initial projective through $F_M(f_1)$, while every other face acts as the identity on the initial projective.

\begin{prop}\label{prop:cpx from module}
For a thin $A$-module $M$, the maps $d_k$ above yield a complex of $A$-modules
\begin{equation*}
\cdots\longrightarrow\mathcal P_2(M)
\xrightarrow{d_2}\mathcal P_1(M)
\xrightarrow{d_1}\mathcal P_0(M)
\xrightarrow{d_0}M\to 0,
\end{equation*}
which is exact at $\mathcal{P}_0(M)$ when $M$ is faithful.
\end{prop}

\begin{proof}
Fix $k\geq1$ and a string $x_\ast$ as in \eqref{eq:x ast string}. For $0\leq r\leq k$, let $\varepsilon_r^{x_\ast}\colon P_{x_0}^{(x_\ast)}\to\mathcal P_{k-1}(M)$ be the $r$-th face map, that is,
\begin{equation}\label{eq:general face map}
\varepsilon_r^{x_\ast}\bigl(u^{(x_\ast)}\bigr):=\begin{cases}\bigl(F_M(f_1)(u)\bigr)^{(x_\ast^{\widehat0})},&r=0,\\[1mm]u^{(x_\ast^{\widehat r})},&r>0,\end{cases}
\qquad\text{so that}\qquad
d_k\big|_{P_{x_0}^{(x_\ast)}}=\sum_{r=0}^k(-1)^r\varepsilon_r^{x_\ast}.
\end{equation}

Functoriality of $F_M$ gives the simplicial identities
\begin{equation*}
\varepsilon_i^{x_\ast^{\widehat j}}
\varepsilon_j^{x_\ast}
=
\varepsilon_{j-1}^{x_\ast^{\widehat i}}
\varepsilon_i^{x_\ast}
\qquad
(0\leq i<j\leq k).
\end{equation*}
The only identity involving two non-identity homomorphisms is the one obtained by deleting $x_0$ and $x_1$, and it follows from $F_M(f_2)F_M(f_1)=F_M(f_2f_1)$.  For $0\leq i<j\leq k$, the equal composites $\varepsilon_i^{x_\ast^{\widehat j}}\varepsilon_j^{x_\ast}$ and $\varepsilon_{j-1}^{x_\ast^{\widehat i}}\varepsilon_i^{x_\ast}$ occur in $d_{k-1}d_k$ with coefficients $(-1)^{i+j}$ and $(-1)^{i+j-1}$, respectively. Hence $d_{k-1}d_k=0$ for $k\geq2$.

For a morphism $f:x_0\to x_1$ and $u\in P_{x_0}$, the equality $\pi_{x_1}\circ F_M(f)=\pi_{x_0}$ gives
\[
d_0d_1\bigl(u^{(f)}\bigr)
=
\pi_{x_1}F_M(f)(u)-\pi_{x_0}(u)=0.
\]
Hence the displayed sequence is a complex.

Suppose now that $M$ is faithful. Then $Q_M=Q$, since every vertex idempotent and every arrow acts non-trivially on $M$. Let $\iota_x:P_x\hookrightarrow\mathcal P_0(M)$ be the canonical inclusion. Let $\overline{d}_0:\operatorname{coker}d_1\to M$ be the map induced by $d_0$. For every non-trivial path $p:x\to y$, we have
\[
d_1\bigl(e_y^{(\mu_p)}\bigr)
=\iota_x(p)-\iota_y(e_y).
\]
Since path classes span each $e_xAe_y$, it follows that $(\operatorname{coker}d_1)e_y$ is spanned by $[\iota_y(e_y)]$. This class maps to $m_y\neq0$, which spans $Me_y$. Hence the restriction $\overline{d}_0|_{(\operatorname{coker}d_1)e_y}:(\operatorname{coker}d_1)e_y\to Me_y$
is an isomorphism for every $y\in Q_0$. Therefore $\overline{d}_0$ is an isomorphism, and $\ker d_0=\operatorname{im}d_1$.
\end{proof}

\begin{eg}\label{eg:DX is C_IX}
For $A=\Sigma_{\mathcal J(P)}$ and $M=I_X$, choose $m_Y\in I_Xe_Y$ to be the functional dual to $b_{YX}$.  Then $\mathcal C_{I_X}$ can be identified with the poset category $\D_X$: there is a non-identity morphism $Z\to T$ precisely when $\underline{Z}<\underline{T}$, and this morphism is represented by the left-multiplication map $\mu_{T,Z}:P_Z\to P_T$ given by $a\mapsto b_{TZ}a$. The set $N_k(\mathcal{C}_M)$ becomes the $k$-chains in $\Delta(\D_X)$.
Note that, when $P\neq\emptyset$, $I_X$ is not faithful, so the faithfulness condition in Proposition \ref{prop:cpx from module} is only sufficient and not necessary.
\end{eg}

\begin{eg}\label{eg:non exact cpx from simple}
For a simple module $S_i$, the category $\mathcal C_{S_i}$ has only the object $i$ and its identity morphism. Hence $\mathcal P_0(S_i)=P_i$ and $\mathcal P_k(S_i)=0$ for $k>0$, so the construction reduces to the projective cover $P_i\to S_i$.  It is a projective resolution only when $S_i$ is itself projective. When $S_i$ is not projective, the complex $\mathcal P_\bullet(S_i)$ is not a projective resolution.
\end{eg}

Now we relate $\mathcal P_\bullet(M)$ to the simplicial bar construction of \cite[Definition~4.2.1]{Rie14}.  For $k\geq0$, put
\[
B_k:=\bigoplus_{x_0\xrightarrow{f_1}\cdots\xrightarrow{f_k}x_k}F_M(x_0)^{(x_\ast)}
 \qquad(\text{identity morphisms allowed among the }f_r),
\]
where the sum is, similar to \eqref{eq:x ast string}, over all strings of $k$ (\emph{possibly identity}) composable maps. 
This comes with face maps $\varepsilon_r\colon B_k\to B_{k-1}$ given by the same formula \eqref{eq:general face map}, and degeneracy maps $s_r\colon B_{k-1}\to B_k$ ($0\leq r\leq k-1$), which insert $\id_{x_r}$ into the string. 
This makes $B_\bullet$ a \emph{simplicial} right $A$-module; the following table gives the detailed translation between our notations and those in \cite{Rie14}.

\begin{center}
\begin{tabular}{p{0.28\textwidth}p{0.64\textwidth}}
\hline
Bar construction setting & Our setting\\
\hline
$\mathcal D$ & $\mathcal C_M$\\[3pt]
$\mathcal M$ & $\operatorname{Mod}A$, the category of right $A$-modules\\[3pt]
$F\colon\mathcal D\to\mathcal M$ & $F_M\colon\mathcal C_M\to\operatorname{proj}A\subseteq\operatorname{Mod}A$\\[3pt]
$G\colon\mathcal D^{\mathrm{op}}\to\mathrm{Set}$ & the constant functor with value a one-point set $\ast$\\[3pt]
$d_0\to\cdots\to d_k$ & a string $x_\ast$, identities allowed\\[3pt]
$Gd_k\otimes Fd_0$ & $\ast\otimes F_M(x_0)\cong F_M(x_0)$\\[3pt]
$B_k(G,\mathcal D,F)$ & $B_k$\\[3pt]
face and degeneracy maps & $\varepsilon_r$ and $s_r$\\
\hline\\[3pt]
\end{tabular}
\end{center}

Having a simplicial module means that taking the alternating sum of the face maps
\[
 d^B_k:=\sum_{r=0}^k(-1)^r\varepsilon_r\colon B_k\longrightarrow B_{k-1}
\]
makes $B_\bullet$ a chain complex: $d^B_{k-1}d^B_k=0$ follows from the simplicial identities, by the same computation as in the proof of Proposition~\ref{prop:cpx from module}. 
We now compare it with $\mathcal P_\bullet(M)$ using the decomposition above.
Partitioning the direct summands of $B_k$ according to whether their labeling strings contain an identity morphism gives
\[
 B_k=\mathcal P_k(M)\oplus D_k,
 \qquad
 D_k:=\bigoplus_{\substack{x_\ast\text{ with }f_r=\id\\ \text{for some }r}}F_M(x_0)^{(x_\ast)}.
\]
\begin{prop}\label{prop:bar-construction-decomposition}
The decomposition $B_k=\mathcal P_k(M)\oplus D_k$ extends to one on the level of chain complexes $B_\bullet=\mathcal P_\bullet(M)\oplus D_\bullet$ with $D_\bullet$ acyclic.
In particular, the inclusion $\mathcal P_\bullet(M)\hookrightarrow B_\bullet$ induces an isomorphism in homology.
\end{prop}

\begin{proof}
The decomposition follows from the fact that face maps of $B_k$'s preserve (1) the property of having $f_r=\id$ for some $r$ in the labeling string, and also (2) the property of not having any identity morphism in the labeling.
Thus the differential restricts to a differential on the two direct summands in the decomposition; in particular, the restriction on $\mathcal{P}_k(M)$ yields exactly the differential of $\mathcal{P}_\bullet(M)$.

(1) We show that $d^B(D_k)\subseteq D_{k-1}$.  Indeed, suppose that $x_\ast$ is a string with $f_{r+1}=\id$. Then both $\varepsilon_r$ and $\varepsilon_{r+1}$ remove this identity and produce the same string and the same map, which appear in $d^B$ with opposite signs $(-1)^r$ and $(-1)^{r+1}$.  Hence, they cancel each other.  Every other face of $x_\ast$ still contains the identity morphism $f_{r+1}$, so $d^B(D_k)\subseteq D_{k-1}$.

(2) Now we show that $d^B|_{\mathcal{P}_k(M)}$ is precisely the differential $d_k:\mathcal{P}_k(M)\to mathcal{P}_{k-1}(M)$.  Indeed, suppose that $x_\ast$ contains no identity morphisms in its letters; then the acyclicity of $Q_M$ implies that none of the compositions $f_rf_{r+1}$ can be identity.  Thus, $d^B(\mathcal{P}_k(M))\subseteq \mathcal{P}_{k-1}(M)$, and the claim follows by direct comparison of the differentials.

(3) Finally, since the image of the degeneracy map $s_r$ is spanned by the strings with $f_{r+1}=\id$, \cite[\S8.3]{Wei94} says that $D_\bullet$ is a subcomplex of $B_\bullet$ that is acyclic.
\end{proof}

The last part of Proposition \ref{prop:bar-construction-decomposition} says that we have
\begin{equation*}
H_0\bigl(\mathcal P_\bullet(M)\bigr)
\cong\colim_{\mathcal C_M}F_M,
\qquad
H_k\bigl(\mathcal P_\bullet(M)\bigr)
\cong L_k\colim_{\mathcal C_M}F_M
\quad(k>0),
\end{equation*}
where $L_k$ denotes the $k$-th left-derived functor. In the general setting of $\mathcal{C}_M$, these higher derived colimits may not vanish and the colimit may not be $M$ itself, whereas Theorem~\ref{thm:exactness} proves that, in the case when $\mathcal{C}_M=\D_X$ (Example \ref{eg:DX is C_IX}), we have $\colim_{\D_X} F_{I_X} \cong I_X$ and $L_k\colim_{\D_X} F_{I_X}=0$ for all $k>0$.


\section{Proofs of Lemma~\ref{lem:lattice-canonical-duality} and Corollary~\ref{cor:relative-ordinary-injective-dimension}}
\label{app:lattice-proofs}
\begin{proof}[Proof of Lemma~\ref{lem:lattice-canonical-duality}]
The functor $D$ is exact on $R$-split short exact sequences, and evaluation gives a natural isomorphism $M\xrightarrow{\sim}D^2(M)$ for every $M\in\lat(A)$. This proves the duality assertion.

For $M\in\lat(A)$, choose a surjection $A^{\oplus m}\twoheadrightarrow M$ with $m$ finite. It splits over $R$, so its kernel is again an $A$-lattice. Thus $\lat(A)$ has enough projectives, and its projective objects are precisely $\add A$. The assertion about enough injective follows by applying $D$ to left projective modules.  Since the syzygy of $M\in \lat(A)$ stays in $\lat(A)$, any projective resolution in $\lat(A)$ is just an ordinary projective resolution; this last assertion follows.
\end{proof}

\begin{proof}[Proof of Corollary~\ref{cor:relative-ordinary-injective-dimension}]
Since $A$ is finitely generated over the noetherian ring $R$, it is noetherian on both sides. By Lemma~\ref{lem:lattice-canonical-duality}, the exact duality $D:\lat(A)^{\operatorname{op}}
 \xrightarrow{\sim}\lat(A^{\operatorname{op}})$
swaps relative injective coresolutions of $A_A$ with projective resolutions of ${}_A D(A)$. Since projective dimension in the lattice category agrees with ordinary projective dimension, we obtain $\idim^{\lat}_A A
 =\pdim_{A^{\operatorname{op}}}D(A).$
Applying the same argument to $A^{\operatorname{op}}$ gives the other relative dimension equality.

Assume now that $\idim_R R<\infty$. 
We claim that $\idim D(A)_A\leq\idim_R R$.  Indeed, for any right $A$-module $N$, so long as $A$ is $R$-projective, any $A$-projective resolution of $N$ is also an $R$-projective resolution.  So the tensor-Hom adjunction $\Hom_A(-,D(A)_A)\cong\Hom_R(-,R)$ gives
\[
 \Ext_A^q(N,D(A)_A)
 \cong \Ext_R^q(N,R)
 \qquad(q\geq0),
\]
and the claim follows.

Put $n:=\pdim_{A^{\operatorname{op}}}D(A)$. The bound for $\idim A_A$ is immediate if $n=\infty$. If $n<\infty$, choose a resolution by finitely generated projective left $A$-modules
\[
 0\longrightarrow Q_n\longrightarrow\cdots
 \longrightarrow Q_0\longrightarrow {}_A D(A)
 \longrightarrow0.
\]
The modules $Q_i$ and $D(A)$ are projective over $R$. Starting with the surjection $Q_0\to D(A)$, we see inductively that every syzygy is projective over $R$ and every successive short exact sequence splits over $R$. Applying $D$ and using $D^2(A)\cong A$ then yields an exact sequence of right $A$-modules
\[
 0\longrightarrow A_A\longrightarrow D(Q_0)
 \longrightarrow\cdots\longrightarrow D(Q_n)
 \longrightarrow0.
\]
For every $0\leq i\leq n$, we have $D(Q_i)\in\add D(A)_A$, and hence $\idim_A D(Q_i)\leq\idim_R R$ by the claim in the previous paragraph. Dimension shifting along this exact sequence yields
\[
 \idim A_A
 \leq
 \max_{0\leq i\leq n}
 \bigl\{\idim_A D(Q_i)+i\bigr\}
 \leq \idim_R R+n.
\]
Applying the same argument to $A^{\operatorname{op}}$ gives
$\idim {}_A A\leq\idim_R R+\pdim_A D(A)$.  Thus, if both $\pdim_{A^{\operatorname{op}}}D(A)$ and $\pdim_{A}D(A)$ are finite, $A$ is Iwanaga--Gorenstein.

Finally, when $R=\Bbbk$, we have $\idim_R R=0$, and relative injective dimension agrees with ordinary injective dimension.  Putting this into the inequalities, the claim follows immediately.
\end{proof}

\bibliographystyle{alpha}
\bibliography{reference.bib}
\Addresses

\end{document}

%% file: pre.tex
\usepackage[margin=30truemm]{geometry}
\usepackage{fancyhdr,appendix,graphicx,layout}
\usepackage{amsmath,amsthm,amssymb}
\usepackage[shortlabels]{enumitem}
\setlist[enumerate]{font=\normalfont}
\usepackage{physics}

\usepackage[numbers]{natbib}
\usepackage{bbm}
\usepackage{bm}
\usepackage{tikz}
\usepackage{tikz-cd}
\usetikzlibrary{positioning,arrows.meta,decorations.markings,calc}
\tikzset{
    cells={font=\everymath\expandafter{\the\everymath\displaystyle}},
}
\usepackage{caption}
\usepackage[pagewise]{lineno}

\usepackage[normalem]{ulem}
\usepackage{float}

\usepackage{color}
\definecolor{darkgreen}{rgb}{0,0.7,0}

\definecolor{darkblue}{rgb}{0,0,0.7} 
\newcommand{\defn}[1]{\textsl{\color{darkblue} #1}}

\usepackage[hyperpageref]{backref}
\usepackage[pagebackref]{hyperref}
\hypersetup{
    colorlinks=true,
    linkcolor=red,  
    filecolor=blue,  
    urlcolor=orange, 
    citecolor=violet   
}
\renewcommand*{\backref}[1]{}  
   \renewcommand*{\backrefalt}[4]{
      \ifcase #1 
         Not cited.
      \or
         Cited on page #2.
      \else
         Cited on pages #2.
      \fi}

\usepackage{cleveref}

\makeatletter
\newcommand\ReDeclareMathOperator[2]{%
    \begingroup \escapechar\m@ne\xdef\@gtempa{{\string#1}}\endgroup
    \expandafter\@ifundefined\@gtempa
        {\@latex@error{Command \string#1 undefined}\@ehc}
        \relax
    \let\@ifdefinable\@rc@ifdefinable
    \DeclareMathOperator#1{#2}}
\makeatother

\counterwithin{equation}{subsection}  
\crefname{thm}{Theorem}{Theorems}
\crefname{dfn}{Definition}{Definitions}
\crefname{dfnprop}{Definition-Proposition}{Definition-Proposition}
\crefname{prop}{Proposition}{Propositions}
\crefname{lem}{Lemma}{Lemmas}
\crefname{cor}{Corollary}{Corollaries}
\crefname{clm}{Claim}{Claims}
\crefname{ass}{Assumption}{Assumption}
\crefname{cond}{Condition}{Condition}
\crefname{nota}{Notation}{Notation}
\crefname{conj}{Conjecture}{Conjecture}
\crefname{fct}{Fact}{Facts}
\crefname{rmk}{Remark}{Remarks}
\crefname{eg}{Example}{Examples}
\crefname{que}{Question}{Questions}
\crefname{figure}{Figure}{Figures}
\crefname{table}{Table}{Tables}
\crefname{section}{Section}{Sections}
\crefname{subsection}{Subsection}{Subsections}
\crefname{appendix}{Appendix}{Appendices}
\crefname{equation}{}{}
\crefname{main}{Theorem}{Theorems}

\theoremstyle{plain}
\newtheorem{thm}{Theorem}[section]
\newtheorem{main}{Theorem}
\newtheorem{prop}[thm]{Proposition}
\newtheorem{lem}[thm]{Lemma}
\newtheorem{cor}[thm]{Corollary}

\theoremstyle{definition}
\newtheorem{dfn}[thm]{Definition}

\newtheorem{eg}[thm]{Example}
\newtheorem{que}[thm]{Question}

\theoremstyle{remark}
\newtheorem{rmk}[thm]{Remark}

\newcommand{\D}{\mathcal D}

\newcommand{\ZZ}{\mathbb{Z}}

\DeclareMathOperator{\add}{add}

\DeclareMathOperator{\rad}{rad}
\ReDeclareMathOperator{\top}{top}

\DeclareMathOperator{\id}{id}

\DeclareMathOperator{\proj}{proj}

\DeclareMathOperator{\colim}{colim}

\DeclareMathOperator{\Hom}{Hom}

\DeclareMathOperator{\Ext}{Ext}

\ReDeclareMathOperator{\l}{\ell}

\newcommand{\Ker}{\operatorname{Ker}}
\newcommand{\Image}{\operatorname{Im}}

\newcommand{\morsethickarrow}{\rightarrow}
\newcommand{\morsedottedarrow}{\rightleftharpoons}

\DeclareMathOperator{\idim}{idim}
\DeclareMathOperator{\pdim}{pdim}

\DeclareMathOperator{\Crit}{Crit}
\DeclareMathOperator{\Ind}{Ind}
\DeclareMathOperator{\sd}{sd}
\DeclareMathOperator{\weight}{wt}

\newcommand{\SDiff}{{\mathbin{\triangle}}}

\DeclareMathOperator{\lat}{lat}

\renewcommand{\leq}{\leqslant}
\renewcommand{\le}{\leqslant}

%% file: reference.bib
@misc{AI24,
  author        = {Asao, Y. and Ivanov, S. O.},
  title         = {Magnitude homology is a derived functor},
  year          = {2024},
  eprint        = {2402.14466},
  archivePrefix = {arXiv},
  note          = {arXiv:2402.14466}
}

@article{AW24,
  author  = {Asao, Y. and Wakatsuki, S.},
  title   = {Minimal projective resolution and magnitude homology of geodetic metric spaces},
  journal = {Homology Homotopy Appl.},
  volume  = {28},
  number  = {2},
  pages   = {123--154},
  year    = {2026},
  doi     = {10.4310/hha.2026.v28.n2.a7}
}

@book{AB69,
  author    = {Auslander, M. and Bridger, M.},
  title     = {Stable module theory},
  series    = {Mem. Amer. Math. Soc.},
  volume    = {94},
  publisher = {American Mathematical Society},
  address   = {Providence, RI},
  year      = {1969}
}

@article{Bas63,
  author  = {Bass, H.},
  title   = {On the ubiquity of {G}orenstein rings},
  journal = {Math. Z.},
  volume  = {82},
  pages   = {8--28},
  year    = {1963},
  doi     = {10.1007/BF01112819}
}

@article{Bir37,
  author  = {Birkhoff, G.},
  title   = {Rings of sets},
  journal = {Duke Math. J.},
  volume  = {3},
  number  = {3},
  pages   = {443--454},
  year    = {1937}
}

@article{BBIKP25,
  author  = {Barthel, T. and Benson, D. J. and Iyengar, S. B. and Krause, H. and Pevtsova, J.},
  title   = {Lattices over finite group schemes and stratification},
  journal = {Compos. Math.},
  volume  = {161},
  number  = {11},
  pages   = {2911--2946},
  year    = {2025},
  doi     = {10.1112/S0010437X25102844}
}

@incollection{Bjo95,
  author    = {Bj{\"o}rner, A.},
  title     = {Topological methods},
  booktitle = {Handbook of combinatorics},
  volume    = {1--2},
  pages     = {1819--1872},
  publisher = {Elsevier},
  address   = {Amsterdam},
  year      = {1995}
}

@book{Buc21,
  author    = {Buchweitz, R.-O.},
  title     = {Maximal {Cohen}--{Macaulay} Modules and {Tate} Cohomology},
  series    = {Mathematical Surveys and Monographs},
  volume    = {262},
  publisher = {American Mathematical Society},
  address   = {Providence, RI},
  year      = {2021},
  doi       = {10.1090/surv/262},
  note      = {Published version of the 1986 preprint}
}

@article{CLZ26,
  author  = {Chen, J. and Liu, Y. and Zhou, G.},
  title   = {Algebraic {M}orse theory via homological perturbation lemma},
  journal = {J. Pure Appl. Algebra},
  volume  = {230},
  number  = {11},
  pages   = {108371},
  year    = {2026},
  doi     = {10.1016/j.jpaa.2026.108371}
}

@article{CKL16,
  author  = {Chuang, J. and King, A. and Leinster, T.},
  title   = {On the magnitude of a finite dimensional algebra},
  journal = {Theory Appl. Categ.},
  volume  = {31},
  pages   = {Paper No. 3, 63--72},
  year    = {2016}
}

@book{CR81,
  author    = {Curtis, C. W. and Reiner, I.},
  title     = {Methods of representation theory: With applications to finite groups and orders},
  volume    = {I},
  publisher = {John Wiley \& Sons},
  address   = {New York},
  year      = {1981}
}

@unpublished{CX26,
  author = {Chan, A. and Xing, B.},
  title  = {On representation-theoretic and homological properties of distance algebras of finite graphs},
  note   = {In preparation},
  year   = {2026}
}

@misc{GIK23,
  author        = {Gnedin, W. and Iyengar, S. B. and Krause, H.},
  title         = {A class of {Gorenstein} algebras and their dualities},
  year          = {2023},
  eprint        = {2303.04893},
  archivePrefix = {arXiv},
  primaryClass  = {math.RT},
  note          = {arXiv:2303.04893}
}

@article{GL89,
  author  = {Gugenheim, V. K. A. M. and Lambe, L. A.},
  title   = {Perturbation theory in differential homological algebra. {I}},
  journal = {Illinois J. Math.},
  volume  = {33},
  number  = {4},
  pages   = {566--582},
  year    = {1989}
}

@article{GN01,
  author  = {Goto, S. and Nishida, K.},
  title   = {Finite modules of finite injective dimension over a {Noetherian} algebra},
  journal = {J. London Math. Soc. (2)},
  volume  = {63},
  number  = {2},
  pages   = {319--335},
  year    = {2001}
}

@book{Gra03,
  author    = {Gr{\"a}tzer, G.},
  title     = {General lattice theory},
  edition   = {2},
  publisher = {Birkh{\"a}user},
  address   = {Basel},
  year      = {2003}
}

@incollection{Hap91,
  author    = {Happel, D.},
  title     = {On {G}orenstein algebras},
  booktitle = {Representation theory of finite groups and finite-dimensional algebras ({B}ielefeld, 1991)},
  series    = {Progr. Math.},
  volume    = {95},
  pages     = {389--404},
  publisher = {Birkh\"{a}user},
  address   = {Basel},
  year      = {1991}
}

@book{Hat02,
  author    = {Hatcher, Allen},
  title     = {Algebraic Topology},
  publisher = {Cambridge University Press},
  year      = {2002},
  url       = {https://pi.math.cornell.edu/~hatcher/AT/AT.pdf}
}

@incollection{Hoc77,
  author    = {Hochster, M.},
  title     = {{Cohen}--{Macaulay} rings, combinatorics,
               and simplicial complexes},
  booktitle = {Ring theory, {II}
               (Proc. Second Conf., Univ. Oklahoma,
               Norman, Okla., 1975)},
  series    = {Lecture Notes in Pure and Applied Mathematics},
  volume    = {26},
  pages     = {171--223},
  publisher = {Marcel Dekker},
  address   = {New York},
  year      = {1977}
}

@book {HH11,
    AUTHOR = {Herzog, J\"{u}rgen and Hibi, Takayuki},
     TITLE = {Monomial ideals},
    SERIES = {Graduate Texts in Mathematics},
    VOLUME = {260},
 PUBLISHER = {Springer-Verlag London, Ltd., London},
      YEAR = {2011},
     PAGES = {xvi+305},
      ISBN = {978-0-85729-105-9},
   MRCLASS = {13D02 (05E40 13D40 13F55 13P10)},
  MRNUMBER = {2724673},
MRREVIEWER = {Rahim\ Zaare-Nahandi},
       DOI = {10.1007/978-0-85729-106-6},
       URL = {https://doi.org/10.1007/978-0-85729-106-6},
}

@article{HW17,
  author  = {Hepworth, R. and Willerton, S.},
  title   = {Categorifying the magnitude of a graph},
  journal = {Homology Homotopy Appl.},
  volume  = {19},
  number  = {2},
  pages   = {31--60},
  year    = {2017},
  doi     = {10.4310/HHA.2017.v19.n2.a3}
}

@article{IK22,
  author  = {Iyengar, Srikanth B. and Krause, Henning},
  title   = {The {Nakayama} functor and its completion for {Gorenstein} algebras},
  journal = {Bulletin de la Soci\'{e}t\'{e} Math\'{e}matique de France},
  volume  = {150},
  number  = {2},
  pages   = {347--391},
  year    = {2022},
  doi     = {10.24033/bsmf.2849},
  eprint  = {2010.05676},
  archivePrefix = {arXiv}
}

@article{IKWY15,
  author  = {Iyama, O. and Kalck, M. and Wemyss, M. and Yang, D.},
  title   = {Frobenius categories, {Gorenstein} algebras and rational surface singularities},
  journal = {Compos. Math.},
  volume  = {151},
  number  = {3},
  pages   = {502--534},
  year    = {2015},
  doi     = {10.1112/S0010437X14007647}
}

@article{Iwa79,
  author  = {Iwanaga, Y.},
  title   = {On rings with finite self-injective dimension},
  journal = {Comm. Algebra},
  volume  = {7},
  number  = {4},
  pages   = {393--414},
  year    = {1979},
  doi     = {10.1080/00927877908822356}
}

@article{JW09,
  author  = {J{\"o}llenbeck, M. and Welker, V.},
  title   = {Minimal resolutions via algebraic discrete {Morse} theory},
  journal = {Mem. Amer. Math. Soc.},
  volume  = {197},
  number  = {923},
  year    = {2009},
  doi     = {10.1090/memo/0923}
}

@article{Lei08,
  author  = {Leinster, T.},
  title   = {The {E}uler characteristic of a category},
  journal = {Doc. Math.},
  volume  = {13},
  pages   = {21--49},
  year    = {2008}
}

@article{Lei19,
  author  = {Leinster, T.},
  title   = {The magnitude of a graph},
  journal = {Math. Proc. Cambridge Philos. Soc.},
  volume  = {166},
  number  = {2},
  pages   = {247--264},
  year    = {2019}
}

@incollection{LM17,
  author    = {Leinster, T. and Meckes, M. W.},
  title     = {The magnitude of a metric space: from category theory to geometric measure theory},
  booktitle = {Measure theory in non-smooth spaces},
  series    = {Partial Differ. Equ. Meas. Theory},
  pages     = {156--193},
  publisher = {De Gruyter Open},
  address   = {Warsaw},
  year      = {2017}
}

@article{LW13,
  author  = {Leinster, T. and Willerton, S.},
  title   = {On the asymptotic magnitude of subsets of {E}uclidean space},
  journal = {Geom. Dedicata},
  volume  = {164},
  pages   = {287--310},
  year    = {2013},
  doi     = {10.1007/s10711-012-9773-6}
}

@inproceedings{LASR24,
  author    = {Limbeck, K. and Andreeva, R. and Sarkar, R. and Rieck, B.},
  title     = {Metric space magnitude for evaluating the diversity of latent representations},
  booktitle = {Advances in Neural Information Processing Systems 37 ({NeurIPS} 2024)},
  year      = {2024},
  doi       = {10.52202/079017-3937}
}

@book {Mun84,
    AUTHOR = {Munkres, James R.},
     TITLE = {Elements of algebraic topology},
 PUBLISHER = {Addison-Wesley Publishing Company, Menlo Park, CA},
      YEAR = {1984},
     PAGES = {ix+454},
      ISBN = {0-201-04586-9},
   MRCLASS = {55-01},
  MRNUMBER = {755006},
MRREVIEWER = {Christopher\ W.\ Stark},
}

@article{NR09,
  author  = {Nagel, U. and Reiner, V.},
  title   = {Betti numbers of monomial ideals and shifted skew shapes},
  journal = {Electron. J. Combin.},
  volume  = {16},
  number  = {2},
  pages   = {Research Paper 3, 59},
  year    = {2009},
  doi     = {10.37236/69}
}

@article{Orl04,
  author  = {Orlov, D. O.},
  title   = {Triangulated categories of singularities and {D}-branes in {L}andau--{G}inzburg models},
  journal = {Proc. Steklov Inst. Math.},
  volume  = {246},
  pages   = {227--248},
  year    = {2004}
}

@book {Rie14,
    AUTHOR = {Riehl, Emily},
     TITLE = {Categorical homotopy theory},
    SERIES = {New Mathematical Monographs},
    VOLUME = {24},
 PUBLISHER = {Cambridge University Press, Cambridge},
      YEAR = {2014},
     PAGES = {xviii+352},
      ISBN = {978-1-107-04845-4},
   MRCLASS = {18G55 (18D20 55U35)},
  MRNUMBER = {3221774},
MRREVIEWER = {David\ A.\ Blanc},
       DOI = {10.1017/CBO9781107261457},
       URL = {https://doi.org/10.1017/CBO9781107261457},
}

@article{Sko06,
  author  = {Sk{\"o}ldberg, E.},
  title   = {{Morse} theory from an algebraic viewpoint},
  journal = {Trans. Amer. Math. Soc.},
  volume  = {358},
  number  = {1},
  pages   = {115--129},
  year    = {2006},
  doi     = {10.1090/S0002-9947-05-04079-1}
}

@article{SP94,
  author  = {Solow, A. and Polasky, S.},
  title   = {Measuring biological diversity},
  journal = {Environ. Ecol. Stat.},
  volume  = {1},
  pages   = {95--107},
  year    = {1994}
}

@book{Wei94,
  author    = {Weibel, C. A.},
  title     = {An introduction to homological algebra},
  series    = {Camb. Stud. Adv. Math.},
  volume    = {38},
  publisher = {Cambridge Univ. Press},
  address   = {Cambridge},
  year      = {1994},
  doi       = {10.1017/CBO9781139644136}
}
